\RequirePackage{fix-cm} 
\documentclass[a4paper,twoside,12pt,reqno]{amsart}

\usepackage{fixltx2e}     

\usepackage[english]{babel}
\usepackage[latin1]{inputenc}

\usepackage{indentfirst,verbatim}

\usepackage{amsmath,amsfonts,amssymb,amsgen,amsbsy,eucal,mathrsfs,dsfont}
\usepackage{stmaryrd}

\usepackage{a4wide,verbatim}

\usepackage{epsfig,rotating,color}
\usepackage{multicol}
\usepackage{tikz,pgfplots}
\usepackage{graphicx}
\usepackage{mathtools}
\usepackage{textcomp}

\usepackage{caption}

\usepackage{hyperref}

\usepackage[square,numbers]{natbib}
\usepackage{systeme}

\hypersetup{
colorlinks=true, 
breaklinks=true, 
urlcolor= blue, 
linkcolor= blue, 
bookmarksopen=true, 
pdftitle={}, 
pdfauthor={}, 
pdfsubject={} 
}

\newcommand{\ud}{\mathrm{d}}

\newcommand{\half}{{\textstyle{1\over2}}}
\newcommand{\third}{{\textstyle{1\over3}}}

\newcommand{\sixth}{{\textstyle{1\over6}}}

\usepackage{amsthm}
\newtheorem{thm}{Theorem}[section]
\newtheorem{definition}[thm]{Definition}
\newtheorem{lem}[thm]{Lemma}
\newtheorem{pro}[thm]{Proposition}

\theoremstyle{remark}
\newtheorem{rem}[thm]{Remark}

\usepackage{enumitem}
\newlist{steps}{enumerate}{1}
\setlist[steps, 1]{label = Step \arabic*:}

\newcommand{\eqdef}{\stackrel{\text{\tiny{def}}}{=}}

\title[Weak solutions of Serre--Green--Naghdi equations]{\bf Global weak solutions of the Serre--Green--Naghdi equations with surface tension}

\author[GUELMAME]{Billel Guelmame}

\newcommand{\nfont}{\fontshape{n}\selectfont}

\address{({\nfont\textbf{Billel Guelmame}})  LJAD,  Inria, CNRS,  Universit\'e C\^ote d'Azur, France.} 
\email{billel.guelmame@univ-cotedazur.fr}

\newcommand\smallO{
  \mathchoice
    {{\scriptstyle\mathcal{O}}}
    {{\scriptstyle\mathcal{O}}}
    {{\scriptscriptstyle\mathcal{O}}}
    {\scalebox{.7}{$\scriptscriptstyle\mathcal{O}$}}
  }

\let\oldtocsection=\tocsection
 
\let\oldtocsubsection=\tocsubsection

\renewcommand{\tocsection}[2]{\hspace{0em}\oldtocsection {#1}{#2}}
\renewcommand{\tocsubsection}[2]{\hspace{2em}\oldtocsubsection{#1}{#2}}

\numberwithin{equation}{section}

\allowdisplaybreaks

\usepackage{lineno}

\begin{document}

\maketitle

\begin{abstract} 
We consider in this paper the Serre--Green--Naghdi equations with surface tension. 
Smooth solutions of this system conserve an $H^1$-equivalent energy.
We prove the existence of global weak dissipative solutions for any relatively small-energy initial data.
We also prove that the Riemann invariants of the solutions satisfy a one-sided Oleinik inequality.
\end{abstract}

\medskip


\medskip


\tableofcontents

\section{Introduction}

The Euler equations are usually used to describe water waves in oceans and Channels.
Due to the difficulties to resolve the Euler equations both numerically and analytically, several simpler approximations have been proposed in the literature for different regimes.
In the shallow water regime, the main assumption is on the ratio of the mean water depth $\bar{h}$ to the wave wave-length $\iota$, the shallowness parameter $\sigma = \bar{h}^2/\iota^2$ is considered to be small. 
Beside the shallowness condition, a restriction on the amplitude of the wave $a$ can be considered assuming that the nonlinearity (or the amplitude) parameter $\epsilon=a/\bar{h}$ is small.
Considering the shallow water regime with the small-amplitude condition \cite{Johnson,Lannes} ($\sigma \ll	 1$, $\epsilon \ll	 1$). Many equations have been derived to model the propagation of the waves, such as the Camassa--Holm equation \cite{CamassaHolm1993}, the Korteweg--deVries (KdV) equation \cite{KdV} and some variants of the Boussinesq equations \cite{Bona,Boussinesq,Whitham}.
Considering shallow water with possibly large-amplitude waves ($\sigma \ll 1$, $\epsilon \approx 1$), by neglecting the terms of order $\mathcal{O}(\sigma)$ in the water waves equations, Saint-Venant obtained the Nonlinear Shallow Water (or Saint-Venant) equations \cite{Wehausen}. Smooth solutions of the Saint-Venant equations have a precision of order $\mathcal{O}(t \sigma)$, where $t$ denotes the time \cite{Lannes}. 
In order to obtain a better precision, one can keep the $\mathcal{O}(\sigma)$ terms in the equations and only neglect the $\mathcal{O}(\sigma^2)$ terms. This leads to the Serre--Green--Naghdi equations. Those equations were firstly derived by Serre \cite{Serre},  rediscovered independently by Su and Gardner \cite{Su} and another time by Green, Laws and Naghdi \cite{Green1,Green2}.
The Serre--Green--Naghdi equations are the most general and most precise, but also the most complicated of the models of shallow water equations presented above.
Of course, one can always keep higher order terms in the equation (keeping terms of order $\mathcal{O}(\sigma^2)$ for example), this will lead to equations with a better precision, but with higher order derivatives. Those equations are not accurate due to the high order derivative terms which make their numerical resolution much slower.

The influence of the surface tension is generally neglected on water waves problems. However, in certain cases, the effect of the surface tension is appreciable. 
Indeed, Longuet-Higgins \cite{Longuet} showed that the surface tension is significant in certain localized regions, and cannot be neglected near the sharp crest of the breaking wave.
Other experimental studies showed the importance of the surface tension on thin layers \cite{Falcon,Myers,Packham}. Those experimentations have been done for different fluids including water and mercury.
Various mathematical studies of water waves equations with surface tension exist in the literature, we refer to \cite{Alazard,Ambrose1,Ambrose2,Beyer,
CDD,CDG,Ming,Schweizer,Shatah,Yosihara}.
\vspace{1cm}
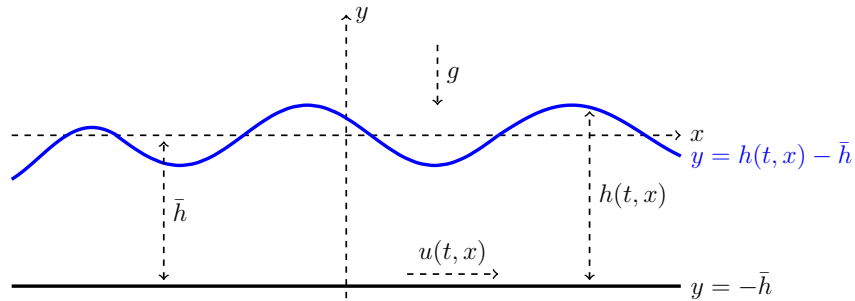
\begin{figure}[!ht]
\begin{tikzpicture}[thick, transform canvas={scale=0.8}, shift={(2.3,-5)}]

\draw[dashed, ->] (-8,4) -- (3,4) node[right]{$x$};

\draw[dashed, ->] (-2.5,1.3) -- (-2.5,6)  node[right]{$y$};

\draw[scale=1, ultra thick, domain=-8:-6.25,smooth,variable=\x,blue] plot ({\x},{0.5*cos(deg(1.8*(\x-0.3)))+3.63});
 \draw[scale=1, ultra thick, domain=-6.29:0,smooth,variable=\x,blue] plot ({\x},{0.5*sin(deg(1.5*\x))+4});
 \draw[scale=1, ultra thick, domain=0:3,smooth,variable=\x,blue] plot ({\x},{0.5*sin(deg(1.3*\x))+4}) node[right, blue]{$y=h(t,x) - \bar{h}$};

\draw[ultra thick] (-8,1.5) -- (3,1.5) node[right]{$y=- \bar{h}$};

\draw[dashed, <->] (-5.5,1.6) -- (-5.5,3.9) node[midway, right]{$\bar{h}$};

\draw[dashed, <->] (1.5,1.6) -- (1.5,4.4) node[midway, right]{$h(t,x)$};

\draw[dashed, ->] (-1.5,1.7) -- (0,1.7) node[midway, above]{$u(t,x)$};
\draw[dashed, ->] (-1,5.5) -- (-1,4.5) node[midway, right]{$g$};

\end{tikzpicture}
\vspace{3cm}
\caption{Fluid domain.}
\label{Sketch}
\end{figure}

Considering a two-dimensional coordinate system $Oxy$ (Figure \ref{Sketch}) and an incompressible fluid layer.
Considering the still fluid level at $y=0$,  the fluid layer is bounded between the flat bottom at $y=-\bar{h}$ and a free surface $y=h(t,x)-\bar{h}$, where $h$ is the total water depth.
Assuming long waves in shallow water with possibly large-amplitude.
The Serre--Green--Naghdi system (without neglecting the surface tension influence) reads
\begin{subequations}\label{SGN0}
\begin{gather}
h_t\ +\,\left[\, h\,u\,\right]_x\ =\ 0, \label{SGNa} \\
\left[\,h\,u\,\right]_t\ +\,\left[\,h\,u^2\,+\, \half\, g\, h^2\, \label{SGN1b}
+\,\mathscr{R}\,\right]_x\ =\ 0,\\ \label{SGNc}
\mathscr{R}\, \eqdef\, \third\, h^3 \left(-u_{tx}\, -\, u\, u_{xx}\, +\, u_x^2 \right)\, 
-\, \gamma \left( h\, h_{xx}\ -\ \half\, h_x^2 \right),
\end{gather}
\end{subequations}
where $u$ denotes the depth-averaged horizontal velocity, $g$ is
the gravitational acceleration and $\gamma>0$ is a constant (the ratio of the surface tension coefficient to the density).
The classical Serre--Green--Naghdi equations (without surface tension) are recovered taking $\gamma=0$.
The Serre--Green--Naghdi (SGN$_{\gamma}$) equations \eqref{SGN0} have been derived in \cite{Dias} as a generalisation of the classical SGN equations ($\gamma=0$).
Due to the appearance of time derivatives in \eqref{SGNc}, it is convenient to apply the inverse of the Sturm--Liouville operator
\begin{equation}\label{Ldef}
\mathcal{L}_h\ \eqdef\ h\ -\ \third\, \partial_x\, h^3\, \partial_x,
\end{equation} 
the system \eqref{SGN0} becomes then
\begin{subequations}\label{SGN1}
\begin{align}
h_t\ +\,\left[\, h\,u\,\right]_x\ &=\ 0, \\ \label{SGN1_2}
u_t\ +\ u\, u_x\ + g\, h_x\ &=\ - \mathcal{L}_{h}^{-1} \partial_x\, \left\{ {\textstyle \frac{2}{3}\, h^3\, u_x^2\ } -\ \left[\gamma\, h\ -\ \third\, g\, h^3 \right]\, h_{xx}\ +\ \half\, \gamma\, h_x^2 \right\}.
\end{align}
\end{subequations}
When $h>0$, the operator $\mathcal{L}_h^{-1}$ is well defined and smoothes two derivatives (see Lemma \ref{Inverseesitimates} below).
This is not enough to control the term containing $h_{xx}$ on the right-hand side of \eqref{SGN1_2}.
To overcome this problem, we use the definition of $\mathcal{L}_h$ to rewrite the system \eqref{SGN1} in the equivalent form 
\begin{subequations}\label{SGN}
\begin{align} 
h_t\ +\,\left[\, h\,u\,\right]_x\ &=\ 0, \\ 
u_t\ +\ u\, u_x\ + 3\, \gamma h^{-2} h_x\ &=\ - \mathcal{L}_{h}^{-1}  \partial_x\, \left\{ {\textstyle \frac{2}{3}\, h^3\, u_x^2\ }  -\ {\textstyle \frac{3}{2}}\, \gamma\, h_x^2\ +\ \half\, g\, h^2  -3\, \gamma \ln(h) \right\}.
\end{align}
\end{subequations}
Smooth solutions of the SGN$_{\gamma}$ equations \eqref{SGN} satisfy the energy equation (see Appendix \ref{App:B})
\begin{equation}\label{ene}
\mathscr{E}_t\ +\ \mathscr{D}_x\ =\ 0,
\end{equation}
where 
\begin{align}\label{TildeEneDef}
\mathscr{E}\ &\eqdef\ \half h\, u^2\ +\ \half\, g\, (h-\bar{h})^2\ +\ \sixth\, h^3\, u_x^2\ +\ \half\, \gamma\, h_x^2, \\
\mathscr{D}\ &\eqdef\ u\, \mathscr{E}\ +\ u\, \left( \mathscr{R}\ +\ \half\, g\, h^2\ -\ \half\, g\, \bar{h}^2 \right)\ +\ \gamma\, h\, h_x\, u_x.
\end{align}
Linearising the SGN$_{\gamma}$ equations \eqref{SGN} around the constant state $(h,u)=(\bar{h},0)$ and looking for travelling waves having the form $\exp \left\{ \left( k x - \omega t \right) i \right\}$ we obtain the dispersion relation $\omega^2 = g \bar{h} k^2  \left(1 + \gamma k^2/g \right)/ \left(1 + \bar{h}^2 k^2/3 \right)$. Defining the Bond number $B = g \bar{h}^2/\gamma$, the SGN$_{\gamma}$ equations are linearly dispersive if and only if $B\neq 3$.
In the dispersionless case ($B=3$), the SGN$_{\gamma}$ equations admit weakly singular peakon travelling wave solutions \cite{Dutykh2018,Mitsotakis}.
A mathematical study of the Serre--Green--Naghdi equations with or without surface tension have been widely studied in the literature. We refer to \cite{Alvarez,Israwi2011,Kazerani,Khorbatly,
Lannes,Li} for the case $\inf h_0 >0$ and to \cite{Lannes18} for the shoreline problem ($\mathrm{sign}(h) = \mathds{1}_{x > x_0}$).
In \cite{Alvarez,Israwi2011,Li}, a proof of the local well-posedness of the SGN equations without surface tension ($\gamma=0$) is given. Kazerani has proved in \cite{Kazerani} the existence of global smooth solutions of the SGN equations with viscosity for small initial data.
A full justification of the model \eqref{SGN} is given in \cite{Khorbatly,Lannes}.
By ``full justification'' we mean local well-posedness of the system and, the solution is close to the solution of the water waves equations with the same initial data.
In a recent work \cite{SGNbu2021}, we have obtained a precise blow-up criterion of \eqref{SGN} (Theorem \ref{thm:buc} below) and we proved that such scenario occurs for a class of small-energy initial data (Theorem \ref{thm:bu} below).
Then, in general, smooth solutions cannot exist globally in time.

This paper investigates the existence of global weak solutions of \eqref{SGN} with $\gamma>0$. To the best of the author's knowledge, the existence of global weak solutions for all the different variants of the inviscid Serre--Green--Naghdi equations has not been established before.
Here, the existence of global weak solutions is established by approximating the system \eqref{SGN} with another system that admits global smooth solutions. We recover weak solutions of \eqref{SGN} by taking the limit. The proof involves several steps. 

We consider initial data satisfying $\int \mathscr{E}_0\, \mathrm{d}x < \sqrt{g \gamma} \bar{h}^2$, which is propagated due to energy conservation \eqref{ene}. Using the fact that the energy is equivalent to $\|(h-\bar{h},u)\|^2_{H^1}$ and a Sobolev-like inequality (essentially $H^1 \hookrightarrow  L^\infty$, see Proposition \ref{pro:ene} below) we obtain a uniform lower bound of $h$.
This is important for ensuring the invertibility of the operator $\mathcal{L}_h$ defined in \eqref{Ldef}.
Smooth solutions of \eqref{SGN} blow-up in finite time due to the presence of quadratic terms in the associated Riccati-type equations.
In order to approximate the SGN$_{\gamma}$ system, we use a cut-off to obtain a linear growth that leads to global smooth solutions (due to Gronwall's inequality).
However, cutting-off directly as in \cite{wave2,wave3} violates the energy conservation \eqref{ene}.
The choice the approximated system is crucial and must conserve the properties of the SGN$_{\gamma}$ system. 
In Section \ref{sec:App_sys} below, we chose carefuly a suitale approximated system that is globaly well-posed and satisfies the energy equation \eqref{energyequationep}.
In order to pass to the limit, some uniform estimates are needed. 
In the previous studies of smooth solutions of the SGN equations, some estimates of the operator $\mathcal{L}^{-1}_h$ have been obtained, those estimates usually depend on the $L^\infty$ norm of $h_x$ which may blow-up for weak solutions. 
In Lemma \ref{Inverseesitimates} below, we present some new estimates of $\mathcal{L}^{-1}_h$ depending only on the $L^\infty$ norm of $h$ and $1/h$.
As in \cite{CH_CPAM,wave2,wave3}, an $L^p_{loc}$ estimate of $(h_x,u_x)$ with $p<3$ is also needed.
In our case and due to the complexity of the SGN$_{\gamma}$ equations, we have to use a change of coordinates to obtain this estimate (see Lemma \ref{lem:alpha+2} below).
We use then some classical compactness arguments 
with Young mesures \cite{Focusing} and a generalised compensated compactness result \cite{PG1991} to pass to the limit.
We follow in this step the techniques developed in \cite{CH_CPAM} for the Camassa-Holm equation and in \cite{wave2,wave3} for the variational wave equation. 
The structure of the SGN$_{\gamma}$ system being more complex, we have to handle the weak limit of some nonlinear terms that do not exist in \cite{wave2,wave3} (see Lemma \ref{BTlim} for example).
Finally, the global weak solutions of \eqref{SGN} are obtained by taking the limit in the approximated system, and are shown to dissipate the energy and satisfy the one-sided Oleinik inequality \eqref{Ol_xi}.

The existence of global solutions to the Boussinesq equations \cite{Boussinesq,Whitham} 
\begin{equation}\label{Bouss}
h_t\ +\ [h\, u]_x\ =\ 0, \qquad \qquad u_t\ +\ u\, u_x\ +\ g\, h_x\ =\ u_{txx}
\end{equation}
have been studied in \cite{Amick,Schonbek}. Schonbek \cite{Schonbek} regularised the conservation of the mass by adding a defusion term, i.e., $h_t + [h u]_x = \varepsilon h_{xx}$, with $\varepsilon >0$. She proved the global well-posedness of the regularised system, and she obtained global weak solutions of \eqref{Bouss} by taking $\varepsilon \to 0$.
In \cite{Amick}, Amick proved that if the initial data, $(h_0,u_0)$, is smooth, then the solution, $(h,u)$, obtained by Schonbek \cite{Schonbek} is also smooth and is the unique smooth solution of the Boussinesq equations \eqref{Bouss}.

The SGN$_{\gamma}$ equations \eqref{SGN0} can be compared with the dispersionless regularised Saint-Venant (rSV) system presented in \cite{ClamondDutykh2018a}. The rSV system can be obtained replacing $\mathscr{R}$ in \eqref{SGNc} by $\varepsilon \mathscr{R}_{rSV}$ with
\begin{equation*}
\mathscr{R}_{rSV}\ \eqdef\ h^3\left(\,u_x^{\,2}\,-\,u_{xt}\, -\, u\, u_{xx}\, \right)\, -\ 
g\, h^2 \left(\,h\,h_{xx}\, +\, \half\, h_x^{\,2}\,\right)
\end{equation*}
and $\varepsilon \geqslant 0$, the classical Saint-Venant system is recovered taking $\varepsilon =0$.
Weakly singular shock profiles of the rSV equations are studied in \cite{PuEtAl2018}. 
In \cite{liu2019well}, Liu et al. proved the local well-posedness of the rSV equations and identified a class of initial data such that the corresponding solutions blow-up in finite time.
The rSV system have been generalised recently to obtain a regularisation of any unidimensional barotropic Euler (rE) system \cite{guelmame2020Euler}.
The system \eqref{SGN0} can also be compared with the modified Serre--Green--Naghdi (mSGN) equations derived in \cite{ClamondEtAl2017} to improve the dispersion relation of the classical SGN system. The mSGN system presented in \cite{ClamondEtAl2017} can be obtained replacing $\mathscr{R}$ in \eqref{SGNc} by
\begin{equation*}
\mathscr{R}_{mSGN}\, \eqdef\, \third\, \left( 1\, +\, {\textstyle \frac{3}{2}}\, \beta \right)\, h^3 \left(-u_{tx}\, -\, u\, u_{xx}\, +\, u_x^2 \right)\, 
-\, \half\, \beta\, g\, h^2\, \left( h\, h_{xx}\ +\ \half\, h_x^2 \right).
\end{equation*}
where $\beta$ is a positive parameter.
The rSV, rE and mSGN conserve $H^1$-equivalent energies and have similar properties as the SGN$_{\gamma}$ system \eqref{SGN0}.
One may can obtain the existence of global weak solutions of those equations following the proof given in this paper.

The study of the classical Serre--Green--Naghdi equations is more challenging. 
Indeed, when $\gamma = 0$, the energy \eqref{TildeEneDef} fails to control the $H^1$ norm of $h-\bar{h}$, then, a lower bound of $h$ cannot be obtained. 
This bound is crucial to obtaining the blow-up result \cite{SGNbu2021} and the global existence in this paper for $\gamma>0$.
To the author's knowledge, the questions of the blow-up of smooth solutions and the existence of global solutions of the SGN equations without surface tension are still open.
However, Bae and Granero-Belinch\'en \cite{Bae} proved recently that for a class of periodic initial data satisfying $\inf h_0=0$, smooth solutions cannot exist globally in time. 
For this class of initial data, it is not known if smooth solutions exist locally in time, but if they do, a singularity must appear in finite time.

This paper is organised as follows. 
In Section \ref{sec:pr} we present the local well-posedness of \eqref{SGN} and some blow-up results. Section \ref{sec:mr} is devoted to define weak solutions of \eqref{SGN} and to present the main result which is the existence of global dissipative weak solutions. 
We discuss in Section \ref{sec:App_sys} the properties needed of the approximated system and we propose a suitable choice. 
Section \ref{sec:UE} is devoted to prove the existence of global smooth solutions of the approximated system and to obtain some uniform estimates.
We obtain strong precompactness results in Section \ref{Sec:Precom}.
The existence of the global weak solutions is proved in Section \ref{sec:gws}.
In Appendix \ref{App:A} we recall some classical lemmas that are used in this paper.
Appendix \ref{App:B} is devoted to obtain the energy equations of the approximated system and of \eqref{SGN}.

\section{Review of previous results}\label{sec:pr}

We consider the Serre--Green--Naghdi equations with surface tension in this form
\begin{subequations}\label{SGNxi}
\begin{gather} \label{SGN2a}
h_t\ +\,\left[\, h\,u\,\right]_x\ =\ 0, \\ \label{SGN2b}
u_t\ +\ u\, u_x\ + 3\, \gamma h^{-2} h_x\ =\ - \mathcal{L}_{h}^{-1}  \partial_x\, \left\{ \mathscr{C}\, +\, F(h) \right\},\\ \label{initialx}
u(0,x)\ =\ u_0(x), \qquad h(0,x)\ =\ h_0(x),
\end{gather}
\end{subequations}
with
\begin{align}\label{Cdef}
\mathscr{C}\ &\eqdef\  {\textstyle \frac{2}{3}}\, h^3\, u_x^2\   -\ {\textstyle \frac{3}{2}}\, \gamma\, h_x^2,\\ \label{Fdef}
F(h)\ &\eqdef\ \half\, g\, h^2 -\ \half\, g\, \bar{h}^2\  -\ 3\, \gamma \ln(h/\bar{h}).
\end{align}
The system \eqref{SGNxi} is locally well-posed in the Sobolev space 
\begin{equation*}
H^s(\mathds{R})\ \eqdef\ \left\{f,\, \|f\|_{H^s(\mathds{R})}^2\ \eqdef\ \int_{\mathds{R}} \left( 1\, +\, |\xi|^2 \right)^s\, |\hat{f}(\xi)|^2\, \mathrm{d}\xi\ <\ \infty \right\}
\end{equation*}
where $s \geqslant 2$ is a real number.
\begin{thm}\label{thm:ex}
Let $\gamma>0$, $\bar{h}>0$ and $s \geqslant 2$, then, for any $(h_0-\bar{h},u_0) \in H^s(\mathds{R})$ satisfying $\inf_{x\, \in\, \mathds{R}}h_0(x)>0$ there exists $T>0$ and $(h-\bar{h},u) \in C([0,T],H^s(\mathds{R})) \cap C^1([0,T],H^{s-1}(\mathds{R}))$ a unique solution of \eqref{SGNxi} such that 
\begin{equation*}
\inf_{(t,x)\, \in\, [0,T] \times \mathds{R}}\ h(t,x)\ >\ 0.
\end{equation*}
Moreover, the solution satisfies the conservation of the energy
\begin{equation}\label{eneC}
\frac{\mathrm{d}}{\mathrm{d}t}\, \int_\mathds{R} \left( \half h\, u^2\ +\ \half\, g\, (h-\bar{h})^2\ +\ \sixth\, h^3\, u_x^2\ +\ \half\, \gamma\, h_x^2 \right)\, \mathrm{d}x\ =\ 0.
\end{equation}
\end{thm}
\begin{rem}
The solution given in Theorem \ref{thm:ex} depends continuously on the initial data, i.e., If $(h_0^1-\bar{h},u_0^1), (h_0^2-\bar{h},u_0^2) \in H^s$, such that $h_0^1,h_0^2 \geqslant h_{min} >0 $, then for all $t \leqslant T$ there exists a constant $C(\|(h^2-\bar{h},u^2)\|_{L^\infty([0,t],H^{s})},\|(h^1-\bar{h},u^1)\|_{L^\infty([0,t],H^{s})})>0$, such that 
\begin{equation*}
\left\|\left(h^1-h^2,u^1-u^2\right)\right\|_{L^\infty([0,t],H^{s-1})}\ \leqslant\ C\, \left\|\left(h_0^1-h_0^2,u_0^1-u_0^2\right)\right\|_{H^s}.
\end{equation*}
\end{rem}
The proof of Theorem \ref{thm:ex} is classic and omitted in this paper. See \cite{guelmame2020Euler,Israwi2011,Lannes,liu2019well} and Theorem 3 of \cite{guelmame2020rSV} for more details.
It is clear from Theorem \ref{thm:ex} that if the solution at time $T$ remains in $H^s$ and $\inf_x h(T,x)>0$ then one can extend the interval of existence. This leads to the blow-up criterion
\begin{equation*}
T_{max}\, <\, \infty \ \implies \ \liminf_{t\to T_{max}}\, \inf_{x\, \in\, \mathds{R}}\, h(t,x)\ =\ 0\quad \mathrm{or}\quad \limsup_{t\to T_{max}}\, \|(h-\bar{h},u)\|_{H^s}\ =\ \infty,
\end{equation*}
where $T_{max}$ is the maximum time existence of the solution.
This criterion has been improved in \cite{SGNbu2021} to 
\begin{thm}(\cite{SGNbu2021})\label{thm:buc}
Let $T_{max}$ be the maximum time existence of the solution given by Theorem \ref{thm:ex}, then
\begin{equation*}
T_{max}\, <\, \infty \ \implies \ \liminf_{t\to T_{max}}\, \inf_{x\, \in\, \mathds{R}}\, h(t,x)\ =\ 0\quad \mathrm{or}\quad
\left\{\begin{aligned} \liminf_{t\to T_{max}}\, \inf_{x\, \in\, \mathds{R}}\, &u_x(t,x)\ =\ -\infty, \\ \mathrm{and}& \\ \limsup_{t\to T_{max}}\, \|h_x&(t,x)\|_{L^\infty}\, =\, \infty, \end{aligned}\sysdelim. . \right.
\end{equation*}
which is equivalent to the second criterion
\begin{equation*}
T_{max}\, <\, \infty \ \implies \ 
\limsup_{t\to T_{max}}\, \|u_x(t,x)\|_{L^\infty}\, =\, \infty\ \mathrm{and}\  
\left\{\begin{aligned} \liminf_{t\to T_{max}}\, \inf_{x\, \in\, \mathds{R}}\, &h(t,x)\ =\ 0, \\ \mathrm{or}& \\ \limsup_{t\to T_{max}}\, \|h_x(t,&x)\|_{L^\infty}\, =\, \infty. \end{aligned}\sysdelim. . \right.
\end{equation*}
\end{thm}
Noting that the energy conserved in \eqref{eneC} is equivalent to the $H^1$ norm of $(h-\bar{h},u)$. Due to the continuous embedding $H^1 \hookrightarrow L^\infty$, we can obtain a uniform (on time) estimate of $\|(h-\bar{h},u)\|_{L^\infty}$, and, if the initial energy is not very large compared to $\bar{h}$, we can obtain a lower bound of $h$.
For that purpose, we present the following proposition.
\begin{pro}\label{pro:ene}
For $\gamma > 0$, $\bar{h}>0$, let $E$ be a positive number such that
\begin{align}\label{Edef}
0\ <\ E\ <\ \sqrt{g\, \gamma}\, \bar{h}^2, 
\end{align}
Defining 
\begin{gather*}
h_{min}\ \eqdef\ \bar{h}\, -\, (g \gamma)^{-\frac{1}{4}} \sqrt{E}, \qquad 
h_{max}\ \eqdef\ \bar{h}\, +\, (g \gamma)^{-\frac{1}{4}} \sqrt{E}, \\ 
u_{max}\ \eqdef\ -u_{min}\ \eqdef\ 3^\frac{1}{4} \sqrt{E}/h_{min} .
\end{gather*}
Then, for any $(h-\bar{h},u) \in H^1$ satisfying $\int \mathscr{E}\, \mathrm{d}x \leqslant E$, we have
\begin{gather}\label{ubound}
0\ <\ h_{min}\ \leqslant\ h\ \leqslant h_{max}\ < 2\, \bar{h}, \qquad u_{min}\ \leqslant\ u\ \leqslant\ u_{max}, 
\end{gather}
\end{pro}
\begin{rem} Taking an initial data satisfying $\int_\mathds{R} \mathscr{E}_0\, \mathrm{d}x \leqslant E$, then, due to the energy conservation \eqref{eneC} and Proposition \ref{pro:ene} the depth $h$ cannot vanish. The blow-up criteria given in Theorem \ref{thm:buc} becomes then
\begin{equation*}
T_{max}\, <\, \infty \ \implies \ 
\inf_{[0,T_{max}) \times \mathds{R}}\, u_x(t,x)\ =\ -\infty \qquad  \mathrm{and} \qquad 
\limsup_{t\to T_{max}}\, \|h_x(t,x)\|_{L^\infty}\ =\ \infty. 
\end{equation*}
\end{rem}
\textit{Proof of Proposition \ref{pro:ene}.}
The Young inequality $\frac{1}{2} a^2 + \frac{1}{2} b^2 \geqslant \pm ab$ implies that
\begin{align*} 
E\ \geqslant\ \int_\mathds{R} \mathscr{E}\, \mathrm{d}y\ 
\geqslant &\ 
\int_\mathds{R} \left( \half\, g\, (h-\bar{h})^2\ +\ \half\, \gamma\, h_x^2 \right)\, \mathrm{d}x\\
\geqslant &\  \sqrt{g\, \gamma}\, \left( \int_{-\infty}^x (h\, -\, \bar{h})\, h_x\, \mathrm{d}y\ -\ \int^{\infty}_x (h\, -\, \bar{h})\, h_x\, \mathrm{d}y  \right)\\
\geqslant &\ \sqrt{g\, \gamma}\, |h\, -\, \bar{h}|^2,
\end{align*}
which implies that $h_{min} \leqslant h \leqslant h_{max}$. Doing the same estimates with $u$ one obtains
\begin{align*} 
E\ \geqslant\ \int_\mathds{R} \mathscr{E}\, \mathrm{d}y\ 
\geqslant &\  \int_\mathds{R} \left( \half h\, u^2\ +\ \sixth\, h^3\, u_x^2 \right)\, \mathrm{d}y\\
\geqslant &\ {\textstyle \frac{1}{\sqrt{3}}}\, h_{min}^{2}\, \left(  \int_{-\infty}^x u\, u_x\, \mathrm{d}y\ -\ \int^{\infty}_x u\, u_x\, \mathrm{d}y \right)\\
\geqslant &\  {\textstyle \frac{1}{\sqrt{3}}}\, h_{min}^{2}\, |u|^2,
\end{align*}
the last inequality ends the proof of $u_{min} \leqslant u \leqslant u_{max}$. \qed

As in \cite{SGNbu2021}, we can build some initial data with small initial data such that the corresponding solutions blow-up in small time.
\begin{thm}(\cite{SGNbu2021})\label{thm:bu}
For any $T>0$ and $E$ satisfying \eqref{Edef},  there exist
\begin{itemize}
\item  $(h_0-\bar{h},u_0) \in C^\infty_c(\mathds{R})$ satisfying $\int_\mathds{R} \mathscr{E}_0\, \mathrm{d}x \leqslant E$ such that the corresponding solution of \eqref{SGNxi} blows-up at finite time $T_{max} \leqslant T$ and 
\begin{equation*}
\inf_{[0,T_{max}) \times \mathds{R}}  u_x(t,x)\ =\ -\infty,
\quad
\sup_{[0,T_{max}) \times \mathds{R}}  h_x(t,x)\ =\ \infty,
\quad
\inf_{[0,T_{max}) \times \mathds{R}}  h_x(t,x)\ >\ -\infty.
\end{equation*}
\item  $(\tilde{h}_0-\bar{h},\tilde{u}_0) \in C^\infty_c(\mathds{R})$ satisfying $\int_\mathds{R} \tilde{\mathscr{E}}_0\, \mathrm{d}x \leqslant E$ such that the corresponding solution of \eqref{SGNxi} blows-up at finite time $\tilde{T}_{max} \leqslant T$ and 
\begin{equation*}
\inf_{[0,\tilde{T}_{max}) \times \mathds{R}}  \tilde{u}_x(t,x)\ =\ -\infty,
\quad
\inf_{[0,\tilde{T}_{max}) \times \mathds{R}}  \tilde{h}_x(t,x)\ =\ -\infty,
\quad
\sup_{[0,\tilde{T}_{max}) \times \mathds{R}}  \tilde{h}_x(t,x)\ <\ \infty.
\end{equation*}
\end{itemize}
\end{thm}

\section{Main results}\label{sec:mr}

Since smooth solutions fail to exist globally in time, even for arbitrary small-energy initial data, we shall define weak solutions of the  SGN$_\gamma$ system \eqref{SGNxi}.
For that purpose, we define the domain $\mathfrak{D} \subset H^1 $
\begin{equation}
\mathfrak{D}\ \eqdef\ \left\{ \left(h-\bar{h},u \right) \in H^1,\ \int_\mathds{R} \mathscr{E}\, \mathrm{d}x < \sqrt{g \gamma} \bar{h}^2 \right\}.
\end{equation}
\begin{definition}\label{WSDef} We say that $(h-\bar{h},u) \in L^\infty(\mathds{R}^+, H^1) \cap \mathrm{Lip}(\mathds{R}^+, L^2)$ is a weak solution of \eqref{SGNxi} if it satisfies the initial condition \eqref{initialx} with \eqref{SGN2a} in $L^2$ and for all $\varphi \in C^\infty_c((0,\infty) \times \mathds{R})$ we have
\begin{equation}\label{WS_xi}
\int_{\mathds{R}^+ \times \mathds{R}} \left\{ \left\{u_t\, +\, u\, u_x\, +\,  3\, \gamma\, (h)^{-2}\, h_x \right\} \mathcal{L}_h\, \varphi\, -\, \varphi_x \left\{ \mathscr{C}\, +\, F(h) \right\} \right\} \mathrm{d}x\, \mathrm{d}t\
=\ 0 .
\end{equation}
Moreover, $(h(t,\cdot)-\bar{h},u(t,\cdot))$ belongs to $\mathfrak{D}$ for all $t \geqslant 0$ and $(h-\bar{h},u) \in C_r(\mathds{R}^+, H^1)$. More precisely, for all $t_0 \geqslant 0$ we have
\begin{equation}
\lim_{\substack{t\to t_0\\ t> t_0}} \left\| \left(h(t,\cdot)\, -\, h(t_0,\cdot),\ u(t,\cdot)\, -\, u(t_0,\cdot) \right) \right\|_{H^1}\ =\ 0.
\end{equation}
\end{definition}
Now we can state the main result of this paper.
\begin{thm}\label{mainthm}
Let $\bar{h},g,\gamma>0$ and $(h_0-\bar{h},u_0) \in \mathfrak{D}$, then there exist a global weak solution $(h-\bar{h},u) \in L^\infty ([0,\infty ), H^1(\mathds{R})) \cap C([0,\infty ) \times \mathds{R})$ of \eqref{SGNxi} in the sense of Definition \ref{WSDef}. Moreover,
\begin{itemize}
\item For any bounded set $\Omega = [t_1,t_2] \times [a,b] \subset (0,\infty) \times \mathds{R}$ and $\alpha \in [0,1)$ there exists $C_{\alpha, \Omega}>0$ such that 
\begin{equation}\label{alpha+2_xi}
\int_\Omega \left[ |h_t|^{2+\alpha}\ +\ |h_x|^{2+\alpha}\ +\ |u_t|^{2+\alpha}\ +\ |u_x|^{2+\alpha} \right] \mathrm{d}x\, \mathrm{d}t\ \leqslant\ C_{\alpha, \Omega} . 
\end{equation}
\item The solution dissipates the energy
\begin{equation}\label{Eneeq_xi}
\int_\mathds{R} \mathscr{E}\, \mathrm{d} x\ \leqslant\ \int_\mathds{R} \mathscr{E}_0\, \mathrm{d} x.
\end{equation}
\item There exists $C>0$ such that the solution satisfies the Oleinik inequality
\begin{equation}\label{Ol_xi}
u_x\, \pm\ \sqrt{3\, \gamma}\, h^{-\frac{3}{2}}\, h_x\ \leqslant\ C\, \left(1+{\textstyle \frac{1}{t}}\right), \quad a.e.\ (t,x) \in (0,\infty) \times \mathds{R}.
\end{equation}
\end{itemize}
\end{thm}
\begin{rem}
The constants $C_{\alpha, \Omega}$ and $C$ depend on $\bar{h},\gamma,g$ and $\int_\mathds{R} \mathscr{E}_0\, \mathrm{d}x$ but not on the initial data.
\end{rem}

In order to obtain global solutions of \eqref{SGNxi}, we use a suitable approximation of the system \eqref{SGNxi} that admits global smooth solutions. Using some compactness arguments and taking the limit we recover a global weak solution of \eqref{SGNxi}.
In the next section we present the choice of the suitable approximated system.

\section{An approximated system}\label{sec:App_sys}

The blow-up of the solutions given in Theorem \ref{thm:bu} is due to the Riccati-type equations. In order to prevent the singularities from appearing, we modify slightly the Riccati-type equations.

\subsection{Riccati-type equations}\label{sec_Riccati}
Defining the Riemann invariants \footnote{Those quantities are constants along the characteristics if the right-hand side of \eqref{SGNxi} is zero.} $R$ and $S$
\begin{gather}
R\ \eqdef\ u\ +\ 2\, \sqrt{3\, \gamma}\, h^{-\frac{1}{2}}, \qquad 
S\ \eqdef\ u\ -\ 2\, \sqrt{3\, \gamma}\, h^{-\frac{1}{2}}, \\
\lambda\ \eqdef\ u\ -\ \sqrt{3\, \gamma}\, h^{-\frac{1}{2}}, \qquad   
\eta\ \eqdef\ u\ +\ \sqrt{3\, \gamma}\, h^{-\frac{1}{2}}.
\end{gather}
The system \eqref{SGNxi} can be rewritten as 
\begin{subequations}\label{SGNR}
\begin{align}
R_t\ +\ \lambda\, R_x\ &=\ - \mathcal{L}_{h}^{-1}\,  \partial_x\, \left\{ \mathscr{C}\, +\, F(h) \right\}, \\
S_t\ +\ \eta\, S_x\ &=\ - \mathcal{L}_{h}^{-1}\,  \partial_x\, \left\{ \mathscr{C}\, +\, F(h) \right\}.
\end{align}
\end{subequations}
Defining
\begin{align*}
P\ \eqdef\ h\, R_x\ =\ h\, u_x\ -\ \sqrt{3\, \gamma}\, h^{-\frac{1}{2}}\, h_x,\\
Q\ \eqdef\ h\, S_x\ =\ h\, u_x\ +\ \sqrt{3\, \gamma}\, h^{-\frac{1}{2}}\, h_x,
\end{align*}
we have 
\begin{equation}\label{uh_x}
u_x\ =\ \frac{P\ +\ Q}{2\, h}, \qquad h_x\ =\ h^{\frac{1}{2}}\, \frac{Q\ -\ P}{2\, \sqrt{3\, \gamma}}.
\end{equation}
From the definition of $\mathcal{L}_{h}$ in \eqref{Ldef}, we obtain that
\begin{equation}\label{Psi}
\partial_x\, \mathcal{L}_{h}^{-1}\, \partial_x\, \Psi\ =\ -3\,  h^{-3}\, \Psi\ +\ 3\, \partial_x\, \mathcal{L}_{h}^{-1}\, \left(h \int_{-\infty}^x h^{-3} \Psi\right)
\end{equation}
for any smooth function $\Psi$ satisfying $\Psi(\pm \infty)=0$. 
Then, 
\begin{align} 
\mathscr{C}\ +\ \third\, h^3\, \partial_x\, \mathcal{L}_{h}^{-1}\,  \partial_x\, \mathscr{C}\ 
&=\  h^3\, \partial_x\, \mathcal{L}_{h}^{-1} \left( h \int_{-\infty}^x h^{-3}\, \mathscr{C} \right).
\end{align}
From \eqref{SGNc} and \eqref{SGN2b} we obtain 
\begin{align}\nonumber
\mathscr{R}\ &=\ - \third\, h^3 \left[ u_t\, +\, u\, u_x\, +\, 3\, \gamma h^{-2} h_x \right]_x\, +\ \mathscr{C}
\\
\label{CR}
&=\ \mathscr{C}\ +\ \third\, h^3\, \partial_x\, \mathcal{L}_{h}^{-1}  \partial_x\, \left\{ \mathscr{C}\, +\, F(h) \right\}\\ \label{RR}
&=\ h^3\, \partial_x\, \mathcal{L}_{h}^{-1} \left(h \int_{-\infty}^x h^{-3}\, \mathscr{C} \right)\ +\ \third\, h^3\, \partial_x\, \mathcal{L}_{h}^{-1}  \partial_x\, F(h).
\end{align}
Let the characteristics $X_{x},Y_{x}$ starting from $x$ defined as the solutions of the ordinary differential equations
\begin{align}
\frac{\mathrm{d}}{\mathrm{d}t}\, X_{x}(t)\ &=\ \eta(t,X_{x}(t)), \qquad X_{x}(0)\ =\ x,\\
\frac{\mathrm{d}}{\mathrm{d}t}\, Y_{x}(t)\ &=\ \, \lambda(t,Y_{x}(t)), \qquad\ Y_{x}(0)\ =\ x.
\end{align}
Differentiation \eqref{SGNR} with respect to $x$, and using \eqref{RR} we obtain the Ricatti-type equations
\begin{subequations}\label{Ricatti}
\begin{align}\label{Ricatti1}
\frac{\mathrm{d}^\lambda}{\mathrm{d}t}\, P\ \eqdef\ P_t\ +\ \lambda\, P_x\ =\ - {\textstyle \frac{1}{8\, h}}\, P^2\ +\ {\textstyle \frac{1}{8\, h}}\, Q^2\ -\ 3\, h^{-2}  \mathscr{R},\\ \label{Ricatti2}
\frac{\mathrm{d}^\eta}{\mathrm{d}t}\, Q\ \eqdef\ Q_t\ +\ \eta\, Q_x\ =\ - {\textstyle \frac{1}{8\, h}}\, Q^2\ +\ {\textstyle \frac{1}{8\, h}}\, P^2\ -\ 3\, h^{-2}  \mathscr{R},
\end{align}
\end{subequations}
where $\frac{\mathrm{d}^\lambda}{\mathrm{d}t}, \frac{\mathrm{d}^\eta}{\mathrm{d}t}$ denote the derivatives along the characteristics with the speed $\lambda,\eta$ respectively.
We prove below that the term $\mathscr{R}$ is bounded. Also, we obtain a bound of the integral of $P^2$ (respectively $Q^2$) on the characteristics $X_x$ (respectively $Y_x$). 
Then, the singularities given in Theorem \ref{thm:bu} appear due to the term $P^2$ in \eqref{Ricatti1} and/or the term $Q^2$ in \eqref{Ricatti2}.

\subsection{The choice of the approximated system}
In order to obtain a system that admits global smooth solutions, we linearise the negative quadratic terms on the right-hand side of \eqref{Ricatti} on the neighbourhood of $- \infty$.
For that purpose, let $\varepsilon > 0$ and we define as in \cite{wave2,wave3}
\begin{equation}\label{chidef}
\chi_\varepsilon (\zeta)\ \eqdef\ \left(\zeta\ +\ \frac{1}{\varepsilon} \right)^2 \mathds{1}_{(-\infty,-\frac{1}{\varepsilon}]} (\zeta)\ =\ 
\begin{cases}
\left(\zeta\ +\ \frac{1}{\varepsilon} \right)^2, & \zeta \leqslant -1/\varepsilon, \\
0, & \zeta > -1/\varepsilon.
\end{cases}
\end{equation}
Noting that \eqref{Ricatti} is like a derivative of \eqref{SGNxi}, then, adding terms to \eqref{Ricatti} will involve some primitive terms in \eqref{SGNxi} which are not uniquely defined and cannot vanish at $\infty$ and $-\infty$. 
That is why the system \eqref{SGNxi} will not be approximated simply by adding $\chi_\varepsilon$ to \eqref{Ricatti} as in \cite{wave2,wave3}.

Our goal is to obtain a system on the form 
\begin{align*}
h_t\ +\,\left[\, h\,u\,\right]_x\ &=\ h^+, \\ 
u_t\ +\ u\, u_x\ + 3\, \gamma h^{-2} h_x\ &=\ - \mathcal{L}_{h}^{-1}  \partial_x\, \left\{ \mathscr{C}\, +\, F(h) \right\}\ +\ u^+,
\end{align*}
where $h^+,u^+$ are suitable terms to be chosen.
As in Section \ref{sec_Riccati}, we obtain 
\begin{align*}
P_t\ +\ \lambda\, P_x\ &=\ - {\textstyle \frac{1}{8\, h}}\, P^2\ +\ {\textstyle \frac{1}{8\, h}}\, \chi_\varepsilon(P)\  +\ {\textstyle \frac{1}{8\, h}}\, Q^2\ -\ 3\, h^{-2}  \mathscr{R}\ +\ \underline{P},\\
Q_t\ +\ \eta\, Q_x\ &=\ - {\textstyle \frac{1}{8\, h}}\, Q^2\ +\ {\textstyle \frac{1}{8\, h}}\, \chi_\varepsilon(Q)\ +\ {\textstyle \frac{1}{8\, h}}\, P^2\ -\ 3\, h^{-2}  \mathscr{R}\ +\ \underline{Q},
\end{align*}
where 
\begin{align*}
\underline{P}\ \eqdef\ h\, (u^+)_x\ +\ u_x\, h^+\ -\ \sqrt{3\, \gamma}\, h^{-1/2}\, (h^+)_x\ +\ \half\, \sqrt{3\, \gamma}\, h^{-3/2}\, h_x\, h^+\ 
 -\  {\textstyle \frac{1}{8\, h}}\, \chi_\varepsilon(P),\\
\underline{Q}\ \eqdef\ h\, (u^+)_x\ +\ u_x\, h^+\ +\ \sqrt{3\, \gamma}\, h^{-1/2}\, (h^+)_x\ -\ \half\, \sqrt{3\, \gamma}\, h^{-3/2}\, h_x\, h^+\ 
 -\  {\textstyle \frac{1}{8\, h}}\, \chi_\varepsilon(Q).
\end{align*}
Due to the definition \eqref{chidef}, when $\zeta$ is near $-\infty$, the term $\chi_\varepsilon(\zeta) - \zeta^2$ behaves as a linear map. 
This prevents singularities from appearing in finite time.
From \eqref{TildeEneDef}, we have 
\begin{equation*}
\mathscr{E}\ =\ \half\, h\, u^2\ +\ \half\, g\, (h-\bar{h})^2\ +\ {\textstyle \frac{1}{12}}\, h\, P^2\ +\ {\textstyle \frac{1}{12}}\, h\, Q^2.
\end{equation*}
Then the energy equation \eqref{ene} becomes
\begin{align}\nonumber
\mathscr{E}_t\ +\ \mathscr{D}_x\ 
=&\ h\, u\, u^+\ +\ \half\, u^2\, h^+\ +\ g\, \left( h\, -\, \bar{h}\right) h^+\ +\ \left( {\textstyle \frac{1}{6}}\, h^2\, u_x^2\, +\, {\textstyle \frac{\gamma}{2\, h}}\, h_x^2 \right) h^+ \\  \nonumber
&+\ {\textstyle \frac{1}{6}}\, h\, P\, \underline{P}\ +\ {\textstyle \frac{1}{6}}\, h\, Q\, \underline{Q}\ +\ {\textstyle \frac{1}{48}}\, P \chi_\varepsilon(P)\ +\ {\textstyle \frac{1}{48}}\, Q \chi_\varepsilon(Q)\\  \nonumber
\leqslant &\ h\, u\, u^+\ +\ \half\, u^2\, h^+\ +\ g\, \left( h\, -\, \bar{h}\right) h^+\ +\ \left( {\textstyle \frac{1}{6}}\, h^2\, u_x^2\, +\, {\textstyle \frac{\gamma}{2\, h}}\, h_x^2 \right) h^+ \\   \label{Choice_ene_inq}
&+\ {\textstyle \frac{1}{6}}\, h\, P\, \underline{P}\ +\ {\textstyle \frac{1}{6}}\, h\, Q\, \underline{Q}.
\end{align}
The goal is to find $h^+$ and $u^+$ such that
\begin{itemize}
\item The right-hand side of \eqref{Choice_ene_inq} is a derivative of some quantity (i.e., $[ \cdots ]_x$), which will insure that $\int_\mathds{R} \mathscr{E} \mathrm{d}x$ is a decreasing function of time.
\item When $P,Q$ are large, we have $\underline{P} = \mathcal{O}(P)$ and $\underline{Q}= \mathcal{O}(Q)$. This insures (with Gronwall inequality) that no singularity will appear in finite time.
\end{itemize}
We can write the right-hand side of \eqref{Choice_ene_inq} as $T_1+T_2$ such that 
\begin{align*}
T_1\ 
=&\ g\, \left(h\, -\, \bar{h} \right) h^+\ +\ \gamma\, h_x\, \left(h^+\right)_x\ +\
{\textstyle \frac{\sqrt{3\, \gamma}}{48\, h^{1/2}}}\, h_x \left( \chi_\varepsilon(P)\, -\, \chi_\varepsilon(Q) \right)\\
=&\  g\, \left(h\, -\, \bar{h} \right) h^+\ +\ \left(h\, -\, \bar{h} \right)_x\, \left[  \gamma \left(h^+\right)_x\ +\
{\textstyle \frac{\sqrt{3\, \gamma}}{48\, h^{1/2}}}\, \left( \chi_\varepsilon(P)\, -\, \chi_\varepsilon(Q) \right) \right].
\end{align*}
Then, a sufficient condition to obtain $T_1 = [\cdots]_x$ is 
\begin{equation}\label{h+def}
g\, h^+\ =\ 
\left[  \gamma \left(h^+\right)_x\ +\
{\textstyle \frac{\sqrt{3\, \gamma}}{48\, h^{1/2}}}\, \left( \chi_\varepsilon(P)\, -\, \chi_\varepsilon(Q) \right) \right]_x.
\end{equation}
On another hand we have 
\begin{align*}
T_2\ 
&=\ \third\, h^3\, u_x\, \left( u^+ \right)_x\ +\ \half\, h^2\, u_x^2\, h^+\ +\ \half\, u^2\, h^+\ +\ h\, u\, u^+\ -\ {\textstyle \frac{1}{48} }\, h \left( \chi_\varepsilon(P)\, +\, \chi_\varepsilon(Q) \right) u_x\\
&=\ \left( \half\, u\, h^+\ +\ h\, u^+ \right) u\ +\ \left[\third\, h^3\, \left( u^+\right)_x\ +\ \half\, h^2\, u_x\, h^+\ -\ {\textstyle \frac{1}{48} }\, h \left( \chi_\varepsilon(P)\, +\, \chi_\varepsilon(Q) \right) \right] u_x
\end{align*}
Then, a sufficient condition to obtain $T_2 = [\cdots]_x$ is 
\begin{equation}\label{u+def}
\half\, u\, h^+\ +\ h\, u^+\ =\  \left[\third\, h^3\, \left( u^+\right)_x\ +\ \half\, h^2\, u_x\, h^+\ -\ {\textstyle \frac{1}{48} }\, h \left( \chi_\varepsilon(P)\, +\, \chi_\varepsilon(Q) \right) \right]_x.
\end{equation}
In next section we prove the global existence of smooth solutions of the approximated system, and we obtain some uniform estimates that do not depend on $\varepsilon$.

\section{Uniform estimates}\label{sec:UE}

In this section, we consider $\gamma>0,h>0$, and $h_0-\bar{h},u_0 \in H^1$ such that $\int_\mathds{R} {\mathscr{E}}_0\, \mathrm{d}x < \sqrt{g \gamma} \bar{h}^2$. Let also $j_\varepsilon$ be a Friedrichs mollifier, we define  $h_0^\varepsilon \eqdef ((h_0-\bar{h}) \ast j_\varepsilon) + \bar{h}$ and $u_0^\varepsilon \eqdef (u_0 \ast j_\varepsilon)$ where $(f \ast g)(x) \eqdef \int_\mathds{R} f(x-x')g(x') \mathrm{d}x'$. Using that $\|(h_0-h_0^\varepsilon, u_0 - u_0^\varepsilon)\|_{H^1} \to 0$ as $\varepsilon \to 0$, we can prove
\begin{equation}\label{Limene}
\lim_{\varepsilon \to 0}\, \int_\mathds{R} {\mathscr{E}}_0^\varepsilon\, \mathrm{d}x\ =\ \int_\mathds{R} {\mathscr{E}}_0\, \mathrm{d}x\ <\ \sqrt{g\, \gamma}\, \bar{h}^2,
\end{equation}
which implies that there exists $\varepsilon_0>0$ such that 
\begin{equation}\label{Eep}
\int_\mathds{R} {\mathscr{E}}_0^\varepsilon\, \mathrm{d}x\ \leqslant\ E\ \eqdef\ \half\, \int_\mathds{R} {\mathscr{E}}_0\, \mathrm{d}x\ +\ \half\, g\, \gamma\, \bar{h}^2  \qquad \forall \varepsilon\ \leqslant\ \varepsilon_0.
\end{equation}
Following the arguments of the previous section (see \eqref{h+def} and \eqref{u+def}), we consider the system
\begin{subequations}\label{Appsys}
\begin{gather}\label{Appsys1}
h_t^\varepsilon\ +\,\left[\, h^\varepsilon\,u^\varepsilon\,\right]_x\ =\ \mathscr{A}^\varepsilon_x, \\ \label{Appsys2}
u^\varepsilon_t\ +\ u^\varepsilon\, u^\varepsilon_x\ + 3\, \gamma (h^\varepsilon)^{-2} h^\varepsilon_x\ =\ - \mathcal{L}_{h^\varepsilon}^{-1}  \partial_x\, \left\{ \mathscr{C}^\varepsilon\, +\, F(h^\varepsilon) \right\}\ +\ \mathscr{B}^\varepsilon,\\ 
u^\varepsilon(0,\cdot)\ =\ u^\varepsilon_0\ \eqdef\ j_\varepsilon \ast u_0, \qquad h^\varepsilon(0,\cdot)\ =\ h^\varepsilon_0\ \eqdef\ j_\varepsilon \ast \left(h_0\, -\, \bar{h}\right) +\, \bar{h},
\end{gather}
\end{subequations}
where 
\begin{align} \nonumber
\mathscr{A}^\varepsilon\ &\eqdef\  \left( g\  -\ \gamma\, \partial_x^2 \right)^{-1}
\left\{ \frac{\sqrt{3\, \gamma}}{48\, (h^\varepsilon)^{1/2}}\, \left(\chi_\varepsilon(P^\varepsilon)\, -\, \chi_\varepsilon(Q^\varepsilon) \right) \right\},\\ \label{Adef}
&=\ \mathfrak{G} \ast \left\{ \frac{\sqrt{3\, \gamma}}{48\, (h^\varepsilon)^{1/2}}\, \left(\chi_\varepsilon(P^\varepsilon)\, -\, \chi_\varepsilon(Q^\varepsilon) \right) \right\},\\ \label{Bdef}
\mathscr{B}^\varepsilon\ &\eqdef\ \mathcal{L}_{h^\varepsilon}^{-1} \left\{ - \half\, u^\varepsilon\, \mathscr{A}^\varepsilon_x\ +\ \partial_x \left\{ \half\, (h^\varepsilon)^2\, u^\varepsilon_x\, \mathscr{A}^\varepsilon_x\ -\ {\textstyle \frac{1}{48}}\, h^\varepsilon \left(\chi_\varepsilon(P^\varepsilon)\, +\, \chi_\varepsilon(Q^\varepsilon) \right) \right\} \right\},
\end{align}
with $\mathfrak{G}$ is defined as
\begin{equation}
\mathfrak{G}\ \eqdef\ {\textstyle \frac{1}{2\, \gamma}}\, \exp\left\{- {\textstyle \frac{g}{\gamma}}\, |\cdot|\right\}.
\end{equation}
Differentiation \eqref{Appsys} with respect to $x$ we obtain
\begin{subequations}\label{Ricatti_ep}
\begin{flalign}\label{Ricatti1_ep}
\frac{\mathrm{d}^\lambda}{\mathrm{d}t} P^\varepsilon \eqdef P^\varepsilon_t + \lambda^\varepsilon P^\varepsilon_x &= - {\textstyle \frac{1}{8\, h^\varepsilon}} (P^\varepsilon)^2 + {\textstyle \frac{1}{8\, h^\varepsilon}} \chi_\varepsilon(P^\varepsilon) + {\textstyle \frac{1}{8\, h^\varepsilon}} (Q^\varepsilon)^2 - {\textstyle \frac{1}{2\, h^\varepsilon}} \mathscr{A}^\varepsilon_x  P^\varepsilon + \mathscr{M}^\varepsilon, \\ \label{Ricatti2_ep}
\frac{\mathrm{d}^\eta}{\mathrm{d}t} Q^\varepsilon \eqdef\ Q^\varepsilon_t + \eta^\varepsilon Q^\varepsilon_x &= - {\textstyle \frac{1}{8\, h^\varepsilon}} (Q^\varepsilon)^2 + {\textstyle \frac{1}{8\, h^\varepsilon}} \chi_\varepsilon(Q^\varepsilon) + {\textstyle \frac{1}{8\, h^\varepsilon}} (P^\varepsilon)^2  - {\textstyle \frac{1}{2\, h^\varepsilon}} \mathscr{A}^\varepsilon_x  Q^\varepsilon + \mathscr{N}^\varepsilon,
\end{flalign}
\end{subequations}
with
\begin{align}
\mathscr{M}^\varepsilon\ &\eqdef\ -\ 3\, (h^\varepsilon)^{-2}  \mathscr{R}^\varepsilon\ +\ \mathscr{V}^\varepsilon_1\ -\ \mathscr{V}^\varepsilon_2, \qquad
\mathscr{N}^\varepsilon\  \eqdef\ -\ 3\, (h^\varepsilon)^{-2}  \mathscr{R}^\varepsilon\ +\ \mathscr{V}^\varepsilon_1\ +\ \mathscr{V}^\varepsilon_2, \\ \label{M1def}
\mathscr{V}^\varepsilon_1\ &\eqdef\ \half\, h^\varepsilon \partial_x \mathcal{L}_{h^\varepsilon}^{-1}\! \left\{\!- u^\varepsilon \mathscr{A}^\varepsilon_x + h^\varepsilon\!\! \int_{-\infty}^x\! \left[3 (h^\varepsilon)^{-1} u^\varepsilon_x\, \mathscr{A}^\varepsilon_x - {\textstyle \frac{1}{8\, (h^\varepsilon)^2}}\! \left(\chi_\varepsilon(P^\varepsilon) + \chi_\varepsilon(Q^\varepsilon) \right)\!  \right]\! \mathrm{d}y\! \right\}\! ,\\ \label{M2def}
\mathscr{V}^\varepsilon_2\ &\eqdef\ {\textstyle \frac{g}{16}}\, (h^\varepsilon)^{-1/2}\,\mathfrak{G} \ast \left\{ (h^\varepsilon)^{-1/2} \left(\chi_\varepsilon(P^\varepsilon)\, -\, \chi_\varepsilon(Q^\varepsilon) \right) \right\}\ =\ {\textstyle \frac{3\, g}{\sqrt{3\, \gamma}}}\, (h^\varepsilon)^{-1/2}\, \mathscr{A}^\varepsilon.
\end{align}
Smooth solutions of \eqref{Appsys} satisfy the energy equation  (see Appendix \ref{App:B})
\begin{equation}\label{energyequationep}
\mathscr{E}^\varepsilon_t\ +\ \tilde{\mathscr{D}}^\varepsilon_x\ =\  {\textstyle \frac{1}{48}}\, P^\varepsilon \chi_\varepsilon(P^\varepsilon)\ +\ {\textstyle \frac{1}{48}}\, Q^\varepsilon \chi_\varepsilon(Q^\varepsilon)\ \leqslant\ 0,
\end{equation}
where 
\begin{equation*}
\tilde{\mathscr{D}}^\varepsilon \eqdef u^\varepsilon \mathscr{E}^\varepsilon\, +\, u^\varepsilon\! \left( \mathscr{R}^\varepsilon + \half\, g\, (h^\varepsilon)^2 - \half\, g\, \bar{h}^2 \right)\, +\, \gamma\, h^\varepsilon h^\varepsilon_x u_x^\varepsilon\, -\, \third\, (h^\varepsilon)^2 u^\varepsilon\, \mathscr{V}^\varepsilon_1\,
  -\,  {\textstyle \frac{\sqrt{3\, \gamma}}{3}}\, (h^\varepsilon)^{1/2}\, \mathscr{V}^\varepsilon_2 \left(h^\varepsilon-\bar{h}\right)\! .
\end{equation*}
The first result on this section is the global well-posedness of \eqref{Appsys}.
\begin{thm}\label{thm:ex:ep}
Let $\bar{h}>0$, $(h_0 - \bar{h}, u_0) \in \mathfrak{D}$ and $\varepsilon \in (0,\varepsilon_0]$ then there exists a global smooth solution $(h^\varepsilon - \bar{h}, u^\varepsilon) \in C(\mathds{R}^+,H^3(\mathds{R})) \cap C^1(\mathds{R}^+,H^2(\mathds{R}))$ of \eqref{Appsys} and for all $t >0$ we have
\begin{equation}\label{energyconservationep}
\int_\mathds{R} \mathscr{E}^\varepsilon\, \mathrm{d}x\ 
-\ \int_0^t \int_\mathds{R} {\textstyle \frac{1}{48}}\, \left( P^\varepsilon\, \chi_\varepsilon(P^\varepsilon)\ +\ Q^\varepsilon\, \chi_\varepsilon(Q^\varepsilon) \right)\, \mathrm{d}x\, \mathrm{d}t\
=\  \int_\mathds{R} \mathscr{E}^\varepsilon_0\, \mathrm{d}x .
\end{equation}
Moreover, there exist $A,B>0$ depending only on $\bar{h},\gamma,g$ and $E$ such that for any $t >0$, $x_2 \in \mathds{R}$,  and for $x_1 \in (-\infty, x_2)$ the solution of $X_{x_1}(t)=Y_{x_2}(t)$ (see Figure \ref{XY}) we have 
\begin{align}\label{PQ^2}
\int_\tau^t \left[P^\varepsilon(s,X_{x_1}(s))\right]^2\, \mathrm{d}s\ +\ \int_\tau^t \left[Q^\varepsilon(s,Y_{x_2}(s))\right]^2\, \mathrm{d}s\ \leqslant\ A\, (t\, -\, \tau)\ +\ B \quad \forall \tau \in [0,t].
\end{align}
\end{thm}
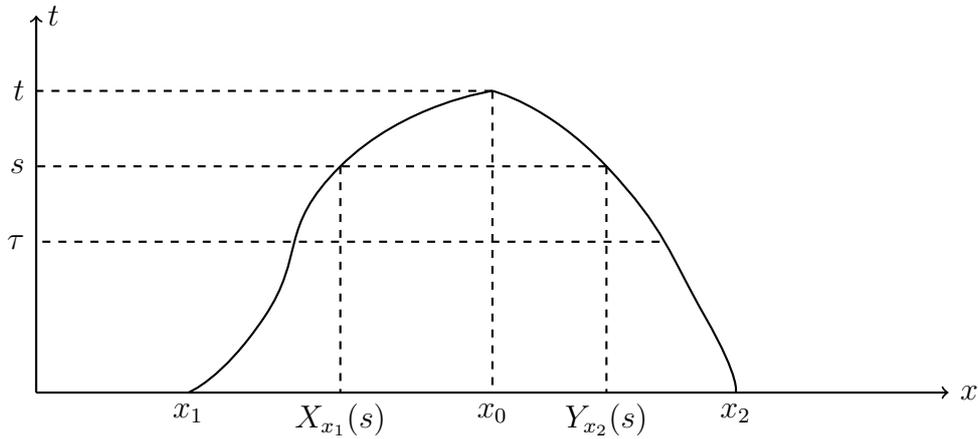
\begin{figure}[!ht]
\begin{tikzpicture}[thick, transform canvas={scale=1}, shift={(0,-5)}]
\draw[->] (-6,0) -- (6,0) node[right]{$x$};
\draw[->] (-6,0) -- (-6,5) node[right]{$t$};

\draw [black] plot [smooth, tension=1] coordinates { (3.2,0) (2.8,1) (1.5,3) (0,4)};

\draw [black] plot [smooth, tension=1] coordinates { (-4,0) (-3,1) (-2,3) (0,4)};

\draw[dashed] (0,4) -- (-6,4) node[left]{$t$};
\draw[dashed] (0,4) -- (0,0) node[below]{$x_0$};
\draw[dashed] (1.5,3) -- (-6,3) node[left]{$s$};
\draw[dashed] (1.5,3) -- (1.5,0) node[below]{$Y_{x_2}(s)$};
\draw[dashed] (-2,3) -- (-2,0) node[below]{$X_{x_1}(s)$};
\draw[dashed] (2.2,2) -- (-6,2) node[left]{$\tau$};


\draw[] (3.2,0) -- (3.2,0) node[below]{$x_2$};
\draw[] (-4,0) -- (-4,0) node[below]{$x_1$};

\end{tikzpicture}
\vspace{5.7cm}
\caption{Characteristics.}
\label{XY}
\end{figure}
In order to prove Theorem \ref{thm:ex:ep}, we need to prove the invertibility of the operator $\mathcal{L}_h$ and to obtain some estimates of its inverse.
\begin{lem}\label{Inverseesitimates}
Let  $0<  h \in H^1(\mathds{R})+\bar{h}$ with $h^{-1} \in L^\infty$. Then the operator $\mathcal{L}_{h}$ is an isomorphism from $H^2$ to $L^2$. Moreover, if $\psi \in {C}_{\mathrm{lim}} \eqdef \{f\in {C}(\mathds{R}),\ f(\infty), f(-\infty) \in \mathds{R} \}$, then $\mathcal{L}_{h}^{-1}\, \psi$ is well defined and there exists constants $C=C\left(\bar{h},\left\| h^{-1} \right\|_{L^\infty}, \|h\|_{L^{\infty}} \right)>0$ 
such that 
\begin{align}\label{estimate1}
\left\|\mathcal{L}_{h}^{-1}\, \psi\right\|_{W^{1,\infty}}\ &\leqslant\ C\, \left\|\psi\right\|_{L^\infty},\\ \label{estimate1.5}
\left|\partial^2_x\, \mathcal{L}_{h}^{-1}\, \psi\right|(x)\ &\leqslant\ C\, \left(1\, +\, \left|h_x(x)\right| \right)\, \left\|\psi\right\|_{L^\infty} \quad  \forall x \in \mathds{R}, \\ \label{estimate2}
\left\|\mathcal{L}_{h}^{-1}\, \psi\right\|_{H^1}\ &\leqslant\ C\, \left\|\psi\right\|_{L^1}, \\ \label{estimate3}
\left\|\mathcal{L}_{h}^{-1}\, \partial_x\, \psi\right\|_{L^\infty}\ &\leqslant\ C\, \left\|\psi\right\|_{L^1}, \\ \label{estimate3.5}
\left\|\partial_x\, \mathcal{L}_{h}^{-1}\, \partial_x\, \psi\right\|_{L^\infty}\ &\leqslant\ C \left[ \left\|\psi\right\|_{L^1}\ +\ \left\|\psi\right\|_{L^\infty} \right],\\ \label{estimate4}
\left\| \mathcal{L}_{h}^{-1}\, \partial_x\, \psi\right\|_{L^2}\ 
&\leqslant\ C\, \left\|\psi\right\|_{L^1}\, \left( \left\| h_x \right\|_{L^2} +\ 1 \right),\\  \label{estimate6}
\left\| \mathcal{L}_{h}^{-1}\, \partial_x\, \psi \right\|_{H^1}\, +\, \left\| \mathcal{L}_{h}^{-1}\, \psi \right\|_{H^1}\, &\leqslant\ C \, \|\psi\|_{L^2},\\ \label{estimate7}
\left\| \mathcal{L}_{h}^{-1}\,  \psi \right\|_{W^{1,\infty}}\, \leqslant\
\left\| \mathcal{L}_{h}^{-1}\,  \psi \right\|_{H^2}\,  &\leqslant\ C \left[1\, +\, \left\|h_x \right\|_{L^2}^2 \right]\left\| \psi \right\|_{L^2}, \\ \label{estimate8}
\left|\partial^2_x\, \mathcal{L}_{h}^{-1}\, \psi\right|(x)\ &\leqslant\ C \left[ \left( 1\, +\, \|h_x\|_{L^2} \right) \left( 1\, +\,  |h_x|(x) \right)  \|\psi\|_{L^2}\ +\ |\psi|(x)  \right].
\end{align}
Also, if $h - \bar{h} \in H^2(\mathds{R})$ we have 
\begin{subequations}\label{estimate5}
\begin{align}
\left\|\mathcal{L}_h^{-1}\, \partial_x\, \varphi \right\|_{H^{3}}\  
&\leqslant\ C \left[ \left(1\, +\, \|h_x\|^2_{L^\infty} \right) \left\|\varphi \right\|_{H^2}\ +\ \left\|h\, -\, \bar{h} \right\|_{H^2}\, 
\left\|\mathcal{L}_h^{-1}\, \partial_x\, \varphi \right\|_{W^{1,\infty}}\right],\\
\left\|\mathcal{L}_h^{-1}\, \psi \right\|_{H^{3}}\  
&\leqslant\ C \left[ \left(1\, +\, \|h_x\|^2_{L^\infty} \right) \left\|\psi \right\|_{H^2}\ +\ \left\|h\, -\, \bar{h} \right\|_{H^2}\, 
\left\|\mathcal{L}_h^{-1}\, \psi \right\|_{W^{1,\infty}} \right].
\end{align}
\end{subequations}
Moreover, there exists a constant $\tilde{C}=\tilde{C}(\gamma,g)$ such that 
\begin{equation}\label{estimate5.5}
\left\| \left( g\  -\ \gamma\, \partial_x^2 \right)^{-1} \psi \right\|_{H^3}\ \leqslant\ \tilde{C} \left\| \psi \right\|_{H^1}, 
\qquad
\left\| \partial_x \left( g\  -\ \gamma\, \partial_x^2 \right)^{-1} \psi \right\|_{H^3}\ \leqslant\ \tilde{C} \left\| \psi \right\|_{H^2}.
\end{equation}
\end{lem}
The proof of \eqref{estimate1}, \eqref{estimate3.5} and  \eqref{estimate5} is inspired by \cite{liu2019well}.
\proof \textbf{Step 0.} Let $(\cdot,\cdot)$ be the scalar product in $L^2$. Defining the bilinear map $a : H^1 \times H^1 \to \mathds{R}$
\begin{equation*}
a(u,v)\ \eqdef\ \left(h\, u, v \right)\, +\, \third \left(h^3\, u_x, v_x \right).
\end{equation*}
It is easy to check that $a$ is continuous and   coercive. Then, Lax--Milgram theorem insures the existence of a continuous bijective linear operator $J : H^1 \to H^{-1}$ satisfying 
\begin{equation*}
a(u,v)\ =\ 	\langle J\, u, v \rangle_{H^{-1} \times H^1} \qquad \forall u,v \in H^1.
\end{equation*}
If $J u \in L^2$, an integration by parts shows that $ \left( h^3 u_x \right)_x = h u -J u \in L^2$ and $J=\mathcal{L}_h$, this implies that $u \in H^2$ which finishes the proof that $\mathcal{L}_h$ is an isomorphism from $H^2$ to $L^2$.

Defining now $C_0 \eqdef \{f\in {C},\ f(\pm \infty)=0 \}$, using that $L^2 \cap {C}_0$ is dense in ${C}_0$ one can define $\mathcal{L}_h^{-1}$ on ${C}_0$.
If $\varphi$ is in ${C}_{\mathrm{lim}}$, we use the change of functions (see Lemma 4.4 in \cite{liu2019well}) 
\begin{equation*}
\varphi_0(x)\ \eqdef\ \varphi(x)\ -\   \mathcal{L}_h {\textstyle \frac{1}{\bar{h}}} \left( \varphi(-\infty)\ +\ \left( \varphi(\infty)\ -\ \varphi(-\infty) \right)\, \frac{\mathrm{e}^x}{1+\mathrm{e}^x} \right) \in C_0,
\end{equation*}
the operator $\mathcal{L}_h^{-1}$ can be defined as
\begin{equation*}
\mathcal{L}_h^{-1}\, \varphi\ \eqdef\ \mathcal{L}_h^{-1}\, \varphi_0\ +\  {\textstyle \frac{1}{\bar{h}}}\, \left( \varphi(-\infty)\ +\ \left( \varphi(\infty)\ -\ \varphi(-\infty) \right)\, \frac{\mathrm{e}^x}{1+\mathrm{e}^x} \right).
\end{equation*}

\textbf{Step 1.}
Let 
\begin{equation}\label{phispi}
\psi\ =\ \mathcal{L}_{h} u =\ h\, u\ -\ \third\, \left(h^3\,  u_x \right)_x,
\end{equation}
using the change of variables 
\begin{equation}\label{cov}
z\ =\ \int \frac{\mathrm{d}x}{h^3}\, ,
\end{equation} 
we obtain 
\begin{equation}\label{maximump}
\psi\ =\ h\, u\ -\ {\textstyle \frac{1}{3\, h^3}}\, u_{zz}.
\end{equation}
The maximum principle insures that $\|u\|_{L^\infty} \leqslant C \|\psi\|_{L^\infty}$ which implies with \eqref{maximump} that 
\begin{equation}\label{zz}
\|u_{zz}\|_{L^\infty} \leqslant C \|\psi\|_{L^\infty}.
\end{equation}
Using Landau--Kolmogorov inequality we obtain $\|u_{z}\|_{L^\infty} \leqslant C \|\psi\|_{L^\infty}$. Using again the change of variables \eqref{cov} we get $\|u_{x}\|_{L^\infty} \leqslant C \|\psi\|_{L^\infty}$ which completes the proof of \eqref{estimate1}.
The estimate \eqref{estimate1.5} follows directly from the change of variables \eqref{cov}, \eqref{zz} and \eqref{estimate1}.

Multiplying \eqref{phispi} by $u$ and using and integration by parts one obtains
\begin{equation}
\left\| u \right\|_{H^1}^2\ \leqslant\ C\, \| \psi\|_{L^1}\, \left\|u\right\|_{L^\infty}.
\end{equation}
The inequality \eqref{estimate2} follows directly using the embedding $H^1 \hookrightarrow L^\infty$.
Using \eqref{estimate1} and 
\begin{equation}
\mathcal{L}_{h}^{-1}\, \partial_x\, \psi\ =\ 
-3\, \int_{-\infty}^x \left(h^{-3}\, \psi\right)\, 
+\ 3\, \mathcal{L}_{h}^{-1}\, \left(h \int_{-\infty}^x h^{-3}\,  \psi \right)
\end{equation}
we obtain \eqref{estimate3} and \eqref{estimate3.5}.
Using the definition of $\mathcal{L}_h$ we obtain
\begin{align} \nonumber
\mathcal{L}_{h}^{-1}\, \partial_x\, \psi\ 
=&\ \mathcal{L}_{h}^{-1}\, \partial_x\, h^3\, \mathcal{L}_h\, \mathcal{L}_h^{-1}\, h^{-3}\, \psi\\ \nonumber
=&\ \mathcal{L}_{h}^{-1}\, \partial_x \left[ h^4\, \mathcal{L}_h^{-1}\, h^{-3}\, \psi\ -\ \third\, h^3\, \partial_x\, h^3\, \partial_x\,\mathcal{L}_h^{-1}\, h^{-3}\, \psi \right]\\ \nonumber
=&\ \mathcal{L}_{h}^{-1}\, \left[4\, h^3\, h_x\, \mathcal{L}_h^{-1}\, h^{-3}\, \psi\ +\ h^4\, \partial_x\, \mathcal{L}_h^{-1}\, h^{-3}\, \psi\ -\ \third\, \partial_x\, h^3\, \partial_x\, h^3\, \partial_x\,\mathcal{L}_h^{-1}\, h^{-3}\, \psi \right]\\ \label{Lpp}
=&\ \mathcal{L}_{h}^{-1}\, \left[4\, h^3\, h_x\, \mathcal{L}_h^{-1}\, h^{-3}\, \psi\right]\, +\ h^3\, \partial_x\,\mathcal{L}_h^{-1}\, h^{-3}\, \psi .
\end{align}
The inequality \eqref{estimate4} follows from \eqref{estimate2} and the Cauchy--Schwarz inequality.

Let $\mathcal{L}_h u = \psi + \varphi_x$, then
\begin{align}
\|u\|_{H^1}^2\ 
&=\ (u,u)\ +\ (u_x,u_x) \nonumber \\ 
&\leqslant\ C \left[ (h\, u,u)\ +\ \third\, (h^3\, u_x,u_x) \right] \nonumber \\
&=\ C\, (\mathcal{L}_h u,u)\ =\ C \left[ (\psi,u)\ -\ (\varphi,u_x) \right]  \nonumber\\
&\leqslant\ C\, \|u\|_{H^1}\, \left( \|\psi\|_{L^2}\ +\ \|\varphi\|_{L^2} \right), \nonumber
\end{align}
which implies that 
\begin{equation}\label{s=0}
\|u\|_{H^1}\  \leqslant\ C \left( \|\psi\|_{L^2}\ +\ \|\varphi\|_{L^2}\right).
\end{equation}
Taking $\psi=0$ (respectively $\varphi=0$) we obtain \eqref{estimate6}.
Replacing $h^{-3} \psi$ by $\psi$ in  \eqref{Lpp}, we multiply by $h^{-3}$ and we differentiate with respect to $x$ to obtain
\begin{align}\label{Lppx}
\partial_x^2\,\mathcal{L}_h^{-1}\, \psi\ 
=&\  -3\, h^{-2}\, h_x \left[\mathcal{L}_{h}^{-1}\, \partial_x\, h^3\, \psi  -
\mathcal{L}_{h}^{-1}\, \left[4\, h^3\, h_x\, \mathcal{L}_h^{-1}\, \psi\right] \right]\\
&+\ h^{-3}\, \partial_x\, \mathcal{L}_{h}^{-1}\, \partial_x\, h^3\, \psi\
-\
h^{-3}\, \partial_x\, \mathcal{L}_{h}^{-1}\, \left[4\, h^3\, h_x\, \mathcal{L}_h^{-1}\, \psi\right].
\end{align}
Using \eqref{estimate6} and the embedding $H^1 \hookrightarrow  L^\infty$ we obtain 
\begin{align*}
\left\| \partial_x^2\, \mathcal{L}_{h}^{-1}\,  \psi \right\|_{L^2}\,  
\leqslant &\  C \left\|h_x \right\|_{L^2} \left[ \left\|\mathcal{L}_{h}^{-1}\, \partial_x\, h^3\, \psi \right\|_{H^1}\ +\ \left\| \mathcal{L}_{h}^{-1}\, \left[4\, h^3\, h_x\, \mathcal{L}_h^{-1}\, \psi\right] \right\|_{H^1} \right]\\
&+\ C \left\| \partial_x\, \mathcal{L}_{h}^{-1}\, \partial_x\, h^3\, \psi\ \right\|_{L^2}\ +\ C \left\|\partial_x\, \mathcal{L}_{h}^{-1}\, \left[4\, h^3\, h_x\, \mathcal{L}_h^{-1} \psi \right]  \right\|_{L^2}\\
\leqslant &\ C \left\|h_x \right\|_{L^2} \left[ \left\| \psi \right\|_{L^2}\ +\ \left\|h_x \right\|_{L^2} \left\|  \mathcal{L}_h^{-1}\, \psi\right\|_{H^1} \right]\\
&+\ C \left\| \psi \right\|_{L^2}\ +\ C \left\|h_x \right\|_{L^2} \left\|  \mathcal{L}_h^{-1}\, \psi\right\|_{H^1}\\
\leqslant &\ C \left[1\, +\, \left\|h_x \right\|_{L^2}^2 \right]\left\| \psi \right\|_{L^2}.
\end{align*}
This with \eqref{estimate6} imply \eqref{estimate7}.

Differentiating now \eqref{Lpp} with respect to $x$, using the definition of $\mathcal{L}_h$ and replacing $h^{-3} \psi$ by $\psi$ we obtain 
\begin{equation*}
\partial_x\, \mathcal{L}_{h}^{-1}\, \partial_x\, h^3\, \psi\ =\ \partial_x\, \mathcal{L}_{h}^{-1}\, \left[4\, h^3\, h_x\, \mathcal{L}_h^{-1}\,  \psi\right]\, 
+\ 3\, h\, \mathcal{L}_h^{-1}\, \psi\ -\ 3\, \psi.
\end{equation*}
Then \eqref{Lppx} becomes
\begin{align*}
\partial_x^2\,\mathcal{L}_h^{-1}\, \psi\ 
=\  -3\, h^{-2}\, h_x \left[\mathcal{L}_{h}^{-1}\, \partial_x\, h^3\, \psi  -
\mathcal{L}_{h}^{-1}\, \left[4\, h^3\, h_x\, \mathcal{L}_h^{-1}\, \psi\right] \right]\,
+\  3\, h^{-2}\, \mathcal{L}_h^{-1}\, \psi\ -\ 3\, h^{-3}\, \psi.
\end{align*}
Then, using \eqref{estimate6} we obtain \eqref{estimate8}.

\textbf{Step 2.}
Using $\mathcal{L}_h u = \psi + \varphi_x$ and the Young inequality $ab \leqslant {\textstyle \frac{1}{2\/ \alpha}} a^2 + {\textstyle \frac{\alpha}{2}} b^2$ with $\alpha>0$ we obtain
\begin{align}
\|u_x\|_{H^1}^2\ 
&=\ (u_x,u_x)\ +\ (u_{xx},u_{xx}) \nonumber \\ 
&\leqslant\ C \left[ (h\, u_x,u_x)\ +\ \third\, (h^3\, u_{xx},u_{xx}) \right]\nonumber \\
&=\ C \left[ -(h\, u, u_{xx})\ -\ (h_x\, u, u_x)\ +\ \third\, \left( \left(h^3\, u_x \right)_x\, -\, (h^3)_x u_x, u_{xx} \right)\right] \nonumber\\
&=\ C \left[  -(\mathcal{L}_h\, u, u_{xx})\ -\ (h_x\, u, u_x)\ -\ \left( h^2\, h_x\, u_x, u_{xx} \right) \right] \nonumber\\
& \leqslant\ C \left[ \alpha\, \|u_{xx}\|_{L^2}^2\ +\ {\textstyle \frac{1}{\alpha}}\, \|\mathcal{L}_h\, u\|_{L^2}^2\ +\ C_\alpha \left(1\, +\, \|h_x\|_{L^\infty}^2 \right) \|u\|_{H^1}^2 \right].\nonumber
\end{align}
Taking $\alpha>0$ small enough we obtain that 
\begin{equation*}
\|u_x\|_{H^1}^2\ \leqslant\ C \left[ \|\mathcal{L}_h\, u\|_{L^2}^2\ +\, \left(1\, +\, \|h_x\|_{L^\infty}^2 \right) \|u\|_{H^1}^2\right],
\end{equation*}
then 
\begin{equation*}
\|u_x\|_{H^1}\ \leqslant\ C \left[ \|\mathcal{L}_h\, u\|_{L^2}\ +\, \left(1\, +\, \|h_x\|_{L^\infty} \right) \|u\|_{H^1}\right].
\end{equation*}
Taking $\psi=0$ (respectively $\varphi=0$) and using \eqref{s=0}, we obtain 
\begin{equation}\label{H2est}
\left\|\mathcal{L}_h^{-1}\, \partial_x\, \varphi \right\|_{H^2} \leqslant\, C \left(1 + \|h_x\|_{L^\infty} \right) \|\varphi\|_{H^1}, \qquad
\left\|\mathcal{L}_h^{-1}\, \psi \right\|_{H^{2}}  \leqslant\, C \left(1 + \|h_x\|_{L^\infty} \right) \|\psi\|_{L^2}.
\end{equation}
Let $\Lambda$ be defined as $\widehat{\Lambda f}=(1+\xi^2)^\frac{1}{2} \hat{f}$. Since $\mathcal{L}_h u = \psi + \varphi_x$, we have
$$ \mathcal{L}_h\, \Lambda^2 u\ =\ [h, \Lambda^2]\, u\ +\ \Lambda^2 \psi\ +\ \partial_x\, \left\{ -\, \third\, [h^3, \Lambda^2]\, u_x\ +\ \Lambda^2 \varphi \right\}. $$ 
Defining $\tilde{u}=\Lambda^2 u$, $\tilde{\psi}=[h, \Lambda^2] u +  \Lambda^2 \psi $ and $\tilde{\varphi}= - \third\, [h^3, \Lambda^2] u_x +  \Lambda^2 \varphi $ and using \eqref{s=0}, \eqref{Commutator} we obtain
\begin{align*}
\|\Lambda^2 u\|_{H^1}\ 
&\leqslant\ C \left[ \left\|[h, \Lambda^2]\, u \right\|_{L^2}\ +\ \left\|[h^3 , \Lambda^2]\, u_x \right\|_{L^2}\ +\ \left\|\psi \right\|_{H^2}\ +\ \left\|\varphi \right\|_{H^2} \right] \\
&\leqslant\ C \left[ \|h_x\|_{L^\infty}\, \| u\|_{H^2}\ +\ \|h\, -\, \bar{h}\|_{H^2}\, \|u\|_{W^{1,\infty}}\ +\ \left\|\psi \right\|_{H^2}\ +\ \left\|\varphi \right\|_{H^2} \right]. 
\end{align*}
Taking $\psi=0$ (respectively $\varphi=0$) and using \eqref{s=0} with \eqref{H2est}, we obtain \eqref{estimate5}.

\textbf{Step 3.} It remains only to prove the inequalities \eqref{estimate5.5}. Since the operator $\left( g\  -\ \gamma\, \partial_x^2 \right)^{-1}$ is nothing but a convolution with the function $\mathfrak{G}$, the result follows directly using the Young inequality.
\qed

\begin{lem}\label{lem:bounded}
Let  $ (h-\bar{h},u) \in H^1(\mathds{R})$ such that $\int_\mathds{R} \mathscr{E}\, \mathrm{d}x \leqslant E < \sqrt{g \gamma} \bar{h}^2$, then, there exists a constant $C=C( \gamma, \bar{h}, E)>0$ independent on $\varepsilon$ and $h$ such that 
\begin{gather} \label{Estimate1}
\left\| \mathcal{L}_{h}^{-1}  \partial_x\, \mathscr{C}  \right\|_{L^\infty\left( \mathds{R}\right)}\ +\ 
\left\| \mathscr{R}  \right\|_{L^\infty\left( \mathds{R}\right)}\ 
\leqslant\ C, \\ \label{Estimate2}
 \int_\mathds{R} \left( \chi_\varepsilon (P)\ +\ \chi_\varepsilon (Q) \right)\, \mathrm{d}x\ 
\leqslant\  C,\\ \label{Estimate3}
\left\| \mathcal{L}_{h}^{-1}  \partial_x\, \left\{   h \left(\chi_\varepsilon(P)\, +\, \chi_\varepsilon(Q) \right) \right\}  \right\|_{L^\infty\left( \mathds{R}\right)}\
\leqslant\   C,\\ \label{Estimate3.5}
\left\| \left( \mathcal{L}_{h}^{-1}  \partial_x\, \left\{ h^2\, u_x\, \mathscr{A}^\varepsilon_x \right\}, \mathcal{L}_{h}^{-1} \left\{u\, \mathscr{A}^\varepsilon_x \right\} \right) \right\|_{L^\infty\left( \mathds{R}\right)}\ \leqslant\  C,\\ \label{Estimate4}
\left\|\mathscr{A}_x^\varepsilon \right\|_{L^2}\ +\, \left\| \left(\mathscr{A}^\varepsilon, \mathscr{A}^\varepsilon_x, \mathscr{B}^\varepsilon, \mathscr{V}^\varepsilon_1, \mathscr{V}^\varepsilon_2 \right) \right\|_{L^\infty\left(\mathds{R}\right)}\
\leqslant\  C,
\end{gather}
where $\mathscr{A}^\varepsilon$, $\mathscr{B}^\varepsilon$, $\mathscr{V}_1^\varepsilon$ and $\mathscr{V}_2^\varepsilon$ are defined as in \eqref{Adef}, \eqref{Bdef}, \eqref{M1def} and \eqref{M2def} by replacing $(h^\varepsilon,u^\varepsilon)$ with $(h,u)$.
\end{lem}
\proof
From $\int_\mathds{R} \mathscr{E}\, \mathrm{d}x \leqslant E$ we have $\|(\mathscr{C}, P^2, Q^2)\|_{L^1} \leqslant C$, then the proof of \eqref{Estimate1} follows from \eqref{RR}, \eqref{estimate1}, \eqref{estimate3} and \eqref{estimate7}.
Since $\chi_\varepsilon(\lambda) \leqslant \lambda^2$ we obtain \eqref{Estimate2}. Then, \eqref{Estimate2} with \eqref{estimate3} imply \eqref{Estimate3}. 
In remains to prove \eqref{Estimate4}. For that purpose, we use Young inequality, \eqref{Adef} and \eqref{M2def} to obtain
\begin{equation}\label{Ne}
\left\|\mathscr{A}\right\|_{L^\infty}\ +\ \left\| \mathscr{A}_x \right\|_{L^2}\ +\ \left\| \mathscr{A}_x\right\|_{L^\infty}\ +\ \left\| \mathscr{V}_2\right\|_{L^\infty}\ \leqslant\ C.
\end{equation}
The estimates \eqref{estimate1}, \eqref{estimate3}, \eqref{Estimate2}, \eqref{M1def}, \eqref{Ne}, \eqref{Bdef} and the Cauchy--Schwarz inequality imply \eqref{Estimate4}.
\qed

\medskip

\textit{Proof of Theorem \ref{thm:ex:ep}.} 
Following \cite{guelmame2020Euler,
guelmame2020rSV,Israwi2011,Lannes,
liu2019well}, we can prove easily the local existence of solutions of \eqref{Appsys}. Integrating the energy equation \eqref{energyequationep} on $[0,t] \times \mathds{R} $, we obtain \eqref{energyconservationep}.

\textbf{Step 1.}
Defining 
\begin{gather*}
U^\varepsilon\ \eqdef\ (h^\varepsilon\, -\, \bar{h},u^\varepsilon)^\top \qquad
A(U^\varepsilon)\ \eqdef\ \begin{pmatrix} 3\, \gamma\, (h^\varepsilon)^{-3} & 0 \\ 0 & h^\varepsilon \end{pmatrix}, \qquad   
B(U^\varepsilon)\ \eqdef\ \begin{pmatrix} u^\varepsilon & h^\varepsilon \\ 3\, \gamma\, (h^\varepsilon)^{-3} & u^\varepsilon \end{pmatrix},
\\
\mathscr{F}^\varepsilon(U^\varepsilon)\ \eqdef\ \begin{pmatrix} \mathscr{A}^\varepsilon_x \\ - \mathcal{L}_{h^\varepsilon}^{-1}  \partial_x\, \left\{ \mathscr{C}^\varepsilon\, +\, F(h^\varepsilon) \right\}\, +\ \mathscr{B}^\varepsilon \end{pmatrix},
\end{gather*}
the system \eqref{Appsys} becomes
\begin{equation}\label{system}
U^\varepsilon_t\ +\ B(U^\varepsilon)\, U^\varepsilon_x\ =\ \mathscr{F}^\varepsilon(U^\varepsilon).
\end{equation}
Let $(\cdot,\cdot)$ be the scalar product in $L^2$ and 
$ E(U^\varepsilon)\ \eqdef\ \left(\Lambda^3\, U^\varepsilon,\ A^\varepsilon\, \Lambda^3\, U^\varepsilon \right)$.
Since $A^\varepsilon B^\varepsilon$ is a symmetric matrix, straightforward calculations with \eqref{system} imply that
\begin{align}\label{energyterms}
E(U^\varepsilon)_t\ 
=&\ -2\, \left( [\Lambda^3, B^\varepsilon]\, U^\varepsilon_x,\, A^\varepsilon\, \Lambda^3 U^\varepsilon \right)\ -\ 2\, \left(B^\varepsilon\, \Lambda^3  U^\varepsilon_x,\, A^\varepsilon\, \Lambda^3 U^\varepsilon \right) \nonumber \\
&\ + 2\, \left( \Lambda^3 \mathscr{F}^\varepsilon,\, A^\varepsilon\, \Lambda^3 U^\varepsilon \right)\ +\ \left( \Lambda^3 U^\varepsilon,\, A^\varepsilon_t\, \Lambda^3 U^\varepsilon \right) \nonumber \\
=&\ -2\, \left( [\Lambda^3, B^\varepsilon]\, U^\varepsilon_x,\, A^\varepsilon\, \Lambda^3 U^\varepsilon \right)\ +\ \left( \Lambda^3  U^\varepsilon,\, (A^\varepsilon\, B^\varepsilon)_x\, \Lambda^3 U^\varepsilon \right) \nonumber \\
&\ + 2\, \left( \Lambda^3 \mathscr{F}^\varepsilon,\, A^\varepsilon\, \Lambda^3 U^\varepsilon \right)\ +\ \left( \Lambda^3 U^\varepsilon,\, A^\varepsilon_t\, \Lambda^3 U^\varepsilon \right).
\end{align}
From the definition of $\chi_\varepsilon$ we have 
\begin{equation*}
\left| \chi_\varepsilon(\xi) \right|\ \leqslant\ \xi^2, \qquad \left| \chi_\varepsilon'(\xi) \right|\ \leqslant\ 2 |\xi|, \qquad \left| \chi_\varepsilon''(\xi) \right|\ \leqslant\ 2.
\end{equation*}
Then, using the Gagliardo--Nirenberg interpolation inequality $\|f_x\|_{L^4}^2 \leqslant C \|f\|_{L^\infty} \|f_{xx}\|_{L^2}$ with \eqref{Algebra}, we obtain
\begin{align}\nonumber
\left\| \chi_\varepsilon(P^\varepsilon) \right\|_{H^2}\ 
&\leqslant\ C \left[ 
\left\| \chi_\varepsilon(P^\varepsilon)  \right\|_{L^2}\ +\ 
\left\| \chi_\varepsilon'(P^\varepsilon)\, P^\varepsilon_x  \right\|_{L^2}\ +\ 
\left\| \chi_\varepsilon'(P^\varepsilon)\, P^\varepsilon_{xx} \right\|_{L^2}\ +\ 
\left\| \chi_\varepsilon''(P^\varepsilon)\, (P_x^\varepsilon)^2 \right\|_{L^2} \right] \\ \label{chiPH}
&\leqslant\ C\, \|P^\varepsilon\|_{L^\infty}\, \left\| P^\varepsilon \right\|_{H^2}\ \leqslant\ C\, \|U^\varepsilon_x\|_{L^\infty}\, \left\| U^\varepsilon \right\|_{H^3}.
\end{align}
The same inequality can be obtained for $Q^\varepsilon$
\begin{equation}\label{chiQH}
\left\| \chi_\varepsilon(Q^\varepsilon) \right\|_{H^2}\ 
\leqslant\ C\, \|U^\varepsilon_x\|_{L^\infty}\, \left\| U^\varepsilon \right\|_{H^3}.
\end{equation}
Using \eqref{estimate5.5} and \eqref{Algebra}, we obtain 
\begin{align}\nonumber
\left\| \left(\mathscr{A}^\varepsilon, \mathscr{A}^\varepsilon_x \right) \right\|_{H^3}\ 
&\leqslant\ C \left[ 
\|P^\varepsilon\|_{L^\infty}^2\, \|h^\varepsilon\, -\, \bar{h}\|_{H^2}\ +\ \left\| \left(\chi_\varepsilon(P^\varepsilon), \chi_\varepsilon(Q^\varepsilon) \right) \right\|_{H^2}
\right]\\ \label{AAx}
&\leqslant\ C \left(1\, +\, \|U^\varepsilon_x\|^2_{L^\infty} \right) \left\| U^\varepsilon \right\|_{H^3}.
\end{align}
Using \eqref{energyconservationep}, \eqref{estimate3}, \eqref{estimate3.5}, \eqref{estimate6} and Lemma \ref{lem:bounded} we obtain that 
\begin{equation}\label{F2W}
\left\| - \mathcal{L}_{h^\varepsilon}^{-1}  \partial_x\, \left\{ \mathscr{C}^\varepsilon\, +\, F(h^\varepsilon) \right\}\, +\ \mathscr{B}^\varepsilon  \right\|_{W^{1,\infty}}\ \leqslant\ C \left(1\, +\, \|U^\varepsilon_x\|^2_{L^\infty} \right).
\end{equation}
Using now \eqref{estimate5}, \eqref{Algebra}, \eqref{Composition2}, \eqref{chiPH}, \eqref{chiQH}, \eqref{AAx} and \eqref{F2W} we obtain 
\begin{equation}
\left\|\mathscr{B}^\varepsilon  - \mathcal{L}_{h^\varepsilon}^{-1}  \partial_x\, \left\{ \mathscr{C}^\varepsilon\, +\, F(h^\varepsilon) \right\}  \right\|_{H^3} \leqslant\ \mathcal{P} \left(\|U_x^\varepsilon\|_{L^\infty} \right) \|U^\varepsilon\|_{H^3},
\end{equation}
where $\mathcal{P}$ is a polynomial function.
The last inequality with \eqref{AAx} imply that 
\begin{equation}\label{I}
\|\mathscr{F}^\varepsilon\|_{H^3}\ \leqslant\  \mathcal{P} \left(\|U^\varepsilon_x\|_{L^\infty} \right) \|U^\varepsilon\|_{H^3}.
\end{equation}
Defining $\bar{B} \eqdef B(\bar{h},0)$, and using \eqref{Commutator} one obtains
\begin{align} \nonumber
\left| \left( [\Lambda^3, B^\varepsilon]\, U^\varepsilon_x, A^\varepsilon\, \Lambda^3 U^\varepsilon \right)\ \right|\ 
&\leqslant\ C\, \|A^\varepsilon\|_{L^\infty} \|U^\varepsilon\|_{H^3} \left( \|B^\varepsilon_x\|_{L^\infty} \|U^\varepsilon_x\|_{H^{2}}\, +\, \|B^\varepsilon-\bar{B}\|_{H^3} \|U^\varepsilon_x\|_{L^\infty}\!  \right) \\ \label{II}
&\leqslant\ C\, \|U^\varepsilon_x\|_{L^\infty}\, \|U^\varepsilon\|_{H^3}^2.
\end{align}
Using \eqref{Appsys1} and \eqref{Estimate4} one obtains that 
\begin{equation}\label{III}
\left| \left( \Lambda^3  U^\varepsilon,\, (A^\varepsilon\, B^\varepsilon)_x\, \Lambda^3 U^\varepsilon \right)\right|\ +\ \left|   \left( \Lambda^3 U^\varepsilon,\, A^\varepsilon_t\, \Lambda^3 U^\varepsilon \right)  \right|\
\leqslant\  C\, \left( \|U^\varepsilon_x\|_{L^\infty}\, +\, 1 \right) \|U^\varepsilon\|_{H^3}^2.
\end{equation}
Summing up \eqref{I}, \eqref{II} and \eqref{III} we obtain
\begin{equation*}
E(U^\varepsilon)_t\ \leqslant\  \mathcal{P} \left(\|U^\varepsilon_x\|_{L^\infty} \right) \|U^\varepsilon\|_{H^3}\ \leqslant\ \mathcal{P} \left(\|U^\varepsilon_x\|_{L^\infty} \right) E(U^\varepsilon),
\end{equation*}
which implies with Gronwall inequality that 
\begin{equation}\label{H3estimate}
\|U^\varepsilon\|_{H^3}\ \leqslant\ C\, E(U^\varepsilon)\ \leqslant\ C\, E(U^\varepsilon_0)\, \mathrm{e}^{\int_0^t \mathcal{P} \left(\|U^\varepsilon_x\|_{L^\infty} \right) \mathrm{d}s}\ \leqslant\ C\, \|U^\varepsilon_0\|_{H^3}\, \mathrm{e}^{\int_0^t \mathcal{P} \left(\|U^\varepsilon_x\|_{L^\infty} \right) \mathrm{d}s}.
\end{equation}
This implies that if $T^\varepsilon_{max}$ is the maximal existence time, then
\begin{equation}\label{BUCep}
T_{max}^\varepsilon\, <\, \infty \ \implies 
\limsup_{t\to T^\varepsilon_{max}}\, \|U^\varepsilon_x(t,\cdot)\|_{L^\infty}\ =\ \infty. 
\end{equation}

\textbf{Step 2.} Defining 
\begin{align*}
\mathscr{H}^\varepsilon_1\ 
\eqdef &\ \half\, \sqrt{3\, \gamma} \left( (h^\varepsilon)^{1/2}\, (u^\varepsilon)^2\ +\ g\, (h^\varepsilon)^{-1/2} \left(h^\varepsilon\, -\, \bar{h}\right)^2\right)\, -\ u^\varepsilon \left(\mathscr{R}^\varepsilon\ +\ \half\, g\, ((h^\varepsilon)^2-\bar{h}^2)\right)\\  
&+\   u^\varepsilon \, \third\, (h^\varepsilon)^2\, \mathscr{V}^\varepsilon_1\
+\ {\textstyle \frac{\sqrt{3\, \gamma}}{3}}\, (h^\varepsilon)^{1/2}\, (h-\bar{h})\, \mathscr{V}_2^\varepsilon,\\
\mathscr{H}_2^\varepsilon\ 
\eqdef &\  \half\, \sqrt{3\, \gamma} \left( (h^\varepsilon)^{1/2}\, (u^\varepsilon)^2\ +\ g\, (h^\varepsilon)^{-1/2} \left(h^\varepsilon\, -\, \bar{h}\right)^2\right)\, +\ u^\varepsilon \left(\mathscr{R}^\varepsilon\ +\ \half\, g\, ((h^\varepsilon)^2-\bar{h}^2)\right)\\  
&-\   u^\varepsilon \, \third\, (h^\varepsilon)^2\, \mathscr{V}^\varepsilon_1\
-\ {\textstyle \frac{\sqrt{3\, \gamma}}{3}}\, (h^\varepsilon)^{1/2}\, (h-\bar{h})\, \mathscr{V}_2^\varepsilon.
\end{align*} 
We note that 
\begin{align*}
\eta^\varepsilon\, \mathscr{E}^\varepsilon\ -\ \tilde{\mathscr{D}}^\varepsilon\ =\ {\textstyle \frac{\sqrt{3\, \gamma}}{6}}\, (h^\varepsilon)^{1/2}\, (P^\varepsilon)^2\ +\ \mathscr{H}^\varepsilon_1, \qquad \quad
\tilde{\mathscr{D}}^\varepsilon\ -\ \lambda^\varepsilon\, \mathscr{E}^\varepsilon\ =\ {\textstyle \frac{\sqrt{3\, \gamma}}{6}}\, (h^\varepsilon)^{1/2}\, (Q^\varepsilon)^2\ +\ \mathscr{H}^\varepsilon_2.
\end{align*}
From Lemma \ref{lem:bounded} we deduce that $\mathscr{H}^\varepsilon_1$ and $\mathscr{H}^\varepsilon_2$ are bounded, then, integrating \eqref{energyequationep} on the set (see Figure \ref{XY}) 
$$\left\{ (s,x),\ s \in [\tau,t], X_{x_1}(s) \leqslant x \leqslant Y_{x_2}(s)  \right\},$$ 
and using the divergence theorem with \eqref{energyconservationep} one obtains \eqref{PQ^2} for all $t \in [0,T^\varepsilon_{max})$. 

Defining $t_1 \eqdef \inf \{t \geqslant 0, P^\varepsilon(t,Y_{x_2}(t)) \geqslant 1 \}$ and $t_2 \leqslant T^\varepsilon_{max}$ be the largest time such that $P^\varepsilon(t,Y_{x_2}(t)) \geqslant 1 $ on  $[t_1,t_2]$.
Dividing \eqref{Ricatti1_ep} by $P^\varepsilon$ and integrating on the characteristics between $t_1$ and $t \in [t_1,t_2]$ we obtain with \eqref{PQ^2} and Lemma \ref{lem:bounded} that
\begin{equation*}
P^\varepsilon(t,Y_{x_2}(t))\ \leqslant\ P^\varepsilon(t_1,Y_{x_2}(t_1))\, \mathrm{e}^{ C \left( 1\, +\, t\right)t} \qquad \forall t \in [t_1,t_2].
\end{equation*}
Using that $P^\varepsilon(t_1,Y_{x_2}(t_1)) = \max \left\{1, P^\varepsilon_0(x_2) \right\}$ and doing the same for $Q^\varepsilon$ we obtain 
\begin{align}\label{Psupbound}
P^\varepsilon(t,Y_{x_2}(t))\ \leqslant\ \max \left\{1, P^\varepsilon_0(x_2) \right\}\, \mathrm{e}^{ C \left( 1\, +\, t\right)t} \qquad \forall (t,x_2) \in [0,T^\varepsilon_{max}) \times \mathds{R},\\ \label{Qsupbound}
Q^\varepsilon(t,X_{x_1}(t))\ \leqslant\ \max \left\{1, Q^\varepsilon_0(x_1) \right\}\, \mathrm{e}^{ C \left( 1\, +\, t\right)t} \qquad \forall (t,x_1) \in [0,T^\varepsilon_{max})\times \mathds{R}.
\end{align}
On another hand, we define $\tilde{t}_1 \eqdef \inf \{t \geqslant 0, P^\varepsilon(t,Y_{x_2}(t)) \leqslant -1/\varepsilon \}$ and $\tilde{t}_2 \leqslant T^\varepsilon_{max}$ be the largest time such that $P^\varepsilon(t,Y_{x_2}(t)) \leqslant -1/\varepsilon $ on  $[\tilde{t}_1,\tilde{t}_2]$.
Using \eqref{Ricatti1_ep} and Lemma \ref{lem:bounded} one obtains 
\begin{equation}
\frac{\mathrm{d}^\lambda}{\mathrm{d}t}\, P^\varepsilon\ \eqdef\ P^\varepsilon_t\ +\ \lambda^\varepsilon\, P^\varepsilon_x\ \geqslant\ C \left( {\textstyle \frac{1}{\varepsilon}\, +\, 1} \right)P^\varepsilon\ -\ C \qquad \forall t \in \left[\tilde{t}_1,\tilde{t}_2\right].
\end{equation}
Using that $P^\varepsilon(\tilde{t}_1,Y_{x_2}(\tilde{t}_1)) = \min \left\{P^\varepsilon_0(x_2),-1/\varepsilon \right\}$ we obtain for all $(t,x_2) \in [0,T^\varepsilon_{max}) \times \mathds{R}$
\begin{equation}\label{Pinfbound}
P^\varepsilon(t,Y_{x_2}(t))\, \geqslant\,
 \min\! \left\{-1/\varepsilon,  \min\! \left\{P^\varepsilon_0(x_2),-1/\varepsilon \right\} \mathrm{e}^{C \left( 1 + 1/\varepsilon \right) t}\, +\, {\textstyle \frac{ \varepsilon}{\varepsilon + 1}}  \left( 1 - \mathrm{e}^{C \left( 1 + 1/\varepsilon \right) t} \right)
 \right\}.
 \end{equation}
Doing the same for $Q^\varepsilon$ we obtain for all $(t,x_1) \in [0,T^\varepsilon_{max}) \times \mathds{R}$
\begin{equation}\label{Qinfbound}
Q^\varepsilon(t,X_{x_1}(t))\, \geqslant\,
 \min\! \left\{-1/\varepsilon,  \min\! \left\{Q^\varepsilon_0(x_1),-1/\varepsilon \right\} \mathrm{e}^{C \left( 1 + 1/\varepsilon \right) t}\, +\, {\textstyle \frac{ \varepsilon}{\varepsilon + 1}}  \left( 1 - \mathrm{e}^{C \left( 1 + 1/\varepsilon \right) t} \right)
 \right\}.
\end{equation}
Finally, using \eqref{BUCep}, \eqref{Psupbound}, \eqref{Qsupbound}, \eqref{Pinfbound} and \eqref{Qinfbound} we deduce that $T^\varepsilon_{max}=  \infty$.
\qed 

The remaining of this section is devoted to obtain some uniform (on $\varepsilon$) estimates on the solution of \eqref{Appsys} given by Theorem \ref{thm:ex:ep}. Those estimate are crucial to obtain the precompactness results in next section.
\begin{lem}\label{lem:bounded2}
Let $(h_0-\bar{h},u_0) \in \mathfrak{D}$ and let $(h^\varepsilon-\bar{h},u^\varepsilon)$ be the solution given by Theorem \ref{thm:ex:ep}, then there exists a constant $C=C(\gamma, \bar{h}, E)>0$ independent on $\varepsilon \leqslant \varepsilon_0$ and  $(h_0-\bar{h},u_0)$ such that 
\begin{gather} \label{Estimate6}
\left\| \mathcal{L}_{h^\varepsilon}^{-1}  \partial_x\, \mathscr{C}^\varepsilon  \right\|_{L^\infty(\mathds{R}^+ \times \mathds{R})} +
\left\| \left( \mathscr{B}^\varepsilon, \mathscr{V}^\varepsilon_1, \mathscr{V}^\varepsilon_2, \mathscr{R}^\varepsilon \right) \right\|_{L^\infty(\mathds{R}^+ \times \mathds{R})} + \left\| \mathscr{A}^\varepsilon \right\|_{L^\infty(\mathds{R}^+, H^1(\mathds{R}))}
\leqslant\, C, \\ \label{Estimate7}
\int_{\mathds{R}^+} \int_\mathds{R} \left( \chi_\varepsilon (P^\varepsilon)\ +\ \chi_\varepsilon (Q^\varepsilon) \right)\, \mathrm{d}x\, \mathrm{d}t\
\leqslant\ \varepsilon\, C,\\ \label{Estimate8}
\int_{\mathds{R}^+} \left\| \mathcal{L}_{h^\varepsilon}^{-1}  \partial_x\, \left\{   h^\varepsilon \left(\chi_\varepsilon(P^\varepsilon)\, +\, \chi_\varepsilon(Q^\varepsilon) \right)  \right\}  \right\|_{L^\infty\left( \mathds{R}\right)}\, \mathrm{d}t\ 
\leqslant\  \varepsilon\, C,\\ \label{Estimate8.5}
\int_{\mathds{R}^+} \left\| \left( \mathcal{L}_{h^\varepsilon}^{-1}  \partial_x\, \left\{ (h^\varepsilon)^2\, u^\varepsilon_x\, \mathscr{A}^\varepsilon_x \right\}, \mathcal{L}_{h^\varepsilon}^{-1} \left\{u^\varepsilon\, \mathscr{A}^\varepsilon_x \right\} \right) \right\|_{L^\infty\left( \mathds{R}^\varepsilon\right)}\, \mathrm{d}t\ \leqslant\  \varepsilon\, C\\ \label{Estimate9}
\int_{\mathds{R}^+} \left\| \left(\mathscr{A}^\varepsilon, \mathscr{A}^\varepsilon_x, \mathscr{B}^\varepsilon, \mathscr{V}^\varepsilon_1, \mathscr{V}^\varepsilon_2 \right) \right\|_{L^\infty\left(\mathds{R}\right)}\, \mathrm{d}t\
\leqslant\  \varepsilon\, C.
\end{gather}
\end{lem}
\proof
The inequality \eqref{Estimate6} follows from \eqref{energyconservationep},  \eqref{Estimate1} and \eqref{Estimate4}. Note that
\begin{align*}
\int_{\mathds{R}^+} \int_\mathds{R} \left( \chi_\varepsilon (P^\varepsilon)\, +\, \chi_\varepsilon (Q^\varepsilon) \right)\, \mathrm{d}x\, \mathrm{d}t\ 
\leqslant &\ - \varepsilon \int_{\{ P^\varepsilon \leqslant -1/\varepsilon \}}\! P^\varepsilon\, \chi_\varepsilon (P^\varepsilon)\, \mathrm{d}x\, \mathrm{d}t\\
&\, -\ \varepsilon \int_{\{ Q^\varepsilon \leqslant -1/\varepsilon \}}\! Q^\varepsilon\, \chi_\varepsilon (Q^\varepsilon)\, \mathrm{d}x\, \mathrm{d}t.
\end{align*}
The last inequality with \eqref{energyconservationep} imply \eqref{Estimate7}.
Finally, we use \eqref{Estimate7} and Lemma \ref{Inverseesitimates} as in the proofs of \eqref{Estimate3}, \eqref{Estimate3.5} and \eqref{Estimate4}. We integrate on $\mathds{R}^+$ with respect to $t$ to obtain \eqref{Estimate8}, \eqref{Estimate8.5} and \eqref{Estimate9}.
\qed

\begin{lem}{\bf [Oleinik inequality]}\label{lem:Oleinik}
There exists $C>0$ that depends only on $\gamma, \bar{h}, g$ and $E$ such that for all $(t,x) \in (0,\infty) \times \mathds{R}$ and $\varepsilon \leqslant \varepsilon_0$ we have 
\begin{equation}
P^\varepsilon(t,x)\ \leqslant\ C\, \left(1\, +\, t^{-1} \right),
   \qquad 
Q^\varepsilon(t,x)\ \leqslant\ C\, \left(1\, +\, t^{-1} \right).
\end{equation}
\end{lem}
\proof
Let $D>0$ be a constant such that $2D^{-1} \leqslant 16h^\varepsilon \leqslant D$, and $A,B>0$ be the constants given in Theorem \ref{thm:ex:ep}. 
Using Lemma \ref{lem:bounded2}, we obtain  a constant $M>0$ depending only on $\gamma,\bar{h}$ and $E$ such that  
$$M\ \geqslant\ \sup_{t,x}\, \left\{ {\textstyle \frac{1}{h^\varepsilon}}\, \left(\mathscr{A}_x^\varepsilon\right)^2\ +\ \mathscr{M}^\varepsilon\right\}\ +\ D\, A.$$ 
Defining 
\begin{equation}
\mathcal{F}(s)\ \eqdef\ \frac{D}{s}\ +\ \sqrt{2\, M\, D}, \qquad \mathcal{G}(s)\ \eqdef\ \mathcal{F}(s)\ +\ B\, D.
\end{equation}
The goal is to prove that for all $t$ and $x$ we have $P^\varepsilon(t,X_x(t)) \leqslant \mathcal{G}(t)$ and $Q^\varepsilon(t,Y_x(t)) \leqslant \mathcal{G}(t)$. Since the proof is the same, we only prove the inequality for $P^\varepsilon$. 

Using the inequality $- \mathscr{A}^\varepsilon_x P^\varepsilon \leqslant 2 (\mathscr{A}^\varepsilon)_x^2 + (P^\varepsilon)^2/8$ and using \eqref{Ricatti1_ep}, we obtain 
\begin{align}\nonumber
\frac{\mathrm{d}^\lambda}{\mathrm{d}t}\, P^\varepsilon\ \eqdef\ P^\varepsilon_t\ +\ \lambda^\varepsilon\, P^\varepsilon_x\ 
&\leqslant\ - {\textstyle \frac{1}{16\, h^\varepsilon}}\, (P^\varepsilon)^2\ +\ {\textstyle \frac{1}{8\, h^\varepsilon}}\, \chi_\varepsilon(P^\varepsilon)\ +\ {\textstyle \frac{1}{8\, h^\varepsilon}}\, (Q^\varepsilon)^2\ +\ {\textstyle \frac{1}{h^\varepsilon}}\, \left(\mathscr{A}_x^\varepsilon \right)^2\ +\ \mathscr{M}^\varepsilon\\ \label{Ricatti1_epOl}
&\leqslant\ - {\textstyle \frac{1}{D}}\, (P^\varepsilon)^2\  +\ D\, (Q^\varepsilon)^2\ +\ {\textstyle \frac{1}{8\, h^\varepsilon}}\, \chi_\varepsilon(P^\varepsilon)\ +\ M\ -\ A\, D.
\end{align}
Let $x \in \mathds{R}$ be fixed, we suppose that there exist $t_1>0$ such that $P^\varepsilon(t_1,X_x(t_1) = \mathcal{F}(t_1)$ and $P^\varepsilon(t,X_x(t) \geqslant \mathcal{F}(t)$ for all $t \in [t_1,t_2]$. Since $P^\varepsilon \geqslant 0$ then $\chi_\varepsilon(P^\varepsilon)=0$.
Integrating \eqref{Ricatti1_epOl} on the characteristics between $t_1$ and $t\in [t_1,t_2]$ we obtain
\begin{align} \nonumber
P^\varepsilon(t,X_x(t)\ 
&\leqslant\ P^\varepsilon(t_1,X_x(t_1)\ - \int_{t_1}^t \frac{\mathcal{F}(s)^2}{D}\, \mathrm{d}s\, +\, A\, D\, (t - t_1)\ +\ B\, D\, +\, (M-A\, D)\, (t - t_1)\\ \nonumber
&=\ \mathcal{F}(t_1)\ +\ \frac{D}{t}\ -\ \frac{D}{t_1}\ -\  M\, (t - t_1)\ -2\, \sqrt{2\, M\, D}\, \ln (t/t_1)\ +\, B\, D\\ \label{Oleinikproof}
& \leqslant\ \mathcal{G}(t).
\end{align}
Since the solution $(h^\varepsilon - \bar{h}, u^\varepsilon) \in L^\infty(\mathds{R}^+,H^3)$, then initially we have $P^\varepsilon(0^+,X_x(0^+) < \mathcal{F}(0^+)=\infty$. The inequality \eqref{Oleinikproof} shows that if $P^\varepsilon$ crosses $\mathcal{F}$ at some $t_1>0$, $P^\varepsilon$ remains always smaller than $\mathcal{G}$ for $t \geqslant t_1$.
This completes the proof of $P^\varepsilon(t,X_x(t)) \leqslant \mathcal{G}(t)$ for all $t>0$. The proof for $Q^\varepsilon$ can be done similarly. \qed

\begin{lem}{\bf [$L^{2+\alpha}$ estimates]}\label{lem:alpha+2}
For any bounded set $\Omega = [t_1,t_2] \times [a,b] \subset (0,\infty) \times \mathds{R}$ and $\alpha \in (0,1)$ there exists $C_{\alpha,\Omega}>0$ such that for all $\varepsilon \leqslant \varepsilon_0$ we have 
\begin{align}\label{alpha+2}
\int_\Omega \left[ |h^\varepsilon_t|^{2+\alpha}\ +\ |h^\varepsilon_x|^{2+\alpha}\ +\ |u^\varepsilon_t|^{2+\alpha}\ +\ |u^\varepsilon_x|^{2+\alpha} \right] \mathrm{d}x\, \mathrm{d}t\ &\leqslant\ C_{\alpha, \Omega}, \\ \label{Balpha+2} 
\left\| \mathcal{L}_{h^\varepsilon}^{-1} \left(h^\varepsilon \int_{-\infty}^x (h^\varepsilon)^{-3}\, \mathscr{C}^\varepsilon\, \mathrm{d}y\, +\, \third\, F(h^\varepsilon)_x \right) \right\|_{L^\infty \left([t_1,t_2],W^{2,2+\alpha}([a,b])\right)}\ &\leqslant\ C_{\alpha, \Omega}.
\end{align}
\end{lem}
\begin{rem}
The constant $C_{\alpha,\Omega}$ depends also on $\bar{h}, \gamma$ and $E$ but not on $\varepsilon$ and the initial data.
\end{rem}

\proof 
\textbf{Step 1.} In order to prove \eqref{alpha+2} we use the change of variables 
\begin{equation*}
\tau\ \eqdef\ t, \qquad 
z\ \eqdef\ \half\, \left( \int_{-\infty}^x - \int_x^{\infty} \right) \left( h^\varepsilon(t,y)\ -\ \bar{h} \right)\,\ud\/y\ +\bar{h}\, x,
\end{equation*}
we obtain with \eqref{Appsys1} that
\begin{gather*}
\partial_x\ =\ h^\varepsilon\,\partial_z, \qquad
\partial_t\ =\ \partial_\tau\ +\, \left( \mathscr{A}^\varepsilon -\ h^\varepsilon\,u^\varepsilon \right)\partial_z, \qquad
\partial_t\ +\, u^\varepsilon\, \partial_x\ =\ \partial_\tau\ +\ \mathscr{A}^\varepsilon\, \partial_z.
\end{gather*}
The map 
$$\Phi : \mathds{R}^+ \times \mathds{R} \to \mathds{R}^+ \times \mathds{R}, \qquad (t,x)\ \mapsto\ \Phi(t,x)=(\tau,z) $$
is bijective.
Then \eqref{Ricatti_ep} becomes 
\begin{subequations}\label{Ricatti_epxi}
\begin{align}\nonumber
P^\varepsilon_\tau\, +\, \left(\mathscr{A}^\varepsilon\, -\, \sqrt{3\, \gamma}\, (h^\varepsilon)^{1/2}\right) P^\varepsilon_z\, 
=&\ - {\textstyle \frac{1}{8\, h^\varepsilon}} (P^\varepsilon)^2\, +\, {\textstyle \frac{1}{8\, h^\varepsilon}}\, \chi_\varepsilon(P^\varepsilon)\, +\, {\textstyle \frac{1}{8\, h^\varepsilon}}\, (Q^\varepsilon)^2\\ \label{Ricatti1_epxi}
&\, -\, {\textstyle \frac{1}{2\, h^\varepsilon}}\, \mathscr{A}^\varepsilon_x\,  P^\varepsilon\, +\, \mathscr{M}^\varepsilon, \\ \nonumber
Q^\varepsilon_\tau\ +\ \left(\mathscr{A}^\varepsilon\, +\, \sqrt{3\, \gamma}\, (h^\varepsilon)^{1/2}\right) Q^\varepsilon_z\, 
=&\ - {\textstyle \frac{1}{8\, h^\varepsilon}}\, (Q^\varepsilon)^2\, +\, {\textstyle \frac{1}{8\, h^\varepsilon}}\, \chi_\varepsilon(Q^\varepsilon)\, +\, {\textstyle \frac{1}{8\, h^\varepsilon}}\, (P^\varepsilon)^2\\ \label{Ricatti2_epxi}
&\,  -\, {\textstyle \frac{1}{2\, h^\varepsilon}}\, \mathscr{A}^\varepsilon_x\,  Q^\varepsilon\, +\, \mathscr{N}^\varepsilon.
\end{align}
\end{subequations}
Without loss of generality, we suppose that $\alpha= 2k /(2k+1)$ with $k \in \mathds{N}$.  
Multiplying \eqref{Ricatti1_epxi} by $(P^\varepsilon)^\alpha$, using $\frac{2}{\alpha + 1} = 1 + \frac{1-\alpha}{\alpha+1}$, $(P^\varepsilon)^\alpha \chi_\varepsilon(P^\varepsilon) \geqslant 0 $ and \eqref{Estimate6} one obtains
\begin{gather*}
{\textstyle \frac{1}{8\, h^\varepsilon}}\! \left\{\! {\textstyle \frac{1\, -\, \alpha}{\alpha\, +\, 1}}\, (P^\varepsilon)^{\alpha + 1}\, (P^\varepsilon-Q^\varepsilon) - (P^\varepsilon)^\alpha\, Q^\varepsilon\, (P^\varepsilon-Q^\varepsilon)\right\}\\
\leqslant\, 
\left( {\textstyle \frac{(P^\varepsilon)^{\alpha+1}}{\alpha+1}} \right)_\tau\!\! +\! \left( {\textstyle\frac{\mathscr{A}^\varepsilon -\sqrt{3\, \gamma}\, (h^\varepsilon)^{1/2} }{\alpha+1}}(P^\varepsilon)^{\alpha+1} \right)_z\!\! +\ C \left( |P^\varepsilon|^{\alpha+1} + (P^\varepsilon)^\alpha \right)\!.
\end{gather*}
Doing the same with \eqref{Ricatti2_ep}, we obtain
\begin{gather*}
{\textstyle \frac{1}{8\, h^\varepsilon}}\! \left\{\! {\textstyle \frac{1\, -\, \alpha}{\alpha\, +\, 1}}\, (Q^\varepsilon)^{\alpha + 1}\, (Q^\varepsilon-P^\varepsilon) - (Q^\varepsilon)^\alpha\, P^\varepsilon\, (Q^\varepsilon-P^\varepsilon)\right\}\\
\leqslant\, 
\left( {\textstyle \frac{(Q^\varepsilon)^{\alpha+1}}{\alpha+1}} \right)_\tau\!\! +\! \left( {\textstyle\frac{\mathscr{A}^\varepsilon +\sqrt{3\, \gamma}\, (h^\varepsilon)^{1/2} }{\alpha+1}}(Q^\varepsilon)^{\alpha+1} \right)_z\!\! +\ C \left( |Q^\varepsilon|^{\alpha+1} + (Q^\varepsilon)^\alpha \right)\!.
\end{gather*}
Adding both equations yields to 
\begin{gather}\nonumber
{\textstyle \frac{1}{8\, h^\varepsilon}} \left\{\! {\textstyle \frac{1\, -\, \alpha}{\alpha\, +\, 1}}\, ((P^\varepsilon)^{\alpha + 1}-(Q^\varepsilon)^{\alpha + 1})\, (P^\varepsilon-Q^\varepsilon)\, +\, (P^\varepsilon)^\alpha (Q^\varepsilon)^\alpha ((P^\varepsilon)^{1-\alpha}\, -\,  (Q^\varepsilon)^{1-\alpha})\, (P^\varepsilon-Q^\varepsilon)\right\}\\ \nonumber 
\leqslant\, 
\left( {\textstyle \frac{(P^\varepsilon)^{\alpha+1}\, +\, (Q^\varepsilon)^{\alpha+1}}{\alpha+1}} \right)_\tau\!\!  +\! \left( {\textstyle\frac{\sqrt{3\, \gamma}\, (h^\varepsilon)^{1/2}\, ((Q^\varepsilon)^{\alpha+1}-(P^\varepsilon)^{\alpha+1})\, +\, \mathscr{A}^\varepsilon\, ((Q^\varepsilon)^{\alpha+1}+(P^\varepsilon)^{\alpha+1}) }{\alpha+1}} \right)_z\\ \label{Lalpha}
+\ C \left( |Q^\varepsilon|^{\alpha+1} + (Q^\varepsilon)^\alpha + |P^\varepsilon|^{\alpha+1} + (P^\varepsilon)^\alpha \right).
\end{gather}
Let $\varphi \in C^\infty_c((t_1/2, t_2+1) \times (a-1,b+1))$ be a non negative function such that $\varphi(t,x)=1$ on $\Omega$. Multiplying \eqref{Lalpha} by $\varphi(\Phi^{-1}(\tau,z))$ and using integration by parts with \eqref{energyconservationep} we obtain
\begin{align*}
& {\textstyle \frac{1\, -\, \alpha}{\alpha\, +\, 1}} \int_{\mathds{R}^+ \times \mathds{R}}\! \varphi\, (P^\varepsilon-Q^\varepsilon)^2((P^\varepsilon)^\alpha+(Q^\varepsilon)^\alpha)\, \mathrm{d}x\, \mathrm{d}t\\
\leqslant &\ \int_{\mathds{R}^+ \times \mathds{R}}\!\!  \left\{ {\textstyle \frac{1\, -\, \alpha}{\alpha\, +\, 1}} ((P^\varepsilon)^{\alpha + 1}-(Q^\varepsilon)^{\alpha + 1})\, (P^\varepsilon-Q^\varepsilon)\right\} \varphi(t,x)\, \mathrm{d}x\, \mathrm{d}t\\
&+\, \int_{\mathds{R}^+ \times \mathds{R}}\!\!  \left\{ (P^\varepsilon)^\alpha (Q^\varepsilon)^\alpha ((P^\varepsilon)^{1-\alpha}\, -\,  (Q^\varepsilon)^{1-\alpha})\, (P^\varepsilon-Q^\varepsilon) \right\}\varphi(t,x)\, \mathrm{d}x\, \mathrm{d}t\\
= &\  \int_{\mathds{R}^+ \times \mathds{R}}\!\!  \left\{ {\textstyle \frac{1\, -\, \alpha}{\alpha\, +\, 1}} ((P^\varepsilon)^{\alpha + 1}-(Q^\varepsilon)^{\alpha + 1})\, (P^\varepsilon-Q^\varepsilon) \right\} \varphi(\Phi^{-1}(\tau,z)) \frac{\mathrm{d}z\, \mathrm{d}\tau}{h}\\
&+\, \int_{\mathds{R}^+ \times \mathds{R}}\!\!  \left\{ (P^\varepsilon)^\alpha (Q^\varepsilon)^\alpha ((P^\varepsilon)^{1-\alpha}\, -\,  (Q^\varepsilon)^{1-\alpha})\, (P^\varepsilon-Q^\varepsilon) \right\}\varphi(\Phi^{-1}(\tau,z)) \frac{\mathrm{d}z\, \mathrm{d}\tau}{h}\\
\leqslant &\ 8 \int_{\mathds{R}^+ \times \mathds{R}}\!\! \left[
C \left( |Q^\varepsilon|^{\alpha+1} + (Q^\varepsilon)^\alpha + |P^\varepsilon|^{\alpha+1} + (P^\varepsilon)^\alpha \right) \varphi\, \mathrm{d}z\, 
- \left( {\textstyle \frac{(P^\varepsilon)^{\alpha+1}\, +\, (Q^\varepsilon)^{\alpha+1}}{\alpha+1}} \right)  \varphi_\tau \right] \mathrm{d}z\, \mathrm{d}\tau\\
&- 8\int_{\mathds{R}^+ \times \mathds{R}}\!  \left( {\textstyle\frac{\sqrt{3\, \gamma}\, (h^\varepsilon)^{1/2}\, ((Q^\varepsilon)^{\alpha+1}-(P^\varepsilon)^{\alpha+1})\, +\, \mathscr{A}^\varepsilon\, ((Q^\varepsilon)^{\alpha+1}+(P^\varepsilon)^{\alpha+1}) }{\alpha+1}} \right) \varphi_z\,
\mathrm{d}z\, \mathrm{d}\tau\\
= &\ 8 \int_{\mathds{R}^+ \times \mathds{R}}\!\! 
C \left( |Q^\varepsilon|^{\alpha+1} + (Q^\varepsilon)^\alpha + |P^\varepsilon|^{\alpha+1} + (P^\varepsilon)^\alpha \right) \varphi\, h^\varepsilon\, \mathrm{d}x\, \mathrm{d}t\\
&- 8\int_{\mathds{R}^+ \times \mathds{R}}\! 
 \left( {\textstyle \frac{(P^\varepsilon)^{\alpha+1}\, +\, (Q^\varepsilon)^{\alpha+1}}{\alpha+1}} \right) h^\varepsilon \left( \varphi_t\, +\, u^\varepsilon\, \varphi_x -\, (h^\varepsilon)^{-1}\, \mathscr{A}^\varepsilon\, \varphi_x \right)\, \mathrm{d}x\, \mathrm{d}t\\
&- 8\int_{\mathds{R}^+ \times \mathds{R}}\!  
 \left( {\textstyle\frac{\sqrt{3\, \gamma}\, (h^\varepsilon)^{1/2}\, ((Q^\varepsilon)^{\alpha+1}-(P^\varepsilon)^{\alpha+1})\, +\, \mathscr{A}^\varepsilon\, ((Q^\varepsilon)^{\alpha+1}+(P^\varepsilon)^{\alpha+1}) }{\alpha+1}} \right) \varphi_x\, \mathrm{d}x\, \mathrm{d}t\\
\leqslant &\ C_{\alpha,\Omega}.
\end{align*}
The last inequality follows from \eqref{energyconservationep} and from the fact that $\varphi$ is compactly supported. Then we have
\begin{equation}\label{PQalpha-}
 \int_{\Omega}\! (P^\varepsilon-Q^\varepsilon)^2((P^\varepsilon)^\alpha+(Q^\varepsilon)^\alpha)\, \mathrm{d}x\, \mathrm{d}t\, \leqslant\, C_{\alpha,\Omega} .
\end{equation}
\textbf{Step 2.}
Multiplying \eqref{Lalpha} by $(h^\varepsilon)^{-1/2}$ we obtain 
\begin{gather*}
{\textstyle \frac{1}{8\, (h^\varepsilon)^{3/2}}}\, \left\{\! {\textstyle \frac{1\, -\, \alpha}{\alpha\, +\, 1}}\, ((P^\varepsilon)^{\alpha + 1}-(Q^\varepsilon)^{\alpha + 1})\, (P^\varepsilon-Q^\varepsilon)\, +\, (P^\varepsilon)^\alpha (Q^\varepsilon)^\alpha ((P^\varepsilon)^{1-\alpha}\, -\,  (Q^\varepsilon)^{1-\alpha})\, (P^\varepsilon-Q^\varepsilon)\right\}\\ 
\leqslant\, \left( {\textstyle \frac{(P^\varepsilon)^{\alpha+1}\, +\, (Q^\varepsilon)^{\alpha+1}}{(\alpha+1)\, (h^\varepsilon)^{1/2}}} \right)_\tau\,  +\! \left( {\textstyle\frac{\sqrt{3\, \gamma}\, ((Q^\varepsilon)^{\alpha+1}-(P^\varepsilon)^{\alpha+1}) }{\alpha+1}} \right)_z\, +\ \left( {\textstyle\frac{ \mathscr{A}^\varepsilon\, ((Q^\varepsilon)^{\alpha+1}+(P^\varepsilon)^{\alpha+1}) }{(\alpha+1)\, (h^\varepsilon)^{1/2}}} \right)_z\\
+\ (h^\varepsilon)^{-1/2}\, C \left( |Q^\varepsilon|^{\alpha+1} + (Q^\varepsilon)^\alpha + |P^\varepsilon|^{\alpha+1} + (P^\varepsilon)^\alpha \right)\ +\ \half\, \mathscr{A}_x^\varepsilon\, {\textstyle \frac{(P^\varepsilon)^{\alpha+1}\, +\, (Q^\varepsilon)^{\alpha+1}}{(\alpha+1)\, (h^\varepsilon)^{3/2}}}\\
-\ {\textstyle \frac{1}{8\, (h^\varepsilon)^{3/2}}}\, {\textstyle \frac{4}{\alpha+1}}\, \left\{\! (P^\varepsilon)^{\alpha+1}\, Q^\varepsilon\, +\, (Q^\varepsilon)^{\alpha+1}\, P^\varepsilon \right\}.
\end{gather*}
Using \eqref{Estimate6} one obtain
\begin{gather*}
{\textstyle \frac{1}{8\, (h^\varepsilon)^{3/2}}} \left\{\! {\textstyle \frac{1\, -\, \alpha}{\alpha\, +\, 1}}\, ((P^\varepsilon)^{\alpha + 1}+(Q^\varepsilon)^{\alpha + 1})\, (P^\varepsilon+Q^\varepsilon)\, +\, (P^\varepsilon)^\alpha (Q^\varepsilon)^\alpha ((P^\varepsilon)^{1-\alpha}\, +\,  (Q^\varepsilon)^{1-\alpha})\, (P^\varepsilon+Q^\varepsilon)\right\}\\ 
\leqslant\
\left( {\textstyle \frac{(P^\varepsilon)^{\alpha+1}\, +\, (Q^\varepsilon)^{\alpha+1}}{(\alpha+1)\, (h^\varepsilon)^{1/2}}} \right)_\tau\,  +\! \left( {\textstyle\frac{\sqrt{3\, \gamma}\, ((Q^\varepsilon)^{\alpha+1}-(P^\varepsilon)^{\alpha+1}) }{\alpha+1}} \right)_z\, +\ \left( {\textstyle\frac{ \mathscr{A}^\varepsilon\, ((Q^\varepsilon)^{\alpha+1}+(P^\varepsilon)^{\alpha+1}) }{(\alpha+1)\, (h^\varepsilon)^{1/2}}} \right)_z\\
C \left( |Q^\varepsilon|^{\alpha+1} + (Q^\varepsilon)^\alpha + |P^\varepsilon|^{\alpha+1} + (P^\varepsilon)^\alpha \right).
\end{gather*}
As in the first step we obtain 
\begin{equation}\label{PQalpha+}
 \int_{\Omega}\! (P^\varepsilon+Q^\varepsilon)^2((P^\varepsilon)^\alpha+(Q^\varepsilon)^\alpha)\, \mathrm{d}x\, \mathrm{d}t\, \leqslant\, C_{\alpha,\Omega} .
\end{equation}
Summing up \eqref{PQalpha-} and \eqref{PQalpha+} one obtains
\begin{equation*}
\int_\Omega ((P^\varepsilon)^{\alpha+2}+(Q^\varepsilon)^{\alpha+2})\, \mathrm{d}x\, \mathrm{d}t\ \leqslant\ C_{\alpha,\Omega} 
\qquad \implies \qquad
\int_\Omega \left[ |u^\varepsilon_x|^{2+\alpha}\ +\ |h^\varepsilon_x|^{2+\alpha} \right]\, \mathrm{d}x\, \mathrm{d}t\ \leqslant\ C_{\alpha,\Omega}.
\end{equation*}
\textbf{Step 3.}
The inequality \eqref{alpha+2} follows directly from \eqref{Appsys} and Lemma \ref{lem:bounded2}.  Finally, using \eqref{estimate1}, \eqref{estimate1.5}, \eqref{estimate7}, \eqref{estimate8} and \eqref{alpha+2} we obtain \eqref{Balpha+2}.
\qed

\section{Precompactness}\label{Sec:Precom}
The goal of this section is to obtain a compactness of the solution.
Due to the non-linear terms of the equations, strong precompactness is needed to pass to the limit $\varepsilon \to 0$. The strong precompactness of $(h^\varepsilon)_\varepsilon$ and $(u^\varepsilon)_\varepsilon$ is easy to obtain. However, the strong precompactness of $(P^\varepsilon)_\varepsilon$ and $(Q^\varepsilon)_\varepsilon$ is more challenging.
Several lemmas in this section are inspired by \cite{Coclite1,CH_CPAM,wave0,wave1,wave2,wave3}.
Along this section, Lemma \ref{lem:fg} is used many times without mentioning it.

We start by strong precompactness of $(h^\varepsilon)_\varepsilon$ and $(u^\varepsilon)_\varepsilon$.
\begin{lem}\label{Strong_conv}
There exist $(h-\bar{h},u) \in L^\infty ([0,\infty), H^1(\mathds{R}))$ and a subsequence of $(h^\varepsilon,u^\varepsilon)_\varepsilon$ such that we have the following convergences
\begin{align*}
(h^\varepsilon- \bar{h},u^\varepsilon)\qquad &\to \qquad (h- \bar{h},u) \qquad \mathrm{in}\  L^{\infty}_{loc}([0,\infty) \times \mathds{R}), \\
(h^\varepsilon- \bar{h},u^\varepsilon)\qquad &\rightharpoonup \qquad (h- \bar{h},u) \qquad \mathrm{in}\   H^1([0,T]\times \mathds{R}),\ \forall T>0.
\end{align*}
\end{lem}
\proof 
From the energy equation \eqref{energyconservationep} we have that $(h^\varepsilon - \bar{h},u^\varepsilon)$ is uniformly bounded in $L^\infty([0,\infty), H^1(\mathds{R}))$. Then, using \eqref{Appsys}, and \eqref{Estimate6} we obtain that 
\begin{equation}\label{L2L2}
\left\|\left(h_t^\varepsilon,u_t^\varepsilon \right) \right\|_{L^2([0,T] \times \mathds{R})}\ \leqslant\ C_T.
\end{equation}
The weak convergence in $H^1([0,T] \times \mathds{R})$ follows directly. Using the inequality 
\begin{equation*}
\left\| \theta(t,\cdot)\, -\, \theta(s,\cdot) \right\|_{L^2(\mathds{R})}^2\ =\ \int_\mathds{R} \left( \int_s^t \theta_t(\tau,x)\, \mathrm{d}\tau \right)^2 \mathrm{d}x\ \leqslant\ |t-s|\, \left\|\theta_t \right\|_{L^2([0,T] \times \mathds{R})}^2,
\end{equation*}
with \eqref{L2L2} we obtain that 
\begin{equation*}
\lim_{t \to s} \left\| u^\varepsilon(t,\cdot)\, -\, u^\varepsilon (s,\cdot) \right\|_{L^2(\mathds{R})}\ +\ \lim_{t \to s} \left\| h^\varepsilon(t,\cdot)\, -\, h^\varepsilon (s,\cdot) \right\|_{L^2(\mathds{R})}\ =\ 0
\end{equation*}
uniformly on $\varepsilon$. Then, using Theorem 5 in \cite{Simon} we can deduce that up to a subsequence,
$(h^\varepsilon, u^\varepsilon)$ converges uniformly to $(h,u)$ on any compact set of $[0,\infty) \times \mathds{R} $ when $\varepsilon \to 0$.
\qed
%

Now, we establish the weak precompactness of $(P^\varepsilon)_\varepsilon$ and $(Q^\varepsilon)_\varepsilon$.
\begin{lem}\label{lemYoungm}
There exist a subsequence of $\left\{P^\varepsilon, Q^\varepsilon \right\}_\varepsilon$ denoted also $\left\{P^\varepsilon, Q^\varepsilon \right\}_\varepsilon$ and families of probability Young measures $\nu^1_{t,x}, \nu^2_{t,x}$ on $\mathds{R}$ and $\mu_{t,x}$ on $\mathds{R}^2$, such that for all functions $f,\phi \in C^\infty_c(\mathds{R})$, $g \in C(\mathds{R}^2)$ with $g(\xi,\zeta) = \mathcal{O} (|\xi|^{2}+|\zeta|^{2})$ at infinity, and for all $\varphi \in C^\infty_c((0,\infty) \times \mathds{R})$ we have
\begin{align}\label{convegrgencefa}
\lim_{\varepsilon \to 0}\, \int_\mathds{R} \phi(x)\, f(P^\varepsilon(t,x))\, \mathrm{d}x\ 
=\ \int_\mathds{R} \phi(x) \int_\mathds{R} f(\xi)\, \mathrm{d} \nu^1_{t,x}(\xi)\, \mathrm{d}x, \\ \label{convegrgencefb}
\lim_{\varepsilon \to 0}\, \int_\mathds{R} \phi(x)\, f(Q^\varepsilon(t,x))\, \mathrm{d}x\ 
=\ \int_\mathds{R} \phi(x) \int_\mathds{R}  f(\zeta)\,  \mathrm{d} \nu^2_{t,x}(\zeta)\, \mathrm{d}x, 
\end{align}
uniformly on any compact set $[0,T] \subset [0,\infty)$, and
\begin{equation}\label{convegrgenceg}
\lim_{\varepsilon \to 0}\, \int_{\mathds{R}^+ \times \mathds{R}} \varphi(t,x)\, g(P^\varepsilon,Q^\varepsilon)\, \mathrm{d}x\, \mathrm{d}t\ 
=\ \int_{\mathds{R}^+ \times \mathds{R}} \varphi(t,x) \int_{\mathds{R}^2} g(\xi,\zeta)\,  \mathrm{d} \mu_{t,x}(\xi,\zeta)\, \mathrm{d}x\, \mathrm{d}t.
\end{equation}
Moreover, the map
\begin{equation}\label{Fatou1}
(t,x)\ \mapsto\  \int_\mathds{R} \xi^2\, \mathrm{d} \nu^1_{t,x} (\xi)\ +\ \int_\mathds{R} \zeta^2\, \mathrm{d} \nu^2_{t,x} (\zeta)
\end{equation}
belongs to $L^\infty(\mathds{R}^+, L^1(\mathds{R})) $, and
\begin{equation}\label{Product}
\mu_{t,x}(\xi,\zeta)\ =\ \nu^1_{t,x}(\xi)\, \otimes \nu^2_{t,x}(\zeta).
\end{equation}
\end{lem}
We define
\begin{equation}\label{bardef}
\overline{g(P,Q)}\ \eqdef\ \int_{\mathds{R}^2} g(\xi,\zeta)\, \mathrm{d} \mu_{t,x}(\xi,\zeta)
\end{equation}
which is from \eqref{convegrgenceg} the weak limit of $g(P^\varepsilon,Q^\varepsilon)$.
\proof 
\textbf{Step 1.}
The pointwise convergence of \eqref{convegrgencefa} is a direct corollary of Lemma \ref{lem:Young} with $\mathscr{O}=\mathds{R}$ and $p=2$ and the energy equation \eqref{energyconservationep}. 
The key point to prove the uniform convergence is to show that the map
\begin{equation}\label{equicontinuousmap}
t \in [0,T]\ \mapsto\ \int_\mathds{R} \phi(x)\, f(P^\varepsilon(t,x))\,  \mathrm{d}x\, \mathrm{d}t
\end{equation}
is equicontinuous. 
Multiplying \eqref{Ricatti1_ep} with $f'(P^\varepsilon)$ one obtains
\begin{gather}\nonumber
f(P^\varepsilon)_t\ +\ \left[ \lambda^\varepsilon\, f(P^\varepsilon) \right]_x\ =\ {\textstyle \frac{1}{4\, h^\varepsilon}} \left(P^\varepsilon +3\, Q^\varepsilon \right) f(P^\varepsilon)\\ \label{rhs0}
+\
\left[ - {\textstyle \frac{1}{8\, h^\varepsilon}}\, \left(P^\varepsilon\right)^2\ +\ {\textstyle \frac{1}{8\, h^\varepsilon}}\, \chi_\varepsilon(P^\varepsilon)\ +\ {\textstyle \frac{1}{8\, h^\varepsilon}}\, \left(Q^\varepsilon\right)^2\ -\ {\textstyle \frac{1}{2\, h^\varepsilon}}\, \mathscr{A}^\varepsilon_x\,  P^\varepsilon\ +\ \mathscr{M}^\varepsilon\right] f'(P^\varepsilon).
\end{gather}
Multiplying by $\phi(x)$ and integrating over $[t_1,t_2] \times \mathds{R}$ we have
\begin{gather*}
\int_\mathds{R} \phi(x) \left[ f(P^\varepsilon(t_2,x))\ -\ f(P^\varepsilon(t_1,x))  \right]\, \mathrm{d}x \\
=\ \int_{t_1}^{t_2}  \int_\mathds{R} \left[ \phi'(x)\, \lambda^\varepsilon\, f(P^\varepsilon)\ +\ {\textstyle \frac{1}{4\, h^\varepsilon}}\, \phi(x)\, \left(P^\varepsilon+3\, Q^\varepsilon \right) f(P^\varepsilon) \right]\, \mathrm{d}x\, \mathrm{d}t \\
+
\int_{t_1}^{t_2}  \int_\mathds{R} \phi(x) \left[ - {\textstyle \frac{1}{8\, h^\varepsilon}}\, \left(P^\varepsilon\right)^2 + {\textstyle \frac{1}{8\, h^\varepsilon}}\, \chi_\varepsilon(P^\varepsilon) + {\textstyle \frac{1}{8\, h^\varepsilon}}\, \left(Q^\varepsilon\right)^2 - {\textstyle \frac{1}{2\, h^\varepsilon}}\, \mathscr{A}^\varepsilon_x\,  P^\varepsilon + \mathscr{M}^\varepsilon\right] f'(P^\varepsilon)\, \mathrm{d}x\, \mathrm{d}t .
\end{gather*}
Using that $f\in C^\infty_c$, the the energy equation \eqref{energyconservationep}, Proposition \ref{pro:ene}, and Lemma \ref{lem:bounded2} we find that the map \eqref{equicontinuousmap} is equicontinuous. This finishes the proof of the uniform convergence of \eqref{convegrgencefa}.
The same proof can be used for \eqref{convegrgencefb}.
Using \eqref{Fatou0} we deduce that the map \eqref{Fatou1} belongs to $L^\infty(\mathds{R}^+, L^1(\mathds{R}))$.

\textbf{Step 2.}
Now, we suppose that $g$ satisfies $g(\xi,\zeta) = \smallO (|\xi|^2+|\zeta|^2)$ at infinity, then, using again Lemma \ref{lem:Young} with $\mathscr{O}= (0,\infty) \times \mathds{R}$ and $p=2$ we obtain \eqref{convegrgenceg}. If $g(\xi,\zeta) = \mathcal{O} (|\xi|^2+|\zeta|^2)$, let $\psi$ be a smooth cut-off function with $\psi(\xi)=1$ for $|\xi| \leqslant 1$ and $\psi(\xi)=0$ for $|\xi| \geqslant 2$, then 
\begin{equation}\label{convegrgencegpsi}
\lim_{\varepsilon \to 0}\, \int_{\mathds{R}^+ \times \mathds{R}} \varphi(t,x)\, g_k(P^\varepsilon,Q^\varepsilon)\,  \mathrm{d}x\, \mathrm{d}t\ 
=\ \int_{\mathds{R}^+ \times \mathds{R}} \varphi(t,x) \int_{\mathds{R}^2} g_\kappa(\xi,\zeta)\,  \mathrm{d} \mu_{t,x}(\xi,\zeta)\, \mathrm{d}x\, \mathrm{d}t,
\end{equation}
where $g_\kappa(\xi,\zeta) \eqdef g(\xi,\zeta) \psi\left( {\textstyle \frac{\xi}{\kappa}}\right) \psi\left( {\textstyle \frac{\zeta}{\kappa}}\right)$ with $\kappa>0$. Using Holder inequality, Lemma \ref{lem:alpha+2} with $\Omega = \mathrm{supp}(\varphi)$ we obtain
\begin{gather*}
\left| \int_{\mathds{R}^+ \times \mathds{R}} \varphi(t,x)\, \left( g(P^\varepsilon,Q^\varepsilon)\, -\, g_\kappa(P^\varepsilon,Q^\varepsilon) \right)\, \mathrm{d}x\, \mathrm{d}t \right|\\
\leqslant\ 
\int_{\mathrm{supp}(\varphi) \cap \{|P^\varepsilon| \geqslant \kappa,\ \mathrm{or}\ |Q^\varepsilon| \geqslant \kappa\}} |\varphi(t,x)|\, \left|g(P^\varepsilon,Q^\varepsilon)\right|\, \mathrm{d}x\, \mathrm{d}t \\
\leqslant\ C \left( \int_{\mathrm{supp}(\varphi)} \left|g(P^\varepsilon,Q^\varepsilon)\right|^{p/2}\, \mathrm{d}x\, \mathrm{d}t \right)^{2/p} \left( \int_{\mathrm{supp}(\varphi) \cap \{|P^\varepsilon| \geqslant \kappa,\ \mathrm{or}\ |Q^\varepsilon| \geqslant \kappa\}} \mathrm{d}x\, \mathrm{d}t \right)^\frac{p-2}{p}
\\
\leqslant\ C\, \Big[ \big|\! \left\{(t,x) \in \mathrm{supp}(\varphi), |P^\varepsilon| \geqslant \kappa \right\}\! \big| \ +\ \big|\! \left\{(t,x) \in \mathrm{supp}(\varphi), |Q^\varepsilon| \geqslant \kappa \right\}\! \big| \Big]^\frac{p-2}{p}\\
\leqslant\ C\, \kappa^{2-p}.
\end{gather*}
where $2<p<3$.
The last inequality with \eqref{convegrgencegpsi} imply that we can interchange the limits $\kappa \to \infty$ and $\varepsilon \to 0$. Using that $|g_\kappa| \leqslant |g|$ and the dominated convergence theorem we obtain \eqref{convegrgenceg}.

\textbf{Step 3.}
It only remains to prove \eqref{Product}, for that purpose, let $f \in C^\infty_c(\mathds{R})$, we rewrite \eqref{rhs0} on the form
\begin{gather}\nonumber
f(P^\varepsilon)_t\ +\ \left[ \lambda\, f(P^\varepsilon) \right]_x\ =\ \left[ (\lambda\, -\, \lambda^\varepsilon)\, f(P^\varepsilon) \right]_x\ +\ {\textstyle \frac{1}{4\, h^\varepsilon}} \left(P^\varepsilon +3\, Q^\varepsilon \right) f(P^\varepsilon)\\ \label{rhs}
+\
\left[ - {\textstyle \frac{1}{8\, h^\varepsilon}}\, \left(P^\varepsilon\right)^2\ +\ {\textstyle \frac{1}{8\, h^\varepsilon}}\, \chi_\varepsilon(P^\varepsilon)\ +\ {\textstyle \frac{1}{8\, h^\varepsilon}}\, \left(Q^\varepsilon\right)^2\ -\ {\textstyle \frac{1}{2\, h^\varepsilon}}\, \mathscr{A}^\varepsilon_x\,  P^\varepsilon\ +\ \mathscr{M}^\varepsilon\right] f'(P^\varepsilon).
\end{gather}
Lemma \ref{Strong_conv} implies that $(\lambda - \lambda^\varepsilon) f(P^\varepsilon) \to 0$ in $L^2_{loc}((0,\infty) \times \mathds{R})$ when $\varepsilon \to 0$. This implies that $\left[ (\lambda - \lambda^\varepsilon) f(P^\varepsilon) \right]_x$ is relatively compact in $H^{-1}_{loc}((0,\infty) \times \mathds{R})$. Since $f \in C^\infty_c(\mathds{R})$, using \eqref{Estimate6} and the energy equation \eqref{energyconservationep} we obtain that the remaining terms of the right-hand side of \eqref{rhs} are uniformly bounded in $L^1_{loc}((0,\infty) \times \mathds{R})$. Then, due to Lemma \ref{lem:Evans} they are relatively compact in $H^{-1}_{loc}((0,\infty) \times \mathds{R})$. Doing the same we can prove that for all $f,g \in C^\infty_c(\mathds{R})$ the sequences 
\begin{equation*}
\left\{ \left[f(P^\varepsilon)\right]_t\ +\ \left[ \lambda\, f(P^\varepsilon) \right]_x \right\}_\varepsilon, \qquad \left\{ \left[g(Q^\varepsilon)\right]_t\ +\ \left[ \eta\, g(Q^\varepsilon) \right]_x \right\}_\varepsilon
\end{equation*}
are relatively compact in $H^{-1}_{loc}((0,\infty) \times \mathds{R})$. Then, using Lemma \ref{lem:PG} (a generalised compensated compactness result) we obtain
\begin{equation}
f(P^\varepsilon)\, g(Q^\varepsilon)\ \rightharpoonup\ \overline{f(P)}\, \overline{g(Q)} \qquad \mathrm{when} \quad \varepsilon \to 0,
\end{equation}
where $\left(\overline{f(P)}, \overline{g(Q)}\right)$ is the weak limit of $\left(f(P^\varepsilon), g(Q^\varepsilon)\right)$ defined in \eqref{bardef}. Then, for any $\varphi \in C^\infty_c((0,\infty) \times \mathds{R})$ , we have 
\begin{gather*}
\int_{\mathds{R}^+ \times \mathds{R}} \int_{\mathds{R}^2} \varphi(t,x)\, f(\xi)\, g(\zeta)\,  \mathrm{d} \mu_{t,x}(\xi,\zeta)\, \mathrm{d}x\, \mathrm{d}t\ 
=\ \lim_{\varepsilon \to 0}\, \int_{\mathds{R}^+ \times \mathds{R}} \varphi(t,x)\, f(P^\varepsilon)\, g(Q^\varepsilon)\, \mathrm{d}x\, \mathrm{d}t\\
=\ \int_{\mathds{R}^+ \times \mathds{R}} \varphi(t,x)\, \overline{f(P)}\, \overline{g(Q)}\, \mathrm{d}x\, \mathrm{d}t\\
=\ \int_{\mathds{R}^+ \times \mathds{R}} \int_{\mathds{R}^2}  \varphi(t,x)\, f(\xi)\, g(\zeta)\,  \mathrm{d} \nu^1_{t,x}(\xi)\, \otimes \nu^2_{t,x}(\zeta)\, \mathrm{d}x\, \mathrm{d}t,
\end{gather*}
which implies \eqref{Product}. The proof of Lemma \ref{lemYoungm} is completed.
\qed

Using \eqref{uh_x}, Lemma \ref{lemYoungm}, \eqref{Appsys1}, \eqref{Estimate9} and Lemma \ref{Strong_conv} we can obtain the identities
\begin{gather}\label{uh_xi}
u_x\ =\ \frac{\overline{P}\ +\ \overline{Q}}{2\, h}, \qquad h_x\ =\ h^{\frac{1}{2}}\, \frac{\overline{Q}\ -\ \overline{P}}{2\, \sqrt{3\, \gamma}}\\ \label{Appsys123}
h_t\ +\ (h\, u)_x\ =\ 0.
\end{gather}
Now, we present some technical lemmas that are needed to obtain the strong precompactness of $(P^\varepsilon)_\varepsilon$ and $(Q^\varepsilon)_\varepsilon$.
\begin{lem} As $t \to 0$ we have
\begin{gather} \label{tauto0}
\left\|(h-h_0,u-u_0) \right\|_{H^1(\mathds{R})}\ \to\ 0,\\ \label{tauto0_2}
\int_\mathds{R} \left( \overline{P^2}\ -\ \overline{P}^2\, \right) \mathrm{d}x\ +\
\int_\mathds{R} \left( \overline{Q^2}\ -\ \overline{Q}^2\, \right) \mathrm{d}x\ \to\ 0.
\end{gather}
\end{lem}
\proof
Defining
\begin{align*}
W^\varepsilon (t,x)\ &\eqdef\ \left( {\textstyle \sqrt{\frac{g}{2}}}(h^\varepsilon - \bar{h}),\,{\textstyle \sqrt{\frac{h^\varepsilon}{2}}} u^\varepsilon,\, {\textstyle \frac{1}{\sqrt{12\, h^\varepsilon}}}P^\varepsilon,\, {\textstyle \frac{1}{\sqrt{12\, h^\varepsilon}}}Q^\varepsilon\right), &t \geqslant 0, \\
W(t,x)\ &\eqdef\ \left({\textstyle \sqrt{\frac{g}{2}}}(h - \bar{h}),\, {\textstyle \sqrt{\frac{h}{2}}} u,\, {\textstyle \frac{1}{\sqrt{12\, h}}}\overline{P},\, {\textstyle \frac{1}{\sqrt{12\, h}}}\overline{Q}\right), &t > 0, \\
\widetilde{W}(t,x)\ &\eqdef\ \left({\textstyle \sqrt{\frac{g}{2}}}(h - \bar{h}),\, {\textstyle \sqrt{\frac{h}{2}}} u,\, {\textstyle \frac{1}{\sqrt{12\, h}}}\sqrt{\overline{P^2}},\, {\textstyle \frac{1}{\sqrt{12\, h}}}\sqrt{\overline{Q^2}}\right), &t > 0, \\
W_0(x)\ &\eqdef\ \left({\textstyle \sqrt{\frac{g}{2}}} (h_0 - \bar{h}),\, {\textstyle \sqrt{\frac{h_0}{2}}} u_0,\, {\textstyle \frac{1}{\sqrt{12\, h_0}}} P_0,\, {\textstyle \frac{1}{\sqrt{12\, h_0}}} Q_0\right). 
\end{align*}
From Lemma \ref{Strong_conv} and Lemma \ref{lemYoungm} we have for all $t > 0$
\begin{gather*}
W^\varepsilon(t,\cdot)\ \rightharpoonup\ W(t,\cdot) \qquad \mathrm{when}\ \varepsilon\, \to\, 0 \qquad \mathrm{in\ } L^2(\mathds{R}),\\
\left(P^\varepsilon(t,\cdot)^2,\, Q^\varepsilon(t,\cdot)^2 \right)\ \rightharpoonup\ \left(\overline{P^2}(t,\cdot),\, \overline{Q^2}(t,\cdot) \right) \qquad \mathrm{when}\ \varepsilon\, \to\, 0 \qquad \mathrm{in\ } \mathcal{D}'(\mathds{R}).
\end{gather*}
This, with Jensen's inequality, \eqref{energyconservationep} and \eqref{Limene} imply that 
\begin{align}\nonumber
\|W(t)\|_{L^2(\mathds{R})}^2\ &\leqslant\ \left\|\widetilde{W}(t)\right\|_{L^2(\mathds{R})}^2\  \leqslant\
\liminf_{\varepsilon \to 0} \|W^\varepsilon(t)\|_{L^2(\mathds{R})}^2\\ \label{Energyfin}
&=\ \liminf_{\varepsilon \to 0}  \int_\mathds{R} {\mathscr{E}}^\varepsilon(t,x)\, \mathrm{d}x\ \leqslant\ \lim_{\varepsilon \to 0} \int_\mathds{R} {\mathscr{E}}_0^\varepsilon\, \mathrm{d}x\ =\ \int_\mathds{R} {\mathscr{E}}_0\, \mathrm{d}x\ =\ \|W_0\|_{L^2}^2. 
\end{align}
The energy inequality \eqref{energyconservationep} with \eqref{Limene} imply that $(u^\varepsilon, P^\varepsilon, Q^\varepsilon)$ is bounded in the space $L^\infty([0,\infty), L^2(\mathds{R}))$.
We multiply \eqref{Appsys1} by $1$, \eqref{Appsys2} by $(h^\varepsilon)^{1/2}$ and \eqref{Ricatti1_ep}, \eqref{Ricatti2_ep} by $(h^\varepsilon)^{-1/2}$ we obtain
\begin{gather*}
h_t^\varepsilon\ +\,\left[\, h^\varepsilon\,u^\varepsilon\,\right]_x\ =\ \mathscr{A}^\varepsilon_x, \\
\left[ \sqrt{h^\varepsilon} u^\varepsilon\right]_t + {\textstyle \frac{\left[\, h^\varepsilon\,u^\varepsilon\,\right]_x\, -\, \mathscr{A}^\varepsilon_x}{2\, (h^\varepsilon)^{1/2}}}\, u^\varepsilon\, +\, (h^\varepsilon)^{1/2}\, u^\varepsilon\, u^\varepsilon_x + 3\, \gamma (h^\varepsilon)^{-3/2} h^\varepsilon_x\, =\, - (h^\varepsilon)^{1/2} \mathcal{L}_{h^\varepsilon}^{-1}  \partial_x\, \mathscr{C}^\varepsilon\, +\, (h^\varepsilon)^{1/2} \mathscr{B}^\varepsilon,\\
\left[{\textstyle \frac{P^\varepsilon}{\sqrt{h^\varepsilon}}}\right]_t\, +\, \left[{\textstyle \frac{\lambda^\varepsilon\, P^\varepsilon}{\sqrt{h^\varepsilon}}}\right]_x\, =\ {\textstyle \frac{1}{8\, (h^\varepsilon)^{3/2}}} \left[ \left(P^\varepsilon\right)^2\ +\  \chi_\varepsilon(P^\varepsilon)\ +\ \left(Q^\varepsilon \right)^2\ +\ 10\, P^\varepsilon\, Q^\varepsilon  -\ 8\, \mathscr{A}^\varepsilon_x\,  P^\varepsilon \right]\
 +\ {\textstyle \frac{\mathscr{M}^\varepsilon}{\sqrt{h^\varepsilon}}},\\
\left[{\textstyle \frac{Q^\varepsilon}{\sqrt{h^\varepsilon}}}\right]_t\, +\, \left[{\textstyle \frac{\eta^\varepsilon\, Q^\varepsilon}{\sqrt{h^\varepsilon}}}\right]_x\, =\ {\textstyle \frac{1}{8\, (h^\varepsilon)^{3/2}}} \left[ \left(Q^\varepsilon\right)^2\ +\  \chi_\varepsilon(Q^\varepsilon)\ +\ \left(P^\varepsilon \right)^2\ +\ 10\, P^\varepsilon\, Q^\varepsilon  -\ 8\, \mathscr{A}^\varepsilon_x\,  Q^\varepsilon \right]\
 +\ {\textstyle \frac{\mathscr{N}^\varepsilon}{\sqrt{h^\varepsilon}}}.
\end{gather*}
Then for all $T>0$ and for all $\varphi \in H^1(\mathds{R})$, the map 
\begin{equation}
t\ \mapsto\ \int_\mathds{R} \varphi(x)\, W^\varepsilon\, \mathrm{d}x
\end{equation}
is uniformly (on $t \in [0,T]$ and $\varepsilon \leqslant \varepsilon_0$) continuous. Then, Lemma \ref{lem:C_w} implies that 
\begin{equation}\label{Wconv2}
W(t,\cdot)\ \rightharpoonup\ W_0\qquad \mathrm{when}\ t\, \to\, 0 \qquad \mathrm{in\ } L^2(\mathds{R}),
\end{equation}
which implies that 
\begin{equation*}
\int_\mathds{R} {\mathscr{E}}_0\, \mathrm{d}x\ =\ \|W_0\|_{L^2}^2\ \leqslant\ \liminf_{t \to 0}  \|W\|_{L^2}^2.
\end{equation*}
On another hand, \eqref{Energyfin} implies
\begin{equation*}
\limsup_{t \to 0}  \|W\|_{L^2}^2\ =\ \limsup_{t \to 0}  \int_\mathds{R} {\mathscr{E}}\, \mathrm{d}x\ \leqslant\ \int_\mathds{R} {\mathscr{E}}_0\, \mathrm{d}x\ =\ \|W_0\|_{L^2}^2,
\end{equation*}
then 
\begin{equation}\label{Normconv}
\lim_{t \to 0}  \|W\|_{L^2}^2\ =\ \|W_0\|_{L^2}^2\ =\ \int_\mathds{R} {\mathscr{E}}_0\, \mathrm{d}x, 
\end{equation}
which implies with \eqref{Wconv2} that 
\begin{equation}\label{Sconv}
W(t,\cdot)\ \to\ W_0\qquad \mathrm{when}\ t\, \to\, 0 \qquad \mathrm{in\ } L^2(\mathds{R}).
\end{equation}
The inequality \eqref{Energyfin} with \eqref{Normconv} imply 
\begin{equation}\label{Normconv2}
\lim_{t \to 0}  \left\|\widetilde{W}\right\|_{L^2}^2\ =\ \lim_{t \to 0}  \|W\|_{L^2}^2\ =\ \|W_0\|_{L^2}^2.
\end{equation}
Then \eqref{tauto0_2} follows directly from \eqref{Normconv2}. 
Using \eqref{Sconv} and \eqref{uh_xi} we obtain the strong convergence 
\begin{align*}
\left\| \left(h\, -\, h_0,\ \sqrt{h}\, u\, -\, \sqrt{h_0}\, u_0,\ h_x/h\, -\, h_0'/h_0,\ \sqrt{h}\, u_x\, -\, \sqrt{h_0}\, u_0' \right) \right\|_{L^2}\ \to\ 0,
\end{align*}
as $t \to 0$.
In order to obtain \eqref{tauto0}, we write 
\begin{align*}
u\, -\, u_0\, &=\, {\textstyle \frac{1}{\sqrt{h}}} \left( \sqrt{h}\, u\, -\, \sqrt{h_0}\, u_0 \right) +\, \sqrt{h_0}\, u_0 \left({\textstyle \frac{1}{\sqrt{h}}}\, -\, {\textstyle \frac{1}{\sqrt{h_0}}} \right), \\
h_x\, -\, h'_0\, &=\, h \left( {\textstyle \frac{h_x}{\sqrt{h}}}\, -\, {\textstyle \frac{h_0'}{\sqrt{h_0}}} \right) +\, {\textstyle \frac{h_0'}{\sqrt{h_0}}} \left( h\, -\, h_0 \right), \\
u_x\, -\, u'_0\, &=\, {\textstyle \frac{1}{\sqrt{h}}} \left( \sqrt{h}\, u_x\, -\, \sqrt{h_0}\, u'_0 \right) +\, \sqrt{h_0}\, u'_0 \left({\textstyle \frac{1}{\sqrt{h}}}\, -\, {\textstyle \frac{1}{\sqrt{h_0}}} \right).
\end{align*}
On the right-hand side of the previous equations, the first term converges to $0$ in $L^2$ as $t \to 0$.
Since $h,1/h \in L^\infty$, $u_0,h_0',u_0' \in L^2$ and $h \in C ([0,+\infty) \times \mathds{R})$, the dominated convergence theorem implies that the second term goes to $0$ as $t \to 0$. This ends the proof of  \eqref{tauto0}.
\qed

For any $\kappa >0$, we define
\begin{equation}\label{Sdef}
S_\kappa(\xi)\ \eqdef\ \half\, \xi^2\ -\ \half\, (\xi\, +\, \kappa)^2\, \mathds{1}_{\xi\leqslant-\kappa}\ -\ \half\, (\xi\, -\, \kappa)^2\, \mathds{1}_{\xi\geqslant \kappa}\ =\ 
\begin{cases}
-\kappa\, (\xi\, +\, \half\, \kappa), & \xi \leqslant -\kappa, \\
\half\, \xi^2, & |\xi| \leqslant \kappa, \\
\kappa\, (\xi\, -\, \half\, \kappa), & \xi \geqslant \kappa.
\end{cases}
\end{equation}
\begin{equation}\label{Tdef}
T_\kappa(\xi)\ \eqdef\ S_\kappa'(\xi)\ =\ \xi\ -\ (\xi\, +\, \kappa)\, \mathds{1}_{\xi\leqslant-\kappa}\ -\ (\xi\, -\, \kappa)\, \mathds{1}_{\xi\geqslant \kappa}\ =\ 
\begin{cases}
-\kappa, & \xi \leqslant -\kappa, \\
\xi, & |\xi| \leqslant \kappa, \\
\kappa, & \xi \geqslant \kappa.
\end{cases}
\end{equation}

\begin{lem}\label{BTlim}
For any $\kappa>0$, there exists a subsequence $\left\{ \mathscr{M}^\varepsilon, \mathscr{N}^\varepsilon,P^\varepsilon, Q^\varepsilon \right\}_\varepsilon$ and $\widetilde{\mathscr{M}} \in L^\infty_{loc}((0,\infty) \times \mathds{R})$ such that, when $\varepsilon \to 0$ we have the limits in the sense of distributions on $(0,\infty) \times \mathds{R}$
\begin{align}\label{tildeB1}
\mathscr{M}^\varepsilon\ \rightharpoonup\ \widetilde{\mathscr{M}}, \qquad \mathrm{and} \qquad \mathscr{M}^\varepsilon\, T_\kappa \left( P^\varepsilon \right)\ \rightharpoonup\  \overline{T_\kappa \left( P \right)}\, \widetilde{\mathscr{M}}, \\ \label{tildeB2}
\mathscr{N}^\varepsilon\ \rightharpoonup\ \widetilde{\mathscr{M}}, \qquad \mathrm{and} \qquad \mathscr{N}^\varepsilon\, T_\kappa \left( Q^\varepsilon \right)\ \rightharpoonup\ \overline{T_\kappa \left( Q \right)}\, \widetilde{\mathscr{M}}.
\end{align}
\end{lem}
\proof 
\textbf{Step 1.} 
We define 
\begin{equation}\label{CC}
\mathscr{P}^\varepsilon\ \eqdef\ \mathcal{L}_{h^\varepsilon}^{-1} \left(h^\varepsilon \int_{-\infty}^x (h^\varepsilon)^{-3}\, \mathscr{C}^\varepsilon\, \mathrm{d}y\, +\, \third\, F(h^\varepsilon)_x \right) .
\end{equation}
From \eqref{Estimate6}, we have that  $\mathscr{P}^\varepsilon$ is bounded in $ L^\infty([t_1,t_2], W^{1,\infty}([a,b]))$ for any $b>a, t_2 > t_1>0$. Thus, there exists $\widetilde{\mathscr{P}} \in L^\infty([t_1,t_2], W^{1,\infty}([a,b]))$ such that, up to a subsequence we have 
\begin{equation}
\mathscr{P}^\varepsilon\ \rightharpoonup\ \widetilde{\mathscr{P}}, \qquad
\partial_x \mathscr{P}^\varepsilon\ \rightharpoonup\  \partial_x \widetilde{\mathscr{P}}, \quad \mathrm{as}\ \varepsilon \to 0
\end{equation}
in $L^p_{loc}((0,\infty) \times \mathds{R})$ for any $p <  \infty$.

\textbf{Step 2.}
For a fixed $\varphi \in C^\infty_c((0,\infty) \times \mathds{R})$, the inequality \eqref{estimate1.5}, Lemma \ref{lem:bounded2} and \eqref{Balpha+2} imply that $\left(1 - \partial_x^2\right) \left\{ \varphi \mathscr{P}^\varepsilon\right\}$ is uniformly bounded in $L_{loc}^{2+\alpha}((0,\infty) \times \mathds{R})$ for all $\alpha \in [0,1).$ Then, up to a subsequence we have 
\begin{equation*}
\left(1 - \partial_x^2\right)\, \left\{ \varphi\, \mathscr{P}^\varepsilon\right\}\ \rightharpoonup\ \left(1 - \partial_x^2\right)\, \left\{ \varphi\, \widetilde{\mathscr{P}}\right\}
\end{equation*}
in $L_{loc}^{2+\alpha}((0,\infty) \times \mathds{R})$.

\textbf{Step 3.}
Since $|T_\kappa(P^\varepsilon)| \leqslant \kappa$, the convergence $T_\kappa(P^\varepsilon) \rightharpoonup \overline{T_\kappa(P)} $ is in $L^p_{loc}((0,\infty) \times \mathds{R})$ for any $p \in (1,\infty)$. Then, for any $\psi \in C^\infty_c((0,\infty) \times \mathds{R})$ we have up to a subsequence
\begin{equation*}
\lim_{\varepsilon \to 0} \int_{(0,\infty) \times \mathds{R}}\! \psi\, \left(1 - \partial_x^2\right)^{-1} \left\{ \partial_x T_\kappa(P^\varepsilon)\right\}\, \mathrm{d}x\, \mathrm{d}t\ =\ \int_{(0,\infty) \times \mathds{R}}\! \psi\, \left(1 - \partial_x^2\right)^{-1} \left\{ \partial_x \overline{T_\kappa(P)}\right\}\,  \mathrm{d}x\, \mathrm{d}t .
\end{equation*}
This limit is stronger. Indeed, replacing $f$ in  \eqref{rhs0} by $T_\kappa$ we obtain
\begin{gather*}
\left[T_\kappa(P^\varepsilon)\right]_t\ +\ \left[ \lambda^\varepsilon\, T_\kappa(P^\varepsilon) \right]_x\ =\ {\textstyle \frac{1}{4\, h^\varepsilon}} \left(P^\varepsilon +3\, Q^\varepsilon \right) T_\kappa(P^\varepsilon)\\
+\
\left[ - {\textstyle \frac{1}{8\, h^\varepsilon}}\, \left(P^\varepsilon\right)^2\ +\ {\textstyle \frac{1}{8\, h^\varepsilon}}\, \chi_\varepsilon(P^\varepsilon)\ +\ {\textstyle \frac{1}{8\, h^\varepsilon}}\, \left(Q^\varepsilon\right)^2\ -\ {\textstyle \frac{1}{2\, h^\varepsilon}}\, \mathscr{A}^\varepsilon_x\,  P^\varepsilon\ +\ \mathscr{M}^\varepsilon\right] T_\kappa'(P^\varepsilon).
\end{gather*}
Then, the sequences $\left\{ \left(1 - \partial_x^2\right)^{-1} \left\{ \partial_x T_\kappa(P^\varepsilon)\right\} \right\}_\varepsilon $ is uniformly bounded in $W^{1,\infty}((0,\infty) \times \mathds{R})$. Arzela--Ascoli theorem implies that, up to a subsequence, we have the convergence 
\begin{equation}
\left(1 - \partial_x^2\right)^{-1} \left\{ \partial_x T_\kappa(P^\varepsilon)\right\}\ \longrightarrow\ \left(1 - \partial_x^2\right)^{-1} \left\{ \partial_x \overline{T_\kappa(P)}\right\}
\end{equation}
is uniform on any compact set of $(0,\infty) \times \mathds{R}$. Doing the same proof again we obtain the uniform convergence
\begin{equation}
\left(1 - \partial_x^2\right)^{-1} \left\{ T_\kappa(P^\varepsilon)\right\}\ \longrightarrow\ \left(1 - \partial_x^2\right)^{-1} \left\{ \overline{T_\kappa(P)}\right\}
\end{equation}
on any compact set of $(0,\infty) \times \mathds{R}$.

\textbf{Step 4.}
Let $\varphi \in C^\infty_c((0,\infty) \times \mathds{R})$, then 
\begin{align*}
\int_{(0,\infty) \times \mathds{R}}\! T_\kappa(P^\varepsilon)\, \varphi\, \mathscr{P}^\varepsilon_x\, \mathrm{d}x\, \mathrm{d}t\ 
=&\
\int_{(0,\infty) \times \mathds{R}}\! T_\kappa(P^\varepsilon)\, \left(1 - \partial_x^2\right)^{-1}\, \left(1 - \partial_x^2\right)\, \left[ (\varphi\, \mathscr{P}^\varepsilon)_x\, -\, \varphi_x\, \mathscr{P}^\varepsilon \right]\, \mathrm{d}x\, \mathrm{d}t\\
=&\ - \int_{(0,\infty) \times \mathds{R}}\! \left(1 - \partial_x^2\right)^{-1}\, \left\{ \partial_x\, T_\kappa(P^\varepsilon)\right\}\, \cdot\, \left(1 - \partial_x^2\right)\, \left\{ \varphi\, \mathscr{P}^\varepsilon\right\}\, \mathrm{d}x\, \mathrm{d}t\\
&- \int_{(0,\infty) \times \mathds{R}}\! \left(1 - \partial_x^2\right)^{-1}\, \left\{ T_\kappa(P^\varepsilon)\right\}\, \cdot\, \left(1 - \partial_x^2\right)\, \left\{ \varphi_x\, \mathscr{P}^\varepsilon\right\}\, \mathrm{d}x\, \mathrm{d}t .
\end{align*}
Taking the limit $\varepsilon \to 0$ and using Step 2, Step 3 and  Lemma \ref{lem:fg} we obtain
\begin{align}\nonumber
\lim_{\varepsilon \to 0}\int_{(0,\infty) \times \mathds{R}}\! T_\kappa(P^\varepsilon)\, \varphi\, \mathscr{P}^\varepsilon_x\, \mathrm{d}x\, \mathrm{d}t\ 
=&\ - \int_{(0,\infty) \times \mathds{R}}\! \left(1 - \partial_x^2\right)^{-1}\, \left\{ \partial_x\, \overline{T_\kappa(P)}\right\}\, \cdot\, \left(1 - \partial_x^2\right)\, \left\{ \varphi\, \widetilde{\mathscr{P}}\right\}\, \mathrm{d}x\, \mathrm{d}t\\ \nonumber
&\ - \int_{(0,\infty) \times \mathds{R}}\! \left(1 - \partial_x^2\right)^{-1}\, \left\{ \overline{T_\kappa(P)}\right\}\, \cdot\, \left(1 - \partial_x^2\right)\, \left\{ \varphi_x\, \widetilde{\mathscr{P}}\right\}\, \mathrm{d}x\, \mathrm{d}t \\ \label{TPx}
=&\ \int_{(0,\infty) \times \mathds{R}}\! \overline{T_\kappa(P)}\, \varphi\, \partial_x\,  \widetilde{\mathscr{P}}\, \mathrm{d}x\, \mathrm{d}t.
\end{align}
\textbf{Step 5.}
Since $|T_\kappa \left( P^\varepsilon \right)| \leqslant \kappa$, then, from \eqref{Estimate9} we have $ \mathscr{V}_1 T_\kappa \left( P^\varepsilon \right) \rightharpoonup 0$ and $ \mathscr{V}_2 T_\kappa \left( P^\varepsilon \right) \rightharpoonup 0$. Then, using \eqref{TPx} we obtain \eqref{tildeB1} with 
$$\widetilde{\mathscr{M}}\ \eqdef\ -3\, h\, \partial_x\, \widetilde{\mathscr{P}}. $$
Following the same proof we obtain \eqref{tildeB2}.
\qed

\begin{lem}
For all $T>0$, we have
\begin{gather}\label{L1conv}
\lim_{\kappa \to \infty} \left\| \overline{T_\kappa(P)}\ -\ T_\kappa\left(\overline{P}\right) \right\|_{L^1([0,T] \times \mathds{R})}\ =\
\lim_{\kappa \to \infty} \left\| \overline{T_\kappa(Q)}\ -\ T_\kappa\left(\overline{Q}\right) \right\|_{L^1([0,T] \times \mathds{R})}\ =\ 0,\\ \label{L1conv2}
\lim_{\kappa \to \infty} \left\| \overline{T_\kappa(P)}\ -\ \overline{P} \right\|_{L^1([0,T] \times \mathds{R})}\ =\
\lim_{\kappa \to \infty} \left\| \overline{T_\kappa(Q)}\ -\ \overline{Q} \right\|_{L^1([0,T] \times \mathds{R})}\ =\ 0.
\end{gather}
Moreover, for all $\kappa >0$ we have
\begin{align}\label{TS}
\half\, \left( \overline{T_\kappa(P)}\ -\ T_\kappa\left(\overline{P}\right) \right)^2\ 
&\leqslant\ \overline{S_\kappa(P)}\ -\ S_\kappa\left(\overline{P}\right),\\ 
\half\, \left( \overline{T_\kappa(Q)}\ -\ T_\kappa\left(\overline{Q}\right) \right)^2\ 
&\leqslant\ \overline{S_\kappa(Q)}\ -\ S_\kappa\left(\overline{Q}\right). 
\end{align}
\end{lem}

\proof 
Since the proof for $P$ and $Q$ is the same, we only do the proof for $P$.
From \eqref{Tdef} we have 
\begin{align*}
\left| T_\kappa\left(\xi\right)\ - \xi\right|\ \leqslant\
\left|\xi\, +\, \kappa \right|\, \mathds{1}_{\xi\leqslant-\kappa}\ +\ \left|\xi\, -\, \kappa \right|\, \mathds{1}_{\xi\geqslant \kappa}\
\leqslant\ 2\, \left|\xi\right|\, \mathds{1}_{\kappa \leqslant \left|\xi\right|}\ 
\leqslant\ {\textstyle \frac{2}{\kappa}}\, \xi^2.
\end{align*}
Then, we have 
\begin{equation}
\left| \overline{T_\kappa\left(P\right)}\ - T_\kappa\left(\overline{P}\right) \right|\ \leqslant\ 
\left| \overline{T_\kappa\left(P\right)}\ - \overline{P} \right|\ +\
\left| T_\kappa\left(\overline{P}\right)\ - \overline{P} \right|\
\leqslant {\textstyle \frac{2}{\kappa}}\, \left( \overline{P^2}\, +\, \overline{P}^2 \right).
\end{equation}
Jenson's inequality imply that $\overline{P}^2 \leqslant \overline{P^2}$. Lemma \ref{lemYoungm} implies that $\overline{P^2} \in L^\infty(\mathds{R}^+ , L^1(\mathds{R}))$. Then \eqref{L1conv} and \eqref{L1conv2} follow directly.

Cauchy--Schwarz inequality implies that $\overline{T_\kappa(P)}^2 \leqslant \overline{T_\kappa(P)^2}$, then, using the definition \eqref{Tdef} we obtain
\begin{align} \nonumber
\left( \overline{T_\kappa(P)}\ -\ T_\kappa\left( \overline{P}\right) \right)^2\ 
\leqslant &\ \overline{T_\kappa(P)^2}\ +\ T_\kappa\left( \overline{P}\right)^2\ -\ 2\, T_\kappa\left( \overline{P}\right)\, \overline{T_\kappa(P)}\\ \nonumber
=&\
\overline{T_\kappa(P)^2}\ +\ T_\kappa\left( \overline{P}\right)^2\ -\ 2\, T_\kappa\left( \overline{P}\right)\, \overline{P}\ +\ 2\, T_\kappa\left( \overline{P}\right)\, \overline{ \left( P\, +\, \kappa \right)\, \mathds{1}_{P \leqslant -\kappa}}\\ \nonumber
&+\ 2\, T_\kappa(\overline{P})\, \overline{ \left( P\, -\, \kappa \right)\, \mathds{1}_{P \geqslant \kappa}}\\ \nonumber
=&\ \overline{T_\kappa(P)^2}\ +\ 2\, T_\kappa\left( \overline{P}\right)\, \left[ \overline{ \left( P\, +\, \kappa \right)\, \mathds{1}_{P \leqslant -\kappa}}\ -\ \left( \overline{P}\, +\, \kappa \right)\, \mathds{1}_{\overline{P} \leqslant - \kappa} \right]\\ \nonumber
& - T_\kappa\left( \overline{P}\right)^2\ +\ 2\, T_\kappa\left( \overline{P}\right)\, \left[ \overline{ \left( P\, -\, \kappa \right)\, \mathds{1}_{P \geqslant \kappa}}\ -\ \left( \overline{P}\, -\, \kappa \right)\, \mathds{1}_{\overline{P} \geqslant \kappa} \right]\\ \nonumber
\leqslant &\ \overline{T_\kappa(P)^2}\ -\ 2\, \kappa\, \left[ \overline{ \left( P\, +\, \kappa \right)\, \mathds{1}_{P \leqslant -\kappa}}\ -\ \left( \overline{P}\, +\, \kappa \right)\, \mathds{1}_{\overline{P} \leqslant - \kappa} \right]\\ \label{TT} 
& - T_\kappa\left( \overline{P}\right)^2\ +\ 2\, \kappa\, \left[ \overline{ \left( P\, -\, \kappa \right)\, \mathds{1}_{P \geqslant \kappa}}\ -\ \left( \overline{P}\, -\, \kappa \right)\, \mathds{1}_{\overline{P} \geqslant \kappa} \right],
\end{align}
where the last inequality follows from Jensen's inequality with the concavity of $\xi \mapsto \left( \xi + \kappa \right) \mathds{1}_{\xi \leqslant - \kappa}$, the convexity of $\xi \mapsto \left( \xi - \kappa \right) \mathds{1}_{\xi \geqslant  \kappa}$ and $-\kappa \leqslant T_\kappa(\xi)\leqslant \kappa$. Since 
\begin{equation*}
S_\kappa(\xi)\ =\ \half\, T_\kappa(\xi)^2\ +\ \kappa\, \left( \xi\, -\, \kappa \right)\, \mathds{1}_{\xi \geqslant  \kappa}\ -\ \kappa\, \left( \xi\, +\, \kappa \right)\, \mathds{1}_{\xi \leqslant - \kappa}
\end{equation*}
we have 
\begin{align*}
\overline{S_\kappa(P)}\ &=\ \half\, \overline{T_\kappa(P)^2}\ +\ \kappa\, \overline{ \left( P\, -\, \kappa \right)\, \mathds{1}_{P \geqslant  \kappa}}\ -\ \kappa\, \overline{ \left( P\, +\, \kappa \right)\, \mathds{1}_{P \leqslant - \kappa}},\\
S_\kappa\left(\overline{P}\right)\ &=\ \half\, T_\kappa(\overline{P})^2\ +\ \kappa\, \left( \overline{P}\, -\, \kappa \right)\, \mathds{1}_{\overline{P} \geqslant  \kappa}\ -\ \kappa\, \left( \overline{P}\, +\, \kappa \right)\, \mathds{1}_{\overline{P} \leqslant - \kappa}.
\end{align*}
The last two identities with \eqref{TT} imply \eqref{TS}.
\qed

Now we state the main result of this section.
\begin{lem}\label{lem:2Dirac}
 The measures $\nu^1, \nu^2$ given in Lemma \ref{lemYoungm} are Dirac measures, and
\begin{equation}
\nu^1_{t,x}(\xi)\ =\ \delta_{\overline{P}(t,x)}(\xi),\ \qquad \nu^2_{t,x}(\zeta)\ =\ \delta_{\overline{Q}(t,x)}(\zeta).
\end{equation}
\end{lem}
\proof 
Since the proof is the same, we present here only the proof of $\nu^1_{t,x}(\xi) = \delta_{\overline{P}(t,x)}(\xi)$.
Note that if $\overline{P^2}=\overline{P}^2$ then $\int_\mathds{R} \left(\overline{P}-\xi\right)^2 \mathrm{d}\nu^1_{t,x}(\xi)=0$ which implies that $\mathrm{supp}(\nu^1_{t,x})= \left\{\overline{P}\right\}$. Since $\nu^1_{t,x}$ is a probability measure, necessarily $\nu^1_{t,x} = \delta_{\overline{P}}$. It remains then only to prove that $\overline{P^2}=\overline{P}^2$.
The goal is to obtain an evolutionary inequality of $\overline{P^2}-\overline{P}^2$, then, since it is equal to zero initially, we prove that it remains zero for all time.
The proof is given in several steps.

\textbf{Step 1.}
Replacing $f$ in  \eqref{rhs0} by $S_\kappa$ one obtains
\begin{gather*}\nonumber
S_\kappa(P^\varepsilon)_t\ +\ \left[ \lambda^\varepsilon\, S_\kappa(P^\varepsilon) \right]_x\ =\ {\textstyle \frac{1}{4\, h^\varepsilon}} \left(P^\varepsilon +3\, Q^\varepsilon \right) S_\kappa(P^\varepsilon)\\
+\
\left[ - {\textstyle \frac{1}{8\, h^\varepsilon}}\, \left(P^\varepsilon\right)^2\ +\ {\textstyle \frac{1}{8\, h^\varepsilon}}\, \chi_\varepsilon(P^\varepsilon)\ +\ {\textstyle \frac{1}{8\, h^\varepsilon}}\, \left(Q^\varepsilon\right)^2\ -\ {\textstyle \frac{1}{2\, h^\varepsilon}}\, \mathscr{A}^\varepsilon_x\,  P^\varepsilon\ +\ \mathscr{M}^\varepsilon\right] T_\kappa(P^\varepsilon).
\end{gather*}
Taking $\varepsilon \to 0$, using \eqref{Estimate9}, Lemma \ref{lemYoungm} and Lemma \ref{BTlim} we obtain
\begin{gather}\nonumber
\overline{S_\kappa(P)}_t\ +\ \left( \lambda\, \overline{S_\kappa(P)} \right)_x\ =\\ \label{Step1eq}
{\textstyle \frac{1}{8\, h}}\, \left\{ 2\, \overline{P\, S_\kappa(P)}\ -\ \overline{P^2\, T_\kappa(P)}\ +\ 6\, \overline{Q}\, \overline{S_\kappa(P)}\ +\ \overline{Q^2}\, \overline{T_\kappa(P)} \right\}\ +\ \overline{T_\kappa(P)}\, \widetilde{\mathscr{M}}.
\end{gather}
\textbf{Step 2.} Replacing $f$ in  \eqref{rhs0} by the identity function and taking $\varepsilon \to 0$ we obtain 
\begin{equation}
\overline{P}_t\ +\ \left(\lambda\, \overline{P}\right)_x\ =\ {\textstyle \frac{1}{8\, h}}\, \left( \overline{P^2}\, +\, 6\, \overline{P}\, \overline{Q}\, +\, \overline{Q^2} \right)\ +\ \widetilde{\mathscr{M}}.
\end{equation}
Let $j_\varepsilon$ be a Friedrichs mollifier, we note $\overline{P}^\varepsilon \eqdef \overline{P} \ast j_\varepsilon$, then we have 
\begin{equation}
\overline{P}^\varepsilon_t\ +\ \left(\lambda\, \overline{P}^\varepsilon\right)_x\ =\ \theta_\varepsilon\ +\ \left\{ {\textstyle \frac{1}{8\, h}}\, \left( \overline{P^2}\, +\, 6\, \overline{P}\, \overline{Q}\, +\, \overline{Q^2} \right)\right\}\, \ast j_\varepsilon \ +\ \widetilde{\mathscr{M}}\, \ast j_\varepsilon,
\end{equation}
where $\theta_\varepsilon \eqdef \left(\lambda\, \overline{P}^\varepsilon\right)_x -\left(\lambda\, \overline{P}\right)_x \ast j_\varepsilon $. 
Multiplying by $T_\kappa\left(\overline{P}^\varepsilon\right)$ and using \eqref{uh_xi}, we obtain
\begin{gather*}
S_\kappa\left(\overline{P}^\varepsilon\right)_t\ +\ \left(\lambda\, S_\kappa\left(\overline{P}^\varepsilon\right)\right)_x\
=\  {\textstyle \frac{1}{4\, h}}\, (3\, \overline{Q}\, +\, \overline{P})\, S_\kappa\left(\overline{P}^\varepsilon\right)\ -\ {\textstyle \frac{1}{4\, h}}\, \left(3\, \overline{Q}\, +\, \overline{P}\right)\, \overline{P}^\varepsilon\, T_\kappa\left(\overline{P}^\varepsilon\right)\\
+\ T_\kappa\left(\overline{P}^\varepsilon\right)   \left\{ {\textstyle \frac{1}{8\, h}}\, \left( \overline{P^2}\, +\, 6\, \overline{P}\, \overline{Q}\, +\, \overline{Q^2} \right)\right\}\, \ast j_\varepsilon \ +\ T_\kappa\left(\overline{P}^\varepsilon\right)\, \left( \widetilde{\mathscr{M}}\, \ast j_\varepsilon \right)\, +\ T_\kappa\left(\overline{P}^\varepsilon\right)\, \theta_\varepsilon.
\end{gather*}
Taking $\varepsilon \to 0$ and using Lemma \ref{lemRenor}, one obtains
\begin{gather}\nonumber
S_\kappa\left(\overline{P}\right)_t\ +\ \left(\lambda\, S_\kappa\left(\overline{P}\right)\right)_x\ =\\ \label{Step2eq}
{\textstyle \frac{1}{8\, h}}\, \left\{2\, \overline{P}\, S_\kappa\left(\overline{P}\right)\ +\ 6\, \overline{Q}\, S_\kappa\left(\overline{P}\right)\ +\ T_\kappa\left(\overline{P}\right) \left(\overline{P^2}\, -\, 2\, \overline{P}^2\, +\, \overline{Q^2} \right) \right\}\ +\ T_\kappa\left(\overline{P}\right)\, \widetilde{\mathscr{M}}.
\end{gather}
\textbf{Step 3.} From \eqref{Step1eq} and \eqref{Step2eq} we obtain
\begin{gather}\nonumber
\left[\overline{S_\kappa(P)}\, -\, S_\kappa \left(\overline{P}\right) \right]_t\ + \left[ \lambda\, \left(\overline{S_\kappa(P)}\, -\, S_\kappa\left(\overline{P}\right) \right) \right]_x\ 
=\ \widetilde{\mathscr{M}} \left( \overline{T_\kappa(P)}\
-\ T_\kappa\left(\overline{P}\right) \right)\\ \nonumber
 {\textstyle \frac{1}{8\, h}}\,  \left\{ 2\, \overline{P\, S_\kappa(P)}\  -\ \overline{P^2\, T_\kappa(P)}\ +\ \overline{P}^2\, T_\kappa\left(\overline{P}\right)\ -\ 2\, \overline{P}\, S_\kappa\left(\overline{P}\right)\ +\ T_\kappa\left(\overline{P}\right) \left( \overline{P}^2\ -\ \overline{P^2}\ \right) 
\right\}\\ \label{Step12eq}
+\ {\textstyle \frac{1}{8\, h}}\,  \left\{ 
6\, \overline{Q}\, \left(\overline{S_\kappa(P)}\, -\,  S_\kappa\left(\overline{P}\right)\right)\ +\ \overline{Q^2}\, \left(\overline{T_\kappa(P)}\, -\, T_\kappa\left(\overline{P}\right) \right)
\right\}.
\end{gather}
From \eqref{Sdef} and \eqref{Tdef} we have
\begin{align*}
\xi^2\, T_\kappa(\xi)\ -\ 2\, \xi\, S_\kappa(\xi)\ 
=&\ 
\xi^2\, T_\kappa(\xi)\ -\ 2\, \xi\, S_\kappa(\xi)\ +\xi^3\ -\xi^3\\
=&\ \xi^2\, \left[T_\kappa(\xi)\, -\, \xi \right]\ +\ \xi\, (\xi\, +\, \kappa)^2\, \mathds{1}_{\xi \leqslant -\kappa}\ +\ \xi\, (\xi-\kappa)^2\, \mathds{1}_{\xi \geqslant \kappa}\\
=&\ \kappa^2\, \left[T_\kappa(\xi)\, -\, \xi \right]\ -\ \left( \xi^2\, -\, \kappa^2\right) \left[ (\xi\, +\, \kappa)\, \mathds{1}_{\xi \leqslant -\kappa}\ +\ (\xi-\kappa)\, \mathds{1}_{\xi \geqslant \kappa} \right]\\
&+\ \xi\, (\xi\, +\, \kappa)^2\, \mathds{1}_{\xi \leqslant -\kappa}\ +\ \xi\, (\xi-\kappa)^2\, \mathds{1}_{\xi \geqslant \kappa}\\
=&\ \kappa^2\, \left[T_\kappa(\xi)\, -\, \xi \right]\ +\ \kappa\, (\xi\, +\, \kappa)^2\, \mathds{1}_{\xi \leqslant -\kappa}\ -\ \kappa\, (\xi-\kappa)^2\, \mathds{1}_{\xi \geqslant \kappa}.
\end{align*}
Then from \eqref{Sdef} we have 
\begin{gather}\nonumber
2\, \overline{P\, S_\kappa(P)}\  -\ \overline{P^2\, T_\kappa(P)}\ +\ \overline{P}^2\, T_\kappa\left(\overline{P}\right)\ -\ 2\, \overline{P}\, S_\kappa\left(\overline{P}\right)\ +\ T_\kappa\left(\overline{P}\right) \left( \overline{P}^2\ -\ \overline{P^2}\ \right)
\\ \nonumber
=\ 
 \left( T_\kappa\left(\overline{P}\right)\, +\, \kappa \right)\, (\overline{P}\, +\, \kappa)^2\, \mathds{1}_{\overline{P} \leqslant -\kappa}\
+\ \left( T_\kappa\left(\overline{P}\right)\, -\, \kappa \right)\, (\overline{P}\, -\, \kappa)^2\, \mathds{1}_{\overline{P} \geqslant \kappa}  \\ \nonumber
 -\left( T_\kappa\left(\overline{P}\right)\, +\, \kappa \right)\, \overline{(P\, +\, \kappa)^2\, \mathds{1}_{P \leqslant -\kappa}}\
-\ \left( T_\kappa\left(\overline{P}\right)\, -\, \kappa \right)\, \overline{(P\, -\, \kappa)^2\, \mathds{1}_{P \geqslant \kappa}}\\ \label{Step3_1}
-\ \kappa^2\, \left( \overline{T_\kappa(P)}\, -\, T_\kappa\left(\overline{P}\right) \right)\ 
-\ 2\, T_\kappa \left(\overline{P} \right) \left( \overline{S_\kappa(P)}\, -\, S_\kappa\left(\overline{P}\right) \right)\ 
\end{gather}
From the definition \eqref{Tdef} we have 
\begin{equation}\label{Step3_15}
\left( T_\kappa\left(\overline{P}\right)\, +\, \kappa \right)\, (\overline{P}\, +\, \kappa)^2\, \mathds{1}_{\overline{P} \leqslant -\kappa}\
=\ \left( T_\kappa\left(\overline{P}\right)\, -\, \kappa \right)\, (\overline{P}\, -\, \kappa)^2\, \mathds{1}_{\overline{P} \geqslant \kappa}\ =\ 0
\end{equation}
Since $T_\kappa\left(\overline{P}\right) \geqslant - \kappa $, then 
\begin{equation}\label{Step3_2}
 -\left( T_\kappa\left(\overline{P}\right)\, +\, \kappa \right)\, \overline{(P\, +\, \kappa)^2\, \mathds{1}_{P \leqslant -\kappa}}\ \leqslant\ 0.
\end{equation}
Let $t_0>0$ and $\kappa \geqslant C(1+ t_0^{-1})$, then from Lemma \ref{lem:Oleinik}, we have for all $t \geqslant t_0$ that $P^\varepsilon \leqslant \kappa $ and $\overline{P} \leqslant \kappa$. Then, using the convexity of $T_\kappa$ on $(-\infty,\kappa)$ and the Jensen's inequality we obtain 
\begin{equation}\label{Step3_3}
-\kappa^2\, \left( \overline{T_\kappa(P)}\, -\, T_\kappa\left(\overline{P}\right) \right)\ \leqslant\ 0, \qquad \forall t \geqslant t_0, \quad   \kappa \geqslant C(1+ t_0^{-1}).
\end{equation}
We take again $t_0>0$ and $\kappa \geqslant C(1+ t_0^{-1})$, then for all $\varphi \in C^\infty_c((t_0,\infty) \times \mathds{R})$ we have 
\begin{equation}\label{Step3_4}
\int \overline{(P\, -\, \kappa)^2\, \mathds{1}_{P \geqslant \kappa}}\, \varphi\, \mathrm{d}x\, \mathrm{d}t\ =\ \lim_{\varepsilon \to 0} \int (P^\varepsilon\, -\, \kappa)^2\, \mathds{1}_{P^\varepsilon \geqslant \kappa}\, \varphi\, \mathrm{d}x\, \mathrm{d}t\ =\ 0.
\end{equation}
Then, using \eqref{Step12eq}, \eqref{Step3_1}, \eqref{Step3_15}, \eqref{Step3_2}, \eqref{Step3_3} and \eqref{Step3_4} we obtain that
\begin{gather}\nonumber
\left[\overline{S_\kappa(P)}\, -\, S_\kappa \left(\overline{P}\right) \right]_t\ + \left[ \lambda\, \left(\overline{S_\kappa(P)}\, -\, S_\kappa\left(\overline{P}\right) \right) \right]_x\ 
\leqslant\
\widetilde{\mathscr{M}} \left( \overline{T_\kappa(P)}\
-\ T_\kappa\left(\overline{P}\right) \right)\\ \label{Step3}
+\ {\textstyle \frac{1}{8\, h}}\,  \left\{ 
\left(6\, \overline{Q}\, -\, 2\, T_\kappa \left(\overline{P} \right) \right) \left(\overline{S_\kappa(P)}\, -\,  S_\kappa\left(\overline{P}\right)\right)\ +\ \overline{Q^2}\, \left(\overline{T_\kappa(P)}\, -\, T_\kappa\left(\overline{P}\right) \right)
\right\},
\end{gather}
for any  $t_0>0$,  $\kappa \geqslant C(1+ t_0^{-1})$ and  $t > t_0$.

\textbf{Step 4.} 
Let $f_\kappa(t,x) \eqdef \overline{S_\kappa(P)} - S_\kappa \left(\overline{P}\right)$ and $f_\kappa^\varepsilon \eqdef f_\kappa \ast j_\varepsilon$ where $j_\varepsilon$ is a Friedrichs mollifier, then, from \eqref{Step3} and Lemma \ref{lemRenor} we obtain 
\begin{equation*}
(f_\kappa^\varepsilon)_t\ +\ \left(\lambda\, f_\kappa^\varepsilon \right)_x\ \leqslant\ {\textstyle \frac{1}{4\, h}}\,  \left\{ 
 3\,  \overline{Q}\, -\, T_\kappa\left(\overline{P} \right)  \right\}\, f_\kappa^\varepsilon\ +\ \left( \widetilde{\mathscr{M}}\, +\, \frac{\overline{Q^2}}{8\, h} \right) \left(\overline{T_\kappa(P)}\, -\, T_\kappa\left(\overline{P}\right) \right)\ +\ \theta_\varepsilon,
\end{equation*} 
where $\theta_\varepsilon \to 0$ in $L^1_{loc}((0,\infty) \times \mathds{R})$. Let $\beta>0$, multiplying by $h^{3/2}\left( h^{3/2} f_\kappa^\varepsilon +  \beta \right)^{-1/2}/2$ and using \eqref{Appsys123} one obtains
\begin{gather*}
\left[ \sqrt{h^{3/2}\, f_\kappa^\varepsilon\, +\,  \beta} \right]_t\ +\ \left[ \lambda\, \sqrt{h^{3/2}\, f_\kappa^\varepsilon\, +\,  \beta} \right]_x\ 
\leqslant\ 
\left( \widetilde{\mathscr{M}}\, +\, {\textstyle\frac{\overline{Q^2}}{8\, h}} \right)
{\textstyle\frac{\overline{T_\kappa(P)}\, -\, T_\kappa\left(\overline{P} \right)}{2\, \sqrt{h^{3/2}\, f_\kappa^\varepsilon\, +\,  \beta}}}\, h^{3/2}\ +\ \tilde{\theta}_\varepsilon\\
+\ {\textstyle \frac{\left(\overline{P}\, -\, T_\kappa \left(\overline{P} \right) \right) h^{1/2}\, f_\kappa^\varepsilon}{8\, \sqrt{h^{3/2}\, f_\kappa^\varepsilon\, +\,  \beta}}}\ +\  {\textstyle \frac{ \beta\, \lambda_x}{\sqrt{h^{3/2}\, f_\kappa^\varepsilon\, +\,  \beta}}},
\end{gather*}
where $\tilde{\theta}_\varepsilon \eqdef \theta_\varepsilon h^{3/2}\left( h^{3/2} f_\kappa^\varepsilon +  \beta \right)^{-1/2}/2 \to 0$ in $L^1_{loc}((0,\infty) \times \mathds{R})$. Taking $\varepsilon \to 0$ we obtain
\begin{gather}\nonumber
\left[ \sqrt{h^{3/2}\, f_\kappa\, +\,  \beta} \right]_t\ +\ \left[ \lambda\, \sqrt{h^{3/2}\, f_\kappa\, +\,  \beta} \right]_x\ 
\leqslant\ 
\left( \widetilde{\mathscr{M}}\, +\, {\textstyle\frac{\overline{Q^2}}{8\, h}} \right)
{\textstyle\frac{\overline{T_\kappa(P)}\, -\, T_\kappa\left(\overline{P} \right)}{2\, \sqrt{h^{3/2}\, f_\kappa\, +\,  \beta}}}\, h^{3/2}\\ \label{fkappa}
+\ {\textstyle \frac{\left(\overline{P}\, -\, T_\kappa \left(\overline{P} \right) \right) h^{1/2}\, f_\kappa}{8\, \sqrt{h^{3/2}\, f_\kappa\, +\,  \beta}}}\ +\ {\textstyle \frac{ \beta\, \lambda_x}{\sqrt{h^{3/2}\, f_\kappa\, +\,  \beta}}}.
\end{gather}
From \eqref{TS} we have 
\begin{equation*}
\left|  \left( \widetilde{\mathscr{M}}\, +\, {\textstyle\frac{\overline{Q^2}}{8\, h}} \right)
{\textstyle\frac{\overline{T_\kappa(P)}\, -\, T_\kappa\left(\overline{P} \right)}{2\, \sqrt{h^{3/2}\, f_\kappa\, +\,  \beta}}}\, h^{3/2} \right| \
\leqslant\ {\textstyle \frac{\sqrt{2}}{2}}
\left|  \widetilde{\mathscr{M}}\, +\, {\textstyle \frac{\overline{Q^2}}{8\, h}}  \right| h^{3/4}.
\end{equation*}
Using that $|T_\kappa(\xi)| \leqslant |\xi|$ and $S_\kappa(\xi) \leqslant \xi^2/2$ we obtain 
\begin{equation*}
\left| {\textstyle \frac{\left(\overline{P}\, -\, T_\kappa \left(\overline{P} \right) \right) h^{1/2}\, f_\kappa}{8\, \sqrt{h^{3/2}\, f_\kappa\, +\,  \beta}}} \right| \ \leqslant\ {\textstyle \frac{ \left|\overline{P}\right| \sqrt{f_\kappa}}{4\, h^{1/4}}}\ \leqslant\ 
{\textstyle \frac{1}{8\, h^{1/4}}} \left( \overline{P}^2\, +\, f_\kappa \right)\ \leqslant\ 
{\textstyle \frac{1}{8\, h^{1/4}}} \left( {\textstyle \frac{3}{2}}\, \overline{P}^2\, +\, \half\, \overline{P^2}\right).
\end{equation*}
Since the $L^1$ convergence implies the pointwise convergence (up to a subsequence), then, using the dominated convergence theorem with \eqref{L1conv} and \eqref{L1conv2}, we obtain
\begin{equation*}
\lim_{\kappa\to \infty}\, \left\|  \left( \widetilde{\mathscr{M}}\, +\, {\textstyle\frac{\overline{Q^2}}{8\, h}} \right)
{\textstyle\frac{\overline{T_\kappa(P)}\, -\, T_\kappa\left(\overline{P} \right)}{2\, \sqrt{h^{3/2}\, f_\kappa\, +\,  \beta}}}\, h^{3/2} \right\|_{L^1(\Omega)}\ +\ 
\lim_{\kappa\to \infty}\, \left\| {\textstyle \frac{\left(\overline{P}\, -\, T_\kappa \left(\overline{P} \right) \right) h^{1/2}\, f_\kappa}{8\, \sqrt{h^{3/2}\, f_\kappa\, +\,  \beta}}} \right\|_{L^1(\Omega)}\
 =\ 0 
\end{equation*}
for any compact set $\Omega \subset (0,\infty) \times \mathds{R}$. Since $|S_\kappa(\xi)| \leqslant \xi^2/2$, then $|f_\kappa| \leqslant \overline{P}^2/2 + \overline{P^2}/2$. Taking $\kappa \to \infty$ in \eqref{fkappa} and using the dominated convergence theorem again we obtain 
\begin{gather*}
\left[ \sqrt{h^{3/2}\, f\, +\,  \beta} \right]_t\ +\ \left[ \lambda\, \sqrt{h^{3/2}\, f\, +\,  \beta} \right]_x\ 
\leqslant\  {\textstyle \frac{ \beta\, \lambda_x}{\sqrt{h^{3/2}\, f\, +\,  \beta}}}, \qquad \quad
f\ \eqdef\ \half \left( \overline{P^2}\, -\, \overline{P}^2 \right).
\end{gather*}
Taking now $\beta \to 0$ we obtain
\begin{equation}\label{Tg}
\left[ \sqrt{h^{3/2}\, f} \right]_t\ +\ \left[ \lambda\, \sqrt{h^{3/2}\, f} \right]_x\ 
\leqslant\ 0 \qquad \mathrm{in}\ (t_0,\infty) \times \mathds{R}.
\end{equation}
\textbf{Step 5.} Following \cite{wave3}, let $g \eqdef \sqrt{h^{3/2}\, f} \in L^\infty ((0,\infty), L^2(\mathds{R}))$. Let also $\varphi \in C^\infty_c(\mathds{R})$ satisfying $\varphi(x)=1$ for $|x| \leqslant 1$ and $\varphi(x)=0$ for $|x| \geqslant 2$. Then, for all $n \geqslant 1$, we have $g \varphi(x/n) \in L^\infty ((0,\infty), L^1(\mathds{R}))$. Then almost all $t>0$ are Lebesgue points of $t \mapsto \int_\mathds{R} g(t,x) \varphi(x/n) \mathrm{d}x,$ $\forall n \geqslant 1$.
Let $\bar{t}>0$ be a Lebesgue point of $t \mapsto \int_\mathds{R} g(t,x) \varphi(x/n) \mathrm{d}x$ and $\delta \in (0,\bar{t}/2)$. Let also $\psi \in C^\infty_c((0,\infty))$ satisfying
\begin{gather*}
\psi(t)\ =\ 0 \quad \mathrm{on}\quad (0,\delta/2) \cup (\bar{t}+\delta, \infty ), \qquad \psi(t)\ =\ 1 \quad \mathrm{on}\quad  (\delta, \bar{t}-\delta ),\\
0\ \leqslant\ \psi'(t)\ \leqslant\ C/\delta, \quad \mathrm{on}\quad (\delta/2,\delta), \qquad 
-\psi'(t)\ \geqslant\ C/\delta, \quad \mathrm{on}\quad (\bar{t}-\delta,\bar{t}+\delta).
\end{gather*}
Multiplying \eqref{Tg} by $\varphi(x/n) \psi(t)$, integrating on $(0,\infty) \times \mathds{R}$ and using integration by parts one obtains
\begin{gather*}
{\textstyle \frac{C}{\delta}}\, \int_{\bar{t}-\delta}^{\bar{t}+\delta} \int_\mathds{R} g(t,x)\, \varphi(x/n)\, \mathrm{d}x\, \mathrm{d}t\ 
\leqslant\ 
-\int_{\bar{t}-\delta}^{\bar{t}+\delta} \int_\mathds{R}  g(t,x)\, \varphi(x/n)\, \psi'(t)\, \mathrm{d}x\, \mathrm{d}t\\
\leqslant\ {\textstyle \frac{C}{\delta}}\, \int_{\delta/2}^{\delta} \int_\mathds{R} g(t,x)\, \varphi(x/n)\, \mathrm{d}x\, \mathrm{d}t\ +\
{\textstyle \frac{1}{n}}\, ||\lambda||_{L^\infty}\, \int_{\delta/2}^{\bar{t}+\delta} \int_\mathds{R} g(t,x)\, \left| \varphi'(x/n) \right|\, \mathrm{d}x\, \mathrm{d}t.
\end{gather*}
From \eqref{tauto0_2}, we have 
\begin{equation*}
\lim_{t \to 0} \int_\mathds{R} g(t,x)\, \varphi(x/n)\, \mathrm{d}x\ =\ 0 \quad \implies \quad 
\lim_{\delta \to 0} {\textstyle \frac{1}{\delta}}\, \int_{\delta/2}^{\delta}  \int_\mathds{R} g(t,x)\, \varphi(x/n)\, \mathrm{d}x\, \mathrm{d}t\ =\ 0.
\end{equation*}
Since $\bar{t}>0$ is a Lebesgue point of $t \mapsto \int_\mathds{R} g(t,x) \varphi(x/n) \mathrm{d}x$, then taking first $\delta \to 0$ and then $n \to \infty$ we obtain 
\begin{equation*}
g(\bar{t},x)\ =\ 0 \quad \mathrm{a.e.}\, (\bar{t},x) \in (0,\infty) \times \mathds{R}.
\end{equation*}
Hence $\overline{P^2} = \overline{P}^2$ almost everywhere, which implies that
$\nu^1_{t,x}(\xi) = \delta_{\overline{P}(t,x)}(\xi)$. The proof of $\nu^2_{t,x}(\zeta) = \delta_{\overline{Q}(t,x)}(\zeta)$ can be done similarly.
\qed

\section{The global weak solutions}\label{sec:gws}
We use in this section the precompactness results given in the previous section to prove that the limit $(h,u)$ given in Lemma \ref{Strong_conv} is a weak solution of \eqref{SGNxi}.

Let $(h^\varepsilon-\bar{h}, u^\varepsilon) $ be the solution given in Theorem \ref{thm:ex:ep}. Then, from Lemma \ref{lemYoungm}, Lemma \ref{lem:2Dirac}, Lemma \ref{lem:alpha+2}, \eqref{uh_x} and \eqref{uh_xi} we have that 
\begin{gather}\label{Wtconv_xi}
\left(P^\varepsilon, Q^\varepsilon, u_x^\varepsilon, h_x^\varepsilon \right)\qquad
\rightharpoonup \qquad \left(\overline{P}, \overline{Q}, u_x, h_x\right) \qquad
\mathrm{in}\ L^p_{loc}((0,\infty) \times \mathds{R}), \\ \nonumber
\left\| \left(P^\varepsilon\right)^2, \left(Q^\varepsilon\right)^2, \left(u_x^\varepsilon\right)^2, \left(h_x^\varepsilon\right)^2\right\|_{L^1(\Omega)}\quad
\to \quad \left\|\overline{P}^2, \overline{Q}^2, u_x^2, h_x^2\right\|_{L^1(\Omega)},
\end{gather}
for any $p \in [2,3)$ and compact set $\Omega \subset (0,\infty) \times \mathds{R} $. This implies  that 
\begin{align}\label{Stconv_xi}
\left(P^\varepsilon, Q^\varepsilon, u_x^\varepsilon, h_x^\varepsilon \right)\ 
&\to \left(\overline{P}, \overline{Q}, u_x, h_x\right) 
&\mathrm{in}\ L^2_{loc}((0,\infty) \times \mathds{R}).
\end{align}
Using Lemma \ref{lem:alpha+2} and Lemma \ref{Strong_conv} we obtain that for all $p \in [2,3)$, we have 
\begin{align}\label{Wtconv_tau}
\left(u_t^\varepsilon, h_t^\varepsilon \right)\ 
&\rightharpoonup \left(u_t, h_t \right)
&\mathrm{in}\ L^p_{loc}((0,\infty) \times \mathds{R}).
\end{align}
Now, using \eqref{Estimate9} and taking the weak limit $\varepsilon \to 0$ in \eqref{Appsys1} we obtain \eqref{SGN2a}.
Applying $\mathcal{L}_{h^\varepsilon}$ on \eqref{Appsys2} and multiplying by $ \varphi \in C^\infty_c((0,\infty) \times \mathds{R})$ we obtain
\begin{gather*}
\int_{\mathds{R}^+ \times \mathds{R}} \left\{  \left\{u^\varepsilon_t\, +\, u^\varepsilon\, u^\varepsilon_x\, +\,  3\, \gamma\, (h^\varepsilon)^{-2}\, h^\varepsilon_x \right\}\left\{h^\varepsilon\, \varphi\, -\,  (h^\varepsilon)^2\, h_x^\varepsilon\, \varphi_x\ -\ \third\, (h^\varepsilon)^3\, \varphi_{x x} \right\} +\, \half\, \varphi\, u^\varepsilon\, \mathscr{A}_x^\varepsilon \right\} \mathrm{d}x\, \mathrm{d}t\\
=\ \int_{\mathds{R}^+ \times \mathds{R}} \varphi_x \left\{ \mathscr{C}^\varepsilon\, +\, F(h^\varepsilon)\, -\, {\textstyle \frac{1}{2}}\, (h^\varepsilon)^2\, u^\varepsilon_x\, \mathscr{A}_x^\varepsilon\ +\ {\textstyle \frac{1}{48}}\, h^\varepsilon \left(\chi_\varepsilon(P^\varepsilon)\, +\, \chi_\varepsilon(Q^\varepsilon) \right) \right\}  \mathrm{d}x\, \mathrm{d}t .
\end{gather*}
From \eqref{Stconv_xi} and Lemma \ref{Strong_conv} we obtain the following convergence as $\varepsilon \to 0$ 
\begin{equation*}
\left\{h^\varepsilon\, \varphi\, -\,  (h^\varepsilon)^2\, h_x^\varepsilon\, \varphi_x\ -\ \third\, (h^\varepsilon)^3\, \varphi_{x x} \right\}\ 
\to\ 
\mathcal{L}_{h}\, \varphi
\quad \mathrm{in}\ L^2_{loc}((0,\infty) \times \mathds{R}).
\end{equation*}
Using again \eqref{Stconv_xi}, Lemma \ref{Strong_conv} and also \eqref{Wtconv_tau}  obtain the convergence 
\begin{equation*}
\left\{u^\varepsilon_t\, +\, u^\varepsilon\, u^\varepsilon_x\, +\,  3\, \gamma\, (h^\varepsilon)^{-2}\, h^\varepsilon_x \right\}\
\rightharpoonup\ 
\left\{u_t\, +\, u\, u_x\, +\,  3\, \gamma\, h^{-2}\, h_x \right\}
\quad \mathrm{in}\ L^2_{loc}((0,\infty) \times \mathds{R}).
\end{equation*}
We suppose that $\mathrm{supp}(\varphi) \subset [t_1,t_2] \times [a,b]$, then using the energy equation \eqref{energyconservationep} and Lemma \ref{lem:bounded2} we obtain
\begin{equation*}
\left| \int_{\mathds{R}^+ \times \mathds{R}}  \varphi_x\, (h^\varepsilon)^2\, u^\varepsilon_x\, \mathscr{A}_x^\varepsilon\, \mathrm{d}x\, \mathrm{d}t \right|\, \leqslant\  C\left\|\varphi_x\, (h^\varepsilon)^2\, u^\varepsilon_x \right\|_{L^\infty([t_1,t_2],L^2([a,b]))}  
\left\|\mathscr{A}_x^\varepsilon \right\|_{L^1([t_1,t_2],L^\infty ([a,b]))} \leqslant\ \varepsilon\, C.
\end{equation*}
Following the same argument we obtain 
\begin{equation*}
\left| \int_{\mathds{R}^+ \times \mathds{R}}  \varphi\,  u^\varepsilon\, \mathscr{A}_x^\varepsilon\, \mathrm{d}x\, \mathrm{d}t \right|\, \leqslant\ \varepsilon\, C. 
\end{equation*}
Then, taking $\varepsilon \to 0$, using Lemma \ref{lem:fg}, \eqref{Estimate7} and \eqref{Stconv_xi} we obtain \eqref{WS_xi}. 
Doing the proof of \eqref{tauto0} for any $t_0$, we obtain that $(h-\bar{h},u) \in C_r(\mathds{R}^+, H^1(\mathds{R})$. Lemma \ref{lem:Oleinik} implies \eqref{Ol_xi}. The inequality \eqref{alpha+2_xi} follows from Lemma \ref{lem:alpha+2}, \eqref{Wtconv_xi} and \eqref{Wtconv_tau}.  
Finally, the energy inequality \eqref{Eneeq_xi} follows from \eqref{Energyfin}.

\appendix
\section{Some classical lemmas}\label{App:A}

Here, we recall simple versions of some classical lemmas that are needed in this paper. 

We start this section by the following lemma on the Young measures.
\begin{lem}(\cite{Focusing})\label{lem:Young}
Let $\mathscr{O}$ be a subset of $\mathds{R}^n$ with a zero-measure boundary. For any bounded family $\{v^\varepsilon \}_\varepsilon \subset L^p(\mathscr{O},\mathds{R}^N)$ with $p>1$ there exists a subsequence denoted also $\{v^\varepsilon \}_\varepsilon$ and a family of probability measures on $\mathds{R}^N$, $\left\{ \mu_{y}, y \in \mathscr{O} \right\}$ such that for all $f \in C^0(\mathds{R}^N)$ with $f(\xi)=\smallO(|\xi|^p)$ at infinity and for all $\phi \in C^\infty_c(\mathscr{O})$ we have 
\begin{equation}
\lim_{\varepsilon \to 0}\, \int_\mathscr{O} \phi(y)\, f(v^\varepsilon(y))\, \mathrm{d}y\ =\ \int_\mathscr{O} \phi(y) \int_\mathds{R} f(\xi)\, \mathrm{d} \mu_{y}(\xi)\, \mathrm{d}y
\end{equation}
with 
\begin{equation}\label{Fatou0}
\int_\mathscr{O} \int_\mathds{R} |\xi|^p\, \mathrm{d} \mu_y(\xi)\, \mathrm{d}y\ \leqslant\ \liminf_{\varepsilon \to 0} \|u^\varepsilon\|_{L^p(\mathscr{O})}^p.
\end{equation}
\end{lem}
Also, some other results on strong and weak precompactness are needed, then we recall.
\begin{lem}(\cite{Evans})\label{lem:fg}
Let $\Omega$ be an open set of $\mathds{R}^n$, assuming that $f_n \to f$ in $L^p(\Omega)$ with $p \in (1,\infty)$, $g_n$ is bounded in $L^q$ with $q \in (1,\infty)$ and $g_n \rightharpoonup  g $ in $L^q(\Omega)$, then for any $\varphi \in L^r(\Omega)$ such that $1/p+1/q+1/r=1$, we have
\begin{equation}
\lim_{n \to  \infty}\int_\Omega f_n\, g_n\, \varphi\, \mathrm{d}x\ =\ \int_\Omega f\, g\, \varphi\, \mathrm{d}x.
\end{equation}
\end{lem}

\begin{lem}(\cite{Evans})\label{lem:Evans} For any $p>2$ we have $L^1_{loc}(\mathds{R}^2) \cap W^{-1,p}_{loc}(\mathds{R}^2) \Subset H^{-1}_{loc}(\mathds{R}^2)$. In other words, for any open, bounded, smooth set $U \subset \mathds{R}^2$, if the sequence $(f_n)_n$ is bounded in $L^1(U)\cap W^{-1,p}(U)$, then $(f_n)_n$ is relatively compact in 
$ H^{-1}(U)$.
\end{lem}
\begin{lem}(Lemma II.1 in \cite{DiPernaL1989})\label{lemRenor}
Let $c \in L^1(\mathds{R}^+,H^1_{loc}(\mathds{R})$ and $f \in L^\infty(\mathds{R}^+,L^2_{loc}(\mathds{R}))$. Let also $j_\varepsilon$ be a Friedrichs mollifier, then 
\begin{equation}
(c\, \partial_x f)\, \ast j_\varepsilon\ -\ c\, (\partial_x f \ast j_\varepsilon)\\\xrightarrow{\varepsilon \to 0} 0, \qquad \mathrm{in} \qquad L^1_{loc}(\mathds{R}^+ \times \mathds{R}).
\end{equation}
\end{lem}
\begin{lem}(Lemma C.1 in \cite{Lions})\label{lem:C_w}
Let $(f_n)_n$ be a bounded sequence in $L^\infty([0,T], L^2(\mathds{R}))$. If $f_n$ belongs to $C([0,T], H^{-1}(\mathds{R}))$ and  for any $\varphi \in H^1(\mathds{R})$, the map 
\begin{equation*}
t\ \mapsto\ \int_\mathds{R} \varphi(x)\, f_n(t,x)\, \mathrm{d}x
\end{equation*}
is uniformly continuous for $t \in [0,T]$ and $n \geqslant 1$, then $(f_n)_n$ is relatively compact in the space $C([0,T],L^2_w(\mathds{R}))$, where $L^2_w$ is the $L^2$ space equipped with its weak topology.
\end{lem}

\begin{lem}\textbf{(Compensated compactness \cite{PG1991})}\label{lem:PG}
Let $\Omega$ be an open set of $\mathds{R}^2$, let $a,b \in C(\Omega,\mathds{R})$ such that for all $(x_1,x_2) \in \Omega$ we have $a(x_1,x_2) \neq b(x_1,x_2)$. Let also $(f_n), (g_n)$ be bounded sequences in $L^2_{loc}(\Omega,\mathds{R})$ such that $f_n \rightharpoonup f$ and $g_n \rightharpoonup g$. If the sequences 
\begin{equation*}
\left\{ \partial_{x_1} f_n\ +\ \partial_{x_2} (a\, f_n)\right\}_n, \qquad \mathrm{and} \qquad 
\left\{ \partial_{x_1} g_n\ +\ \partial_{x_2} (b\, g_n)\right\}_n,
\end{equation*}
are relatively compact in $H^{-1}_{loc}(\Omega)$, then $f_n g_n \rightharpoonup f g $ in the sense of distributions. 
\end{lem}

Let $\Lambda$ be defined such that $\widehat{\Lambda f}=(1+\xi^2)^\frac{1}{2} \hat{f}$ and let $[A,B] \eqdef AB-BA$ be the commutator of the operators $A$ and $B$.
We recall now some estimates of the $H^s$ norm of the product, the commutator and the composition of functions.
\begin{lem}(\cite{kato1988commutator})
If $r\, \geqslant\, 0$, then $\exists C>0$ such that
\begin{align}
\|f\, g\|_{H^r}\ &\leqslant\ C\, \left( \|f\|_{L^\infty}\, \|g\|_{H^r}\ +\ \|f\|_{H^r}\, \|g\|_{L^\infty}\right), \label{Algebra} \\
\left\| \left[ \Lambda^r,\, f \right]\, g  \right\|_{L^2}\ &\leqslant\ C\, \left( \|f_x\|_{L^\infty}\, \|g\|_{H^{r-1}}\ +\ \|f\|_{H^r}\, \|g\|_{L^\infty} \right). \label{Commutator}
\end{align}
\end{lem}

\begin{lem}(\cite{constantin2002initial})
Let $F \in {C}^{\infty}(\mathds{R})$ with $F(0)=0$, then for any $m \in \mathds{N}$ there exists a continuous function $\tilde{F}$, such that for all $f \in H^m$ we have
\begin{equation}\label{Composition2}
\|F(f)\|_{H^m}\ \leqslant\ \tilde{F} \left( \|f\|_{L^{\infty}} \right)\, \|f\|_{H^m}.
\end{equation}
\end{lem}

\section{The energy equation}\label{App:B}
The goal of this section is to prove that smooth solutions of \eqref{SGNxi} (respectively \eqref{Appsys}) satisfy the energy equation \eqref{ene} (respectively \eqref{energyequationep}).
Taking $\varepsilon=0$, we notice that \eqref{ene} is a particular case of \eqref{energyequationep}.
We consider $\varepsilon \geqslant 0$ and $(h^\varepsilon,u^\varepsilon)$ smooth solutions of \eqref{Appsys}. Then we have
\begin{align*}
\half\, g \left[ \left( h^\varepsilon - \bar{h} \right)^2 \right]_t\, &=\, -\, g\, (h^\varepsilon - \bar{h}) \left(h^\varepsilon\, u^\varepsilon \right)_x\, +\, g\, (h^\varepsilon - \bar{h})\, \mathscr{A}_x^\varepsilon,\\
\half\, \gamma \left[ (h_x^\varepsilon)^2 \right]_t\, &=\, -\, \gamma\, h^\varepsilon_x \left(h^\varepsilon\, u^\varepsilon \right)_{xx}\, +\, g\, h_x^\varepsilon\, \mathscr{A}^\varepsilon\, -\, {\textstyle \frac{\sqrt{3\, \gamma}}{48\, (h^\varepsilon)^{1/2}}} h^\varepsilon_x \left( \chi_\varepsilon(P^\varepsilon) - \chi_\varepsilon(Q^\varepsilon) \right).
\end{align*}
Summing up we obtain 
\begin{gather}\nonumber
\half\, g \left[ \left( h^\varepsilon - \bar{h} \right)^2 \right]_t\, +\, \half\, \gamma \left[ (h_x^\varepsilon)^2 \right]_t\, =\,  g \left[ (h^\varepsilon - \bar{h})\, \mathscr{A}^\varepsilon \right]_x\\ \label{hene}
-\, g\, (h^\varepsilon - \bar{h}) \left(h^\varepsilon\, u^\varepsilon \right)_x\, -\, \gamma\, h^\varepsilon_x \left(h^\varepsilon\, u^\varepsilon \right)_{xx}\, -\, {\textstyle \frac{1}{96}} \left(Q^\varepsilon - P^\varepsilon \right)  \left( \chi_\varepsilon(P^\varepsilon) - \chi_\varepsilon(Q^\varepsilon) \right).
\end{gather}
Defining $\mathscr{X}^\varepsilon \eqdef \mathscr{C}^\varepsilon + F(h^\varepsilon)$, from \eqref{Appsys} we have
\begin{gather}\nonumber
\sixth\! \left[ (h^\varepsilon)^3\, (u^\varepsilon_x)^2 \right]_t\ 
=\ \half\,  (h^\varepsilon)^2\, (u^\varepsilon_x)^2\, \mathscr{A}^\varepsilon_x\, -\, \half\, (h^\varepsilon)^2\, (u_x^\varepsilon)^2 \left( h^\varepsilon\, u^\varepsilon \right)_x -\, \third\, (h^\varepsilon)^3\, u_x^\varepsilon \left( u^\varepsilon\, u_x^\varepsilon \right)_x\\ \label{uxene}
\, -\, \gamma\,  (h^\varepsilon)^3\, u_x^\varepsilon \left( (h^\varepsilon)^{-2}\, h_x^\varepsilon \right)_x -\, \third u^\varepsilon_x\, (h^\varepsilon)^3\, \partial_x\, \mathcal{L}^{-1}_{h^\varepsilon}\, \partial_x\, \mathscr{X}^\varepsilon\, +\, \third\, u^\varepsilon_x\, (h^\varepsilon)^3\, \partial_x\, \mathscr{B}^\varepsilon.
\end{gather}
Using again \eqref{Appsys} and the definition of $\mathcal{L}_{h^\varepsilon}$ we obtain
\begin{align}\nonumber
\half\! \left[ h^\varepsilon\, (u^\varepsilon)^2\right]_t\, 
=&\ \half\, (u^\varepsilon)^2\, \mathscr{A}^\varepsilon_x\, -\, \half\, (u^\varepsilon)^2 \left(h^\varepsilon\, u^\varepsilon\right)_x -\, h^\varepsilon\, (u^\varepsilon)^2\, u^\varepsilon_x\, -\, 3\, \gamma\, (h^\varepsilon)^{-1}\, u^\varepsilon\, h^\varepsilon_x\\ \nonumber
&\, -\, h^\varepsilon\, u^\varepsilon\, \mathcal{L}_{h^\varepsilon}^{-1} \mathscr{X}^\varepsilon_x\, +\, h^\varepsilon\, u^\varepsilon\, \mathscr{B}^\varepsilon\\ \nonumber
=&\, -\, \half\, (u^\varepsilon)^2 \left(h^\varepsilon\, u^\varepsilon\right)_x -\, h^\varepsilon\, (u^\varepsilon)^2\, u^\varepsilon_x\, -\, 3\, \gamma\, (h^\varepsilon)^{-1}\, u^\varepsilon\, h^\varepsilon_x\, -\, u^\varepsilon\, \mathscr{X}^\varepsilon_x\\ \nonumber
&\, -\, \third\, u^\varepsilon\, \partial_x\, (h^\varepsilon)^3\, \partial_x\, \mathcal{L}^{-1}_{h^\varepsilon}\, \partial_x\, \mathscr{X}^\varepsilon\, +\, \third\, u^\varepsilon\, \partial_x\, (h^\varepsilon)^3\, \partial_x\, \mathscr{B}^\varepsilon\, +\, \half\, u^\varepsilon \left((h^\varepsilon)^2\, u^\varepsilon_x\, \mathscr{A}^\varepsilon_x\right)_x\\ \nonumber
&\, - {\textstyle \frac{1}{48}} u^\varepsilon \left[h^\varepsilon \left(\chi_\varepsilon (P^\varepsilon)\, +\, \chi_\varepsilon(Q^\varepsilon) \right) \right]_x\\ \nonumber
=&\, -\, \half \left[ (u^\varepsilon)^3 h^\varepsilon \right]_x  -\, 3\, \gamma\, (h^\varepsilon)^{-1}\, u^\varepsilon\, h^\varepsilon_x\, - \left[ u^\varepsilon\, \mathscr{X}^\varepsilon\right ]_x\, -\, \third \left[ u^\varepsilon\, (h^\varepsilon)^3\, \partial_x\, \mathcal{L}^{-1}_{h^\varepsilon}\, \partial_x\, \mathscr{X}^\varepsilon \right]_x\\ \nonumber
&\, +\, \third \left[ u^\varepsilon\, (h^\varepsilon)^3\, \partial_x\, \mathscr{B}^\varepsilon\right]_x\, +\, \half \left[ u^\varepsilon\, (h^\varepsilon)^2\, u^\varepsilon_x\, \mathscr{A}^\varepsilon_x\right]_x\, -\, {\textstyle \frac{1}{48}}  \left[u^\varepsilon\, h^\varepsilon \left(\chi_\varepsilon (P^\varepsilon)\, +\, \chi_\varepsilon(Q^\varepsilon) \right) \right]_x\\ \nonumber
&\, +\, u^\varepsilon_x\, \mathscr{X}^\varepsilon\, +\, \third u^\varepsilon_x\, (h^\varepsilon)^3\, \partial_x\, \mathcal{L}^{-1}_{h^\varepsilon}\, \partial_x\, \mathscr{X}^\varepsilon\, -\, \third\, u^\varepsilon_x\, (h^\varepsilon)^3\, \partial_x\, \mathscr{B}^\varepsilon\, -\, \half\,  (h^\varepsilon)^2\, (u^\varepsilon_x)^2\, \mathscr{A}^\varepsilon_x\\ \label{uene}
&\, +\, {\textstyle \frac{1}{96}} \left(P^\varepsilon\, +\, Q^\varepsilon \right) \left(\chi_\varepsilon (P^\varepsilon)\, +\, \chi_\varepsilon(Q^\varepsilon) \right).
\end{align}
Summing up \eqref{uxene} and \eqref{uene} one obtains 
\begin{gather}\nonumber
\half\! \left[ h^\varepsilon\, (u^\varepsilon)^2\right]_t\, +\, \sixth\! \left[ (h^\varepsilon)^3\, (u^\varepsilon_x)^2 \right]_t\ =\ 
  -\, 3\, \gamma\, (h^\varepsilon)^{-1}\, u^\varepsilon\, h^\varepsilon_x\, +\, u^\varepsilon_x\, \mathscr{X}^\varepsilon\, -\, \half\, (h^\varepsilon)^2\, (u_x^\varepsilon)^2 \left( h^\varepsilon\, u^\varepsilon \right)_x\\ \nonumber
 -\, \third\, (h^\varepsilon)^3\, u_x^\varepsilon \left( u^\varepsilon\, u_x^\varepsilon \right)_x -\, \half \left[ (u^\varepsilon)^3 h^\varepsilon \right]_x - \left[ u^\varepsilon\, \mathscr{X}^\varepsilon\right ]_x\, -\, \third \left[ u^\varepsilon\, (h^\varepsilon)^3\, \partial_x\, \mathcal{L}^{-1}_{h^\varepsilon}\, \partial_x\, \mathscr{X}^\varepsilon \right]_x\\ \nonumber
+\, \third \left[ u^\varepsilon\, (h^\varepsilon)^3\, \partial_x\, \mathscr{B}^\varepsilon\right]_x\, 
 +\, \half \left[ u^\varepsilon\, (h^\varepsilon)^2\, u^\varepsilon_x\, \mathscr{A}^\varepsilon_x\right]_x -\, {\textstyle \frac{1}{48}}  \left[u^\varepsilon\, h^\varepsilon \left(\chi_\varepsilon (P^\varepsilon)\, +\, \chi_\varepsilon(Q^\varepsilon) \right) \right]_x\\ \label{uene2}
-\, \gamma\,  (h^\varepsilon)^3\, u_x^\varepsilon \left( (h^\varepsilon)^{-2}\, h_x^\varepsilon \right)_x\,
   +\, {\textstyle \frac{1}{96}} \left(P^\varepsilon\, +\, Q^\varepsilon \right) \left(\chi_\varepsilon (P^\varepsilon)\, +\, \chi_\varepsilon(Q^\varepsilon) \right).
\end{gather}
Using \eqref{Psi} and \eqref{Bdef} we obtain 
\begin{align*}
\third\, (h^\varepsilon)^3\, \partial_x\, \mathscr{B}^\varepsilon \, 
=&\, \third\, (h^\varepsilon)^3\, \partial_x\, \mathcal{L}^{-1}_{h^\varepsilon} \left\{-\half\, u^\varepsilon\, \mathscr{A}_x^\varepsilon \right\} -\, \half\, (h^\varepsilon)^2\, u_x^\varepsilon\, \mathscr{A}_x^\varepsilon\, +\, {\textstyle \frac{1}{48}}\, h^\varepsilon \left(\chi_\varepsilon (P^\varepsilon)\, +\, \chi_\varepsilon(Q^\varepsilon) \right)\\
&\, +\, (h^\varepsilon)^3\, \partial_x\, \mathcal{L}_{h^\varepsilon}^{-1} 
\left\{ h^\varepsilon \int_{-\infty}^x \left( \half\, (h^\varepsilon)^{-1}\, u_x^\varepsilon\, \mathscr{A}_x^\varepsilon\, -\, {\textstyle \frac{1}{48\, (h^\varepsilon)^2}} \left(\chi_\varepsilon (P^\varepsilon)\, +\, \chi_\varepsilon(Q^\varepsilon) \right) \right) \mathrm{d}y \right\}\\
=&\, \third\, (h^\varepsilon)^2\, \mathscr{V}^\varepsilon_1\, -\, \half\, (h^\varepsilon)^2\, u_x^\varepsilon\, \mathscr{A}_x^\varepsilon\, +\, {\textstyle \frac{1}{48}}\, h^\varepsilon \left(\chi_\varepsilon (P^\varepsilon)\, +\, \chi_\varepsilon(Q^\varepsilon) \right).
\end{align*}
Using now \eqref{CR} and \eqref{uene2} we obtain 
\begin{gather}\nonumber
\half\! \left[ h^\varepsilon\, (u^\varepsilon)^2\right]_t\, +\, \sixth\! \left[ (h^\varepsilon)^3\, (u^\varepsilon_x)^2 \right]_t\ =\ 
-\, 3\, \gamma\, (h^\varepsilon)^{-1}\, u^\varepsilon\, h^\varepsilon_x\, +\, u^\varepsilon_x\, \mathscr{X}^\varepsilon\, -\, \half\, (h^\varepsilon)^2\, (u_x^\varepsilon)^2 \left( h^\varepsilon\, u^\varepsilon \right)_x\\ \nonumber
 -\, \third\, (h^\varepsilon)^3\, u_x^\varepsilon \left( u^\varepsilon\, u_x^\varepsilon \right)_x -\, \half \left[ (u^\varepsilon)^3 h^\varepsilon \right]_x - \left[ u^\varepsilon\, \mathscr{R}^\varepsilon\right ]_x\,  -\, \left[u^\varepsilon \left(\half\, g\, \left( (h^\varepsilon)^2 - \bar{h}^2 \right) -\, 3\, \gamma \ln \left(h^\varepsilon/\bar{h} \right) \right) \right]_x\\ \label{uene3}
+\, \third \left[ u^\varepsilon\, (h^\varepsilon)^2\, \mathscr{V}_1^\varepsilon\right]_x -\, \gamma\,  (h^\varepsilon)^3\, u_x^\varepsilon \left( (h^\varepsilon)^{-2}\, h_x^\varepsilon \right)_x
 +\, {\textstyle \frac{1}{96}} \left(P^\varepsilon\, +\, Q^\varepsilon \right) \left(\chi_\varepsilon (P^\varepsilon)\, +\, \chi_\varepsilon(Q^\varepsilon) \right).
\end{gather}
Forward calculations lead to 
\begin{gather}\nonumber
g\, (h^\varepsilon - \bar{h}) \left(h^\varepsilon\, u^\varepsilon \right)_x\, +\, \gamma\, h^\varepsilon_x \left(h^\varepsilon\, u^\varepsilon \right)_{xx}\, +\, 3\, \gamma\, (h^\varepsilon)^{-1}\, u^\varepsilon\, h^\varepsilon_x\, -\, u^\varepsilon_x\, \mathscr{X}^\varepsilon\, +\, \half\, (h^\varepsilon)^2\, (u_x^\varepsilon)^2 \left( h^\varepsilon\, u^\varepsilon \right)_x\\ \nonumber
+\, \third\, (h^\varepsilon)^3\, u_x^\varepsilon \left( u^\varepsilon\, u_x^\varepsilon \right)_x\, +\, \gamma\,  (h^\varepsilon)^3\, u_x^\varepsilon \left( (h^\varepsilon)^{-2}\, h_x^\varepsilon \right)_x\
=\ \half\, g\, \left[u^\varepsilon \left(h^\varepsilon - \bar{h} \right)^2 \right]_x\\ \label{FCLT}
+\, \sixth\, \left[ (h^\varepsilon)^3\, u^\varepsilon\, (u^\varepsilon_x)^2 \right]_x +\, 3\, \gamma \left[u^\varepsilon\, \ln \left( h^\varepsilon/\bar{h}\right) \right]_x + \half\, \gamma \left[ u^\varepsilon\, (h_x^\varepsilon)^2 \right]_x +\, \gamma \left[ h^\varepsilon\, h^\varepsilon_x\, u^\varepsilon_x \right]_x.
\end{gather}
Summing up \eqref{hene}, \eqref{uene3} and \eqref{FCLT} we obtain
\begin{gather*}
\left[\half\, h^\varepsilon\, (u^\varepsilon)^2\, +\ \half\, g\, \left( h^\varepsilon - \bar{h} \right)^2\, +\, \sixth\, (h^\varepsilon)^3\, (u^\varepsilon_x)^2\,  +\, \half\, \gamma\, (h_x^\varepsilon)^2 \right]_t\ + \\
\Big[\half (u^\varepsilon)^3 h^\varepsilon\, -\, g\, (h^\varepsilon - \bar{h})\, \mathscr{A}^\varepsilon\,   +\,  u^\varepsilon\, \mathscr{R}^\varepsilon\, +\, \half\, g\, u^\varepsilon  \left( (h^\varepsilon)^2 - \bar{h}^2 \right) 
 -\, \third\, u^\varepsilon\, (h^\varepsilon)^2\, \mathscr{V}_1^\varepsilon\\
+\, \half\, g\, u^\varepsilon \left(h^\varepsilon - \bar{h} \right)^2\, +\, \sixth\, (h^\varepsilon)^3\, u^\varepsilon\, (u^\varepsilon_x)^2\, +\, \half\, \gamma\, u^\varepsilon\, (h_x^\varepsilon)^2\, +\, \gamma\, h^\varepsilon\, h^\varepsilon_x\, u^\varepsilon_x
\Big]_x \\
=\ {\textstyle \frac{1}{48}}\, P^\varepsilon \chi_\varepsilon(P^\varepsilon)\ +\ {\textstyle \frac{1}{48}}\, Q^\varepsilon \chi_\varepsilon(Q^\varepsilon).
\end{gather*}
This is \eqref{energyequationep}.

%

\end{document}